\documentclass[12pt,a4paper]{article}
\usepackage{amssymb,amsmath,amsthm,verbatim,url}
\usepackage{graphicx,color}
\newtheorem{theorem}{\sc Theorem.}[section]
\newtheorem{lemma}[theorem]{\sc Lemma.}
\newtheorem{remark}[theorem]{\sc Remark.}

\newenvironment{AMS}%
{{\upshape\bfseries AMS subject classifications. }\ignorespaces}{}
\newenvironment{keywords}{{\upshape\bfseries Key words. }\ignorespaces}{}

\newcommand{\bRplus}{{\mathbb R}_{>0}}
\newcommand{\bRgeq}{{\mathbb R}_{\geq 0}}

\newcommand{\bR}{{\mathbb R}}
\newcommand{\bN}{{\mathbb N}}

\newcommand{\bS}{{\mathbb S}}

\DeclareMathOperator*{\argmin}{arg\,min}

\newcommand{\dH}[1]{\;{\rm d}{\mathcal{H}}^{#1}} 

\newcommand{\Id}{{\rm Id}}

\newcommand{\PM}{P_{\mathcal M}}
\newcommand{\HM}{H_{\mathcal M}}
\newcommand{\Mnu}{{\rm n}}
\newcommand{\Mmu}{{\rm m}}
\newcommand{\cPsi}{c_\Psi}
\newcommand{\uD}{w_D}
\newcommand{\nabs}{\nabla_{\!s}}

\newcommand{\tanspace}{{\rm T}}        

\renewcommand{\vec}{}
\def\epsilon{\varepsilon}

\def\hat{\widehat}

\begin{document}
\title{
A finite element method for anisotropic crystal growth on surfaces 
}

\author{Harald Garcke\footnotemark[2]\ \and 
        Robert N\"urnberg\footnotemark[3]}

\renewcommand{\thefootnote}{\fnsymbol{footnote}}
\footnotetext[2]{Fakult{\"a}t f{\"u}r Mathematik, Universit{\"a}t Regensburg, 
93040 Regensburg, Germany \\ {\tt harald.garcke@ur.de}}
\footnotetext[3]{Dipartimento di Mathematica, Universit\`a di Trento,
38123 Trento, Italy \\ {\tt robert.nurnberg@unitn.it}}

\date{}

\maketitle

\begin{abstract}
Phase transition problems  on curved surfaces can lead to a panopticon of fascinating patterns.
In this paper we consider finite element approximations of phase field models with a
spatially inhomogeneous and anisotropic surface energy density. 
The problems are either posed in $\bR^3$ or on a two-dimensional
hypersurface in $\bR^3$. In the latter case, a fundamental choice regarding
the anisotropic energy density has to be made.
One possibility is to use a density defined in the ambient space
$\bR^3$. However, we propose and advocate for an alternative,
where a density is defined on a fixed chosen tangent space, and is then 
moved along geodesics to the other tangent spaces.
Our numerical method can be employed in all of the above situations, 
where for the problems on hypersurfaces the algorithm uses 
parametric finite elements.
We prove an unconditional stability result for our schemes and present
several numerical experiments, including for the modelling of ice 
crystal growth on a sphere.
\end{abstract} 



\renewcommand{\thefootnote}{\arabic{footnote}}

\setcounter{equation}{0}
\section{Introduction} 
\begin{figure}
\center
\includegraphics[angle=0,width=0.5\textwidth]{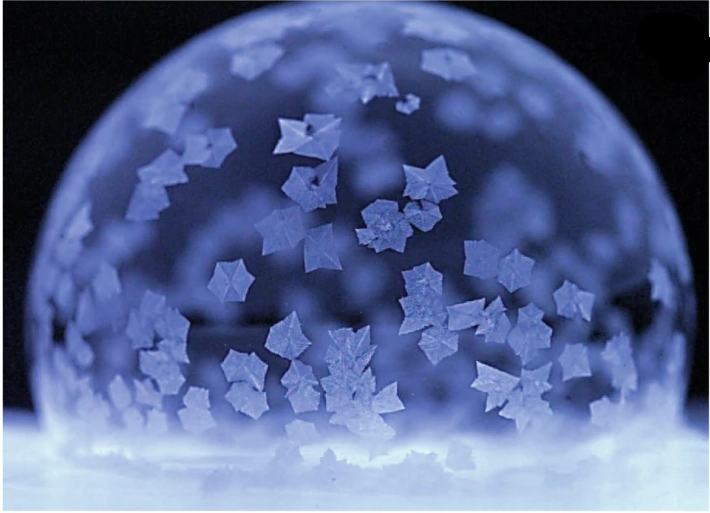}
\caption{Freezing of a soap bubble deposited on an ice disk.
Taken from \cite[Fig.\ 2]{AhmadiNKYB19}. This figure is used under a Creative Commons Attribution 4.0 license (CC BY 4.0). No changes have been made to the original photograph.}
\label{fig:soapbubble}
\end{figure}

Crystal growth on curved surfaces can lead to a multitude of interesting patterns. This phenomenon is one example of a phase change problem on a surface. 
Other applications involve phase separation on surfaces, the formation of two phases in vesicles or in lipid raft formation, see \cite{AhmadiNKYB19, nsns2phase,ElliottS10a,GarckeKRR16,GeraS17,LeeYPKLJKKK22,OlshanskiiPQ23,OrtelladoG20,Ratz16,SaxenaL99,YoonPWLK20}.
In this paper we numerically approximate interface evolutions on manifolds governed by an inhomogeneous, anisotropic interfacial energy  
by using a phase field approach. A particular application we have in mind is dendritic ice crystal growth on surfaces like soap bubbles, see, e.g., the fascinating pictures in \cite{AhmadiNKYB19}, partly reproduced in Figure~\ref{fig:soapbubble}.
Other possible applications are dendritic growths 
on aircraft bodies or metal shaped bodies, see \cite{ThomasCM96}, and phase separation on surfaces, see \cite{SaxenaL99,RatzV06, DuJT11,Ratz16,GeraS17}. 
In the latter case one solves a Cahn--Hilliard equation on a surface with either an isotropic or an anisotropic surface energy.

Although phase field models in the Euclidean space have received a lot of attention, see
\cite{Kobayashi93,GarckeNS04,DeckelnickDE05,TorabiLVW09,Steinbach09,SteinbachS23,DuF19handbook,BanschDGP23}, 
 not much is known for (anisotropic) phase field approaches for interface evolution problems on surfaces. 
Similarly, while for anisotropic phase field models in the Euclidean space 
a lot is known for the analysis and numerical analysis, see
\cite{ElliottS96,TorabiLVW09,GraserKS13,eck,pfsi,gknz,GarckeKW23}
and the references therein, not much is known for anisotropic models 
on surfaces. 

There have been some numerical computations for phase field models describing crystal growth on surfaces, see 
\cite{OrtelladoG20,YoonPWLK20,LeeYPKLJKKK22}. However, no numerical analysis has been performed so far. In what follows, we 
will first introduce the governing equations leading to an anisotropic phase field  model on a surface. This system reduces to an anisotropic Cahn--Hilliard equation in situations in which some terms involving time derivatives are neglected. On a surface the interfacial energy is defined on the tangent spaces of the surface.
Here, a reasonable choice on how to choose the anisotropy when the tangent space changes has to be taken, in order to model physically realistic situations. We basically consider two cases. In the first case we fix an anisotropy in $\bR^3$ and then restrict the anisotropy to the respective tangent spaces.  
The advantage is that the anisotropic density need not depend on
space, and existing physical models and numerical methods can be easily
extended from flat domains to surfaces. However, it
will turn out that this choice has certain undesirable properties. A second choice is obtained by moving an anisotropy given on one fixed tangent space along geodesics to the other tangent spaces. This allows, for example, to choose a six-fold anisotropy on all tangent spaces, see Section~\ref{sec:cons2d}.

For the numerical analysis, we generalize anisotropies introduced by
Barrett, Garcke and N\"urnberg in \cite{triplejANI,eck,vch} to the surface case. 
This enables us to show stability bounds as well as to prove existence and uniqueness results for fully discrete approximations. These results are shown both for the case of a smooth potential, as well as for the case of an obstacle 
potential. We remark that the anisotropies can also depend on space, i.e., they 
can be inhomogeneous. We will then demonstrate the effect of an inhomogeneous energy with the help of numerical simulations. In addition, we will show the effect of different choices of the anisotropy. In particular, in the case where a 3d-anisotropy 
is restricted to the tangent spaces, the form of the anisotropy can change heavily from tangent space 
to tangent space. We observe a change from a six-fold anisotropy to a four-fold anisotropy.
Moreover, a convergence experiment is presented using an explicit solution constructed by R\"atz in \cite{Ratz16}. Finally, computations for spinodal decomposition and crystal growth are also presented, where the latter leads to snow crystal growth on manifolds.

Let us now discuss literature related to this work.
Basic information on 
parametric methods for curvature flow and its anisotropic variants  can be found in \cite{DeckelnickDE05,Pozzi07,DziukE07a,curves3d,hypbol,bgnreview,finsler,eqdproc,BanschDGP23}. Related are also the works
\cite{MikulaS04,MikulaS06} for 
curvature flows on graphs. In the computer graphics literature also 
anisotropies on manifolds have been used and we refer to 
\cite{SeongJC08,SeongJC09,ZhuangZCJ14} for details.
Phase separation on manifolds using the Cahn--Hilliard model on surfaces has been studied numerically in \cite{SaxenaL99,DuJT11,Ratz16,GeraS17,CaetanoE21,OlshanskiiPQ23}. 
Analytic results for inhomogeneous anisotropies can be found in 
\cite{BellettiniP96,AlfaroGHMS10,finsler}. 

The outline of this paper is as follows. In Section~\ref{sec2} we present the governing phase field equations 
on a surface in its strong and weak formulation. In Section~\ref{sec:anen} we present the different choices of the anisotropy and state and prove certain qualitative properties. Section~\ref{sec:fea} is devoted to a fully discrete finite element approximation of the anisotropic phase field equations on surfaces. We also state existence, uniqueness and stability results. 
Their proofs are mostly straightforward extensions of the results in
\cite{vch} to surfaces and to spatially dependent anisotropies.
Finally in Section~\ref{sec:nr} we present several numerical computations
which demonstrate convergence as well as several qualitative properties of solutions.

\section{The mathematical model and its weak formulation}\label{sec2}

\subsection{The underlying anisotropic interfacial energy}
Let $\mathcal M \subset \bR^3$ be a given stationary smooth manifold, with or without
boundary. We let $\Mnu$ denote a continuous unit normal field on 
$\mathcal M$, and $\Mmu$ the outer conormal on $\partial \mathcal M$.
Let $(\Gamma(t))_{t\in[0,T]}$ be a family of evolving curves on
$\mathcal M$. Then we consider a general anisotropic energy of the form
\begin{equation*} 
\mathcal{E}(\Gamma) = \int_{\Gamma} \gamma(z, \nu_{\mathcal M}(z)) \dH1(z)
= \int_{\Gamma} \gamma(\cdot, \nu_{\mathcal M}) \dH1,
\end{equation*}
where $\gamma : \mathcal M \times \bR^3 \to \bRgeq$ is a given anisotropic
energy density that is spatially inhomogeneous and absolutely
one-homogeneous in the second argument. 
In addition, $\nu_{\mathcal M}$ denotes a geodesic normal field
to $\Gamma$ on $\mathcal M$, i.e., it lies in the normal space of $\Gamma$ 
and in the tangent space of $\mathcal M$. We always denote by $\dH{d}$, $d\in \bN$, integration with respect to the $d$-dimensional Hausdorff measure. For more information about the geometry
of curves and surfaces we refer to \cite{bgnreview,curves3d}.

For what follows, we assume that $\Gamma(t)$ separates the manifold into two
regions: $\mathcal M_\pm(t)$ with 
$\mathcal M = \overline{\mathcal M_+(t)} \cup \overline{\mathcal M_-(t)}$
and $\Gamma(t) = \overline{\mathcal M_+(t)} \cap \overline{\mathcal M_-(t)}$.
{From} now on we assume that $\nu_{\mathcal M}$ points into $\mathcal M_+(t)$.
See, e.g., \cite{curves3d}.

The sharp interface problem for anisotropic crystal growth on a manifold is the
surface Stefan problem with surface tension and kinetic undercooling.
To define the flow we introduce $\kappa_\gamma$, 
the anisotropic geodesic curvature of $\Gamma(t)$, see below for a 
precise definition,
and the velocity $\mathcal V$ of $\Gamma(t)$ in the direction of the normal 
$\nu_{\mathcal M}$. We then seek $w$ defined on $\mathcal M$,
which depending on the setting can be either a temperature field
or a concentration field, such that
\begin{subequations} \label{eq:SI}
\begin{alignat}{2}
\vartheta w_t - 
\nabs \cdot (\mathcal{K}\nabs w)&=0 \quad&&\text{in}\quad {\mathcal M_+(t)} \cup {\mathcal M_-(t)},\\
a	w&= \alpha \kappa_\gamma - \rho \mathcal V/\beta(\cdot,\nu_{\mathcal M})\quad&&\text{on}\quad \Gamma(t) ,\\
\lambda \mathcal V& = - [\mathcal{K} \nabs w \cdot \nu_{\mathcal M}] \quad&&\text{on}\quad \Gamma(t),
\end{alignat}
\end{subequations}
where $[.] $ denotes the jump of a quantity across the interface, 
$\beta$ is a kinetic mobility and 
$\vartheta$, $\mathcal{K}$, $a$, $\alpha$, $\rho$, $\lambda$ are
physical parameters, which for simplicity are assumed to be constant. 
The system then needs to be closed with boundary conditions for $w$ on
$\partial\mathcal{M}$, as well as initial conditions for $\Gamma(0)$ 
and possibly $w(0)$.
We note that in the case $\vartheta=0$ we obtain the
surface Mullins--Sekerka problem.
In addition, $\nabs$ denotes the surface gradient on $\mathcal M$, and we
define similarly the surface divergence and the surface Laplacian,
$\Delta_s = \nabs \cdot \nabs$, see, e.g., \cite{bgnreview}. 
Finally, the anisotropic geodesic curvature is defined 
as the first variation of the energy $\mathcal{E}$, so that
$\kappa_\gamma = - \nabs\cdot(\PM \gamma_p(\cdot,\nu_{\mathcal M}))$,
see, e.g., \cite{curves3d}, where 
$\gamma_p = (\gamma_{p_1}, \gamma_{p_2}, \gamma_{p_3})^T$ denotes the first
derivatives of $\gamma$ with respect to the second argument, and where
\begin{equation*} 
\PM = \Id - \Mnu \otimes \Mnu
\end{equation*}
is the projection on the tangent space of $\mathcal M$.
Similarly to \cite{pfsi}, we now introduce a corresponding phase field model,
where in a first step we replace $\mathcal{E}(\Gamma)$ by an analogue
for a diffuse interface.

Following \cite{McFaddenWBCS93,WheelerM96,CahnT94,ElliottS97,Elliott97}, 
we define 
\begin{equation*} 
A(z, p) = \tfrac12 \gamma^2(z, p), \qquad z \in \mathcal M, p \in \bR^3.
\end{equation*} 
On introducing a phase field parameter $\varphi : \mathcal M \to \bR$,
where later on $\mathcal M_\pm(t) \approx \{ z \in \mathcal M : 
\pm \varphi(z, t) > 0 \}$, we consider the anisotropic Ginzburg--Landau energy
\begin{equation} \label{eq:GL}
\mathcal{E}_\epsilon(\varphi) =
\int_{\mathcal M} \epsilon A(\cdot,\nabs \varphi) 
+ \epsilon^{-1} \Psi(\varphi) \dH2,
\end{equation}
where $\epsilon > 0$ is an interfacial parameter and
$\Psi : \bR \to [0, \infty]$ is a suitable potential function. 
In this work we consider either the smooth double-well potential
\begin{equation}
\Psi(s)=\tfrac{1}{4}(1-s^2)^2, \quad s \in \bR, \label{eq:quartic}
\end{equation}
or the double-obstacle potential
\begin{equation} \label{eq:obs}
\Psi(s) = \Psi_0(s) + \mathbb{I}_{[-1,1]}(s), \quad s \in \bR,
\end{equation}
where $\Psi_0(s)=\frac{1}{2}(1-s^2)$ and $\mathbb{I}_{[-1,1]}$ stands for the 
indicator function of the interval ${[-1,1]}$, see \cite{BloweyE91}, i.e.,  $\mathbb{I}_{[-1,1]}$ is zero on ${[-1,1]}$ and $\infty$ outside of the interval ${[-1,1]}$. Note that \eqref{eq:GL} is the natural
generalization of well-known phase field free energies in the flat case to the
case of a smooth manifold studied here. We refer to 
\cite{Elliott97,ElliottS96,DeckelnickDE05,eck,vch} for more details in the flat
case.

In the following presentation, for ease of exposition, we assume that $\Psi$ 
is smooth. 
Rigorous arguments involving the obstacle potential \eqref{eq:obs} involve
subdifferentials and variational inequalities. We leave these details to the
reader, see also \cite{eck,vch}, and make them more precise when we state 
our numerical approximations.

\subsection{Strong formulation} 

Given the Ginzburg--Landau energy \eqref{eq:GL}, for a given
$\epsilon > 0$ and $\Psi$, 
in this paper we want to study the general anisotropic 
surface phase field equations
\begin{subequations} \label{eq:general}
\begin{align} \label{eq:generala}
\vartheta w_t + \tfrac12\lambda\varphi_t - \nabs \cdot (\mathcal{K}\nabs w) &
= 0, \\
\tfrac12\cPsi a w & = 
\rho\epsilon\mu(\cdot, \nabs\varphi)\varphi_t
-\alpha\epsilon\nabs\cdot (\PM A_p(\cdot,\nabs\varphi))
+ \alpha\epsilon^{-1}\Psi'(\varphi)  \label{eq:generalb}
\end{align}
\end{subequations}
on $\mathcal{M}$, together with suitable boundary and initial conditions for 
$w$ and $\varphi$. 
Here
\begin{align*} 
\nabs \cdot (\PM A_p(\cdot,\nabs\varphi) )&
= \nabs \cdot A_p(\cdot,\nabs\varphi)
- \nabs \cdot ((A_p(\cdot,\nabs\varphi) \cdot \Mnu) \Mnu) \nonumber \\ &
=  \nabs \cdot A_p(\cdot,\nabs\varphi) - A_p(\cdot,\nabs\varphi) \cdot \Mnu
\nabs \cdot \Mnu \nonumber \\ &
= \nabs\cdot A_p(\cdot,\nabs\varphi) + \HM A_p(\cdot,\nabs\varphi) \cdot \Mnu,
\end{align*}
with
\begin{equation*} 
\HM = - \nabs \cdot \Mnu
\end{equation*}
denoting the mean curvature of $\mathcal M$.

In addition, we define 
\[
\cPsi = \int_{-1}^1 \sqrt{2\Psi(s)}\;{\rm d}s,
\]
which is needed to relate the phase field approach to the sharp interface 
limit \eqref{eq:SI}, see \cite{pfsi,Ratz16}.
Moreover, $\mu : \mathcal{M} \times \bR^3 \to \bRplus$ is given 
by $\gamma/\beta$ to ensure that \eqref{eq:SI} is recovered in the
sharp interface limit, see \cite{pfsi}. 

We note that the general model \eqref{eq:general} includes several special
cases, amongst which are the 
anisotropic viscous Cahn--Hilliard equation ($\vartheta = 0$),
the anisotropic Cahn--Hilliard equation ($\vartheta = \rho = 0$)
and the anisotropic Allen--Cahn equation ($a = 0$).

We remark that in applications $\gamma$ is often chosen spatially homogeneous,
so that
\begin{equation} \label{eq:gamma0}
\gamma(z, p) = \gamma_0(p) \quad \forall z \in \mathcal{M}, p \in \bR^3.
\end{equation}

Using the approach of \cite{BellettiniP96, BellettiniPV96} it can be shown that 
for the sharp interface limit,
$\epsilon\to0$, in the sense of $\Gamma$-limits it holds that 
\[
\frac1{\cPsi} \mathcal{E}_\epsilon
\to
\mathcal{E}_0,
\]
where for a 
subset $E$ of $\mathcal M$ with finite perimeter one defines
\[ \mathcal{E}_0 (\chi_E ) =\int_{\partial E} \gamma(\cdot, \nu_{\partial E}) \dH1,
\]
where $\nu_{\partial E}$ is the outward geodesic normal to
$E \subset \mathcal M$. In a similar way, the sharp interface limit of the
evolution equations \eqref{eq:general}
can also be investigated. 
In the planar case this sharp interface limit has been studied in 
\cite{Caginalp86,CaginalpC98,McFaddenWBCS93,BellettiniP96,GarckeSN99,AlfaroGHMS10},
and using ideas from \cite{ElliottS10a,GarckeKRR16} it is possible to translate
these ideas to surfaces.

\subsection{Weak formulation}

The natural weak formulation of \eqref{eq:general}, with the associated
boundary conditions
\[
A_p(\cdot, \nabs\varphi) \cdot \Mmu = 0 \quad \mbox{on}\ \partial \mathcal{M},
\qquad \nabs w \cdot \Mmu  = 0 \quad \mbox{on}\ \partial_N\mathcal{M},
\qquad w = \uD \quad \mbox{on}\ \partial_D\mathcal{M},
\]
where $\partial\mathcal{M} = \overline{\partial_N\mathcal{M} }\cup \overline{\partial_D\mathcal{M}}$ with relatively open subsets $\partial_N\mathcal{M}$ and $\partial_D\mathcal{M}$ of $\partial\mathcal{M}$ such that  $\partial_N\mathcal{M} \cap \partial_D\mathcal{M} =\emptyset$,
 and where $\uD \in \bR$ is a
fixed constant, is then given as follows.
For $t\in(0,T)$ find $(\varphi,w) \in H^1(\mathcal M) \times H^1(\mathcal M)$
with $w = \uD$ on $\partial_D\mathcal{M}$ such that
\begin{subequations} \label{eq:weak}
\begin{align}
&
\vartheta \int_{\mathcal M} w_t \chi \dH2
+ \tfrac12\lambda \int_{\mathcal M} \varphi_t \chi \dH2
+ \int_{\mathcal M} \mathcal{K} \nabs w \cdot \nabs \eta \dH2 = 0
\qquad \forall\ \chi\in H^1_0(\mathcal{M}),
\label{eq:weaka} \\
&
\tfrac12 \cPsi a \int_{\mathcal M} w \eta \dH2
= \rho \epsilon \int_{\mathcal M} \mu(\cdot, \nabs\varphi) \varphi_t \eta \dH2
+ \alpha \epsilon \int_{\mathcal M} A_p(\cdot, \nabs\varphi) \cdot \nabs \eta 
\dH2 \nonumber \\ & \hspace{4cm}
+\alpha  \epsilon^{-1} \int_{\mathcal M} \Psi'(\varphi) \eta \dH2 \label{eq:weakb}
\qquad \forall\ \eta\in H^1(\mathcal{M}).
\end{align} 
\end{subequations}
Above, we used the function space
$$ H^1_0(\mathcal{M}) =\{\chi \in  H^1(\mathcal{M}) : \chi = 0\ \mbox{ on $\partial_D\mathcal{M}$} \}.
$$
\setcounter{equation}{0}
\section{Properties of anisotropic energies} \label{sec:anen}
\subsection{Minimal energy directions on the sphere}
If one wants to split the unit sphere $\bS^2$ into two parts with the same area with a minimal isotropic interfacial energy, then each great circle solves this variational problem. In case that the energy is anisotropic a great circle can
lower its interfacial energy by rotating the great circle on $\bS^2$. The following lemma discusses the problem of minimizing anisotropic interfacial energy
for the problem of splitting the sphere  $\bS^2$ into two parts with the same area.
\begin{lemma} \label{lem:greatcircle}
	Let $\mathcal M = \bS^2$ and let $\gamma(z, \cdot) = \gamma_0(\cdot)$ for all 
	$z \in \mathcal M$. Then the isoperimetric problem
	\begin{equation} \label{eq:isoperimetric}
		\min \left\{ \int_{\partial E} \gamma_0(\nu_{\partial E}) \dH1 :
		E \subset \mathcal{M}, \mathcal{H}^2(E) = \tfrac12 \mathcal{H}^2(\mathcal M)
		\right\}
	\end{equation}
	is solved by a half sphere $E$ whose boundary $\partial E$ is a great circle
	with constant outward geodesic normal 
	$\nu_{\partial E} \in \argmin_{\nu \in \mathcal{M}} \gamma_0(\nu)$.
\end{lemma}
\begin{proof}
As $\gamma_0$ is continuous and $\mathcal{M}$ is compact, we deduce
that $\gamma_0$ attains its minimum on $\mathcal{M}$.
	Let $\nu_{\min} \in \argmin_{\nu \in \mathcal{M}} \gamma_0(\nu)$, and let $H \subset
	\bR^3$ be the hyperplane orthogonal to $\nu_{\min}$, with $H_{\leq0}$ denoting
	the half space such that $\nu_{\min}$ is its outer normal on 
	$H = \partial H_{\leq0}$. Then the half sphere
	$E_{\min} = \mathcal M \cap H_{\leq0}$ solves
	\eqref{eq:isoperimetric}. To see this, we first 
	note that it is an admissible candidate due to 
	$\mathcal{H}^2(\mathcal M) = 4\pi$ for the unit sphere.
	In addition, it holds for any admissible 
	candidate $E$ in \eqref{eq:isoperimetric} that
$\mathcal{H}^1(\partial E) \geq \mathcal{H}^1(\partial E_{\min}) =
2\pi$ and hence
	\[
	\int_{\partial E} \gamma_0(\nu_{\partial E}) \dH1 
	\geq \mathcal{H}^1(\partial E) \gamma_0(\nu_{\min}) 
	\geq \mathcal{H}^1(\partial E_{\min}) \gamma_0(\nu_{\min}) 
	= 2\pi\gamma_0(\nu_{\min}),
	\]
	with equality for $E = E_{\min}$, since $\nu_{\partial E_{\min}} = \nu_{\min}$
	everywhere on $\partial E_{\min}$. 
\end{proof}

\subsection{Consistent 2d anisotropies on the unit sphere} \label{sec:cons2d}
In this section we discuss the idea of extending an anisotropy, 
which is initially defined on just one tangent space of a surface, 
in a consistent way to all the tangent spaces of the surface. 
We remark that we only need to define the anisotropy on the tangent spaces and 
not for all $\vec z \in \bR^3$.
Assume that at one reference point on the surface 
an anisotropy density function is given and in certain directions the anisotropy takes a minimal value.
Then we want to move the anisotropy along geodesics
to other tangent spaces. 

This is done in such a way that unit tangent vectors to the geodesic have the same anisotropic density. 
In fact, this idea can be used for general surfaces. However, we discuss the
idea for the unit sphere with the north pole as the reference point, since 
in this case explicit expressions can be more easily stated.
If, for example, on the north pole, we have a six-fold anisotropy, this will remain
the case on all the other tangent spaces as well. The procedure will lead to a spatially dependent anisotropy. 
We point out that for spatially homogeneous anisotropies of the form 
\eqref{eq:gamma0}, 
the strength and shape of the anisotropy may vastly differ 
between tangent spaces based at different points.
Hence the property that the anisotropy remains 
of the same structure (for example six-fold) will in general not hold true 
for these anisotropies. We will visualize this phenomenon with some 
numerical evidence later on.

Due to the hairy ball theorem 
it is not possible to move an anisotropic energy density in a continuous way along the unit sphere. 
In fact, we cannot choose tangents at all  points on the surface
for which the energy density is minimal in such a way that the tangent depends continuously on the point of the surface.
We thus have to consider a subset of the sphere.
Without loss of generality, assume that $\mathcal M \subset \bS^2$ 
with $\partial\mathcal M \not=\emptyset$ and such that 
$\vec e_3 \in \mathcal M$ but $-\vec e_3 \not\in \mathcal M$, so the open
manifold is a part of the unit sphere that contains the north pole but not 
the south pole.

Then for any $\vec z \in \mathcal M$ we will construct a rotation matrix
$Q(\vec z)$ that rotates the plane spanned by $\vec e_3$ and $\vec z$ around
the axis $\vec e_3 \times \vec z$ by the angle between $\vec e_3$ and $\vec z$,
so that $Q(\vec z) \vec z = \vec e_3$. The same transformation then also
rotates the tangent space $\tanspace_{\vec z} \mathcal{M}$, and we can evaluate
the anisotropy in $Q(\vec z) \tanspace_{\vec z} \mathcal{M} = 
\tanspace_{\vec e_3} \mathcal{M} = \bR^2 \times \{0\}$.
Then we consider
\begin{equation} \label{eq:gammasphere}
	\gamma(z, p) = \gamma_0 (Q(\vec z) \vec p).
\end{equation}
Here we have defined an anisotropic density 
$\gamma_0 : \bR^2 \times \{0\} \to \bR$, via
$\gamma_0(p_1,p_2,p_3) = \hat\gamma(p_1,p_2)$ for an
$\hat\gamma : \bR^2 \to \bR$ that is one-homogeneous.

It is easy to check that the following matrix has the desired properties for 
$Q(\vec z)$:
\begin{equation} \label{eq:Q}
	Q(\vec z) =
	\begin{pmatrix}
		z_3 + \frac{z_2^2}{1+z_3} & -\frac{z_1 z_2}{1+z_3} & -z_1 \\
		-\frac{z_1 z_2}{1+z_3} & z_3 + \frac{z_1^2}{1+z_3} & -z_2 \\
		z_1 & z_2 & z_3 
	\end{pmatrix} \qquad \vec z \in \bS^2 \setminus\{-\vec e_3\}.
\end{equation}
In particular, it is easy to check that
\[
Q(\vec z) Q(\vec z)^T = \Id,\quad
Q(\vec z) \vec z = \vec e_3,\quad
Q(\vec z) (\vec e_3 \otimes \vec z) = \vec e_3 \otimes \vec z.
\]

\subsection{BGN-type anisotropies}
For the numerical analysis, we restrict ourselves to the following anisotropies of BGN-type
\begin{equation} \label{eq:g1}
	\gamma(z, \vec{p}) = \sum_{\ell=1}^L \gamma_\ell(z,\vec p)
	= \sum_{\ell=1}^L [{\vec{p}\cdot G_{\ell}(z)\vec{p}}]^\frac12,
	\qquad \forall\ \vec p \in \bR^3, z \in \mathcal{M},
\end{equation}
where $G_{\ell} : \mathcal{M} \to \bR^{3\times 3}$, for $\ell=1,\ldots,L$, 
are functions that map to symmetric and positive definite matrices.
For the case that $G_\ell$ are constant matrices, these anisotropies
have first been introduced by Barrett, Garcke and N\"urnberg in 
\cite{triplejANI,ani3d}.

\begin{lemma}\label{lem:1}
	Let $\gamma$ be of the form (\ref{eq:g1}). Then $\gamma(z,\cdot)$ is convex 
	for every $z \in \mathcal{M}$ and the anisotropic operator $A$ satisfies
\begin{subequations}
	\begin{alignat}{2} \label{eq:lem1a}
		A_p(z,\vec p) \cdot (\vec p-\vec q) & \geq 
		\gamma(z,\vec p)[\gamma(z,\vec p)-\gamma(z,\vec q) ]
		\qquad && \forall\ \vec p \in \bR^3\setminus\{\vec 0\}, \vec q \in \bR^3, \\
		A(z,\vec p) & \leq \tfrac12\gamma(z,\vec q)
		\sum_{\ell=1}^L [\gamma_\ell(z,q)]^{-1}[\gamma_\ell(z,p)]^2 
		\qquad && \forall\ \vec p \in \bR^3, \vec q \in \bR^3\setminus\{\vec 0\}
		\label{eq:lem1b}
	\end{alignat}
\end{subequations}
	for every $z \in \mathcal{M}$.
\end{lemma}
\begin{proof}
	The result follows pointwise from \cite[Lemma~2.1]{eck}.
\end{proof}

Following \cite{eck,vch}, we let for every $z \in \mathcal{M}$
\begin{equation*} 
	B(z, \vec q) := \begin{cases}
		\gamma(z, \vec q)\displaystyle\sum_{\ell = 1}^L 
		[\gamma_\ell(z,\vec q)]^{-1}G_\ell(z)
		& \vec q \not= \vec 0, \\
		L\displaystyle\sum_{\ell = 1}^L G_\ell(z) & \vec q = \vec 0.
	\end{cases}
\end{equation*}
As $A(z, p) = \tfrac12 \gamma^2(z, p)$, it clearly holds that
$$
B(z, \vec p)\vec p = A_p(z,\vec p) 
\qquad \forall\ \vec p \in \bR^3\setminus\{\vec 0\},
$$
and it turns out that approximating $A_p(z,\vec p)$ with $B(z,\vec q)\vec p$
maintains the monotonicity property \eqref{eq:lem1a}.

\begin{lemma} \label{lem:B}
	Let $\gamma$ be of the form (\ref{eq:g1}). Then it holds that
	\begin{equation} \label{eq:Bmono}
		[B(z,\vec q)\vec p]\cdot(\vec p-\vec q) 
		\geq \gamma(z,p)\left[\gamma(z,p)-\gamma(z,q)\right] 
		\qquad \forall\ \vec p , \vec q \in \bR^3 ,
	\end{equation}
	for every $z \in \mathcal{M}$.
	
\end{lemma}
\begin{proof}
	The result follows pointwise from \cite[Lemma~2.2]{eck}.
\end{proof}

We observe that on letting
\begin{equation} \label{eq:hatg1}
	\hat\gamma(\vec{p}) 
	= \sum_{\ell=1}^L [{\vec{p}\cdot \hat G_{\ell}\vec{p}}]^\frac12,
	\qquad \forall\ \vec p \in \bR^2,
\end{equation}
the anisotropy \eqref{eq:gammasphere} with \eqref{eq:Q} falls into the category
\eqref{eq:g1} with the special choices
\begin{equation*} 
	G_{\ell}(z) = Q^T(z) \binom{\hat G_\ell\ 0}{0\ 1} Q(z).
\end{equation*}

\setcounter{equation}{0}
\section{Finite element approximation} \label{sec:fea}

Let $\mathcal{M}^h$ be a polyhedral hypersurface approximating $\mathcal{M}$,
and let $\{{\cal T}^h\}_{h>0}$ be a family of open triangles with
$\mathcal{M}^h=\cup_{\sigma\in{\cal T}^h}\overline{\sigma}$. We refer to \cite{bgnreview} for more details on polyhedral approximations of surfaces and finite element spaces on polyhedral surfaces.
Associated with ${\cal T}^h$ is the finite element space
\begin{equation*} 
 S^h = \{\chi \in C(\mathcal{M}^h) : \chi \mid_{\sigma} \mbox{ is affine }
 \forall\ \sigma \in {\cal T}^h\}.
\end{equation*}
Let $J$ be the set of nodes of ${\cal T}^h$ and $\{p_{j}\}_{j \in J}$ the
coordinates of these nodes.
Let $\{\chi_{j}\}_{j\in J}$ be the standard basis
functions for $S^h$; that is $\chi_{j} \in S^h$ and $\chi_j(p_{i})=\delta_{ij}$
for all $i,j \in J$.
A discrete semi-inner product 
for functions that are piecewise continuous on $\mathcal{T}^h$ 
can be defined by
\begin{equation} \label{eq:dip}
\langle \eta_1, \eta_2 \rangle^h_{\mathcal{M}^h}=
\sum_{\sigma \in\mathcal{T}^h} \langle \eta_1, \eta_2 \rangle^h_\sigma
= \sum_{\sigma \in\mathcal{T}^h} \tfrac13 |\sigma|
\sum_{k=0}^{2} (\eta_1\eta_2) ( ({p}_{j_k})^-) ,
\end{equation}
with $\{{p}_{j_k}\}_{k=0}^{2}$ denoting the vertices of $\sigma$,
and where we define $\eta(({p}_{j_k})^-)=
\underset{\sigma\ni {q}\to {p}_{j_k}}{\lim} \eta({q})$,
$k=0,\ldots,2$.
We note that (\ref{eq:dip}) induces the discrete semi-norm 
$|\eta|_h := [\langle\eta,\eta\rangle^h_{\mathcal{M}^h}]^{\frac{1}{2}}$ 
on $L^\infty(\mathcal{M}^h)$, that becomes a norm on $S^h$.

We introduce also
\begin{align*} 
 K^h & = \{\chi \in S^h : |\chi| \leq 1 \}, \nonumber\\
 S^h_0 & = \{\chi \in S^h : \chi = 0 \ \mbox{ on $\partial_D\mathcal{M}^h$} \}  
\quad \mbox{and} \quad
 S^h_D = \{\chi \in S^h : \chi = \uD\ \mbox{ on $\partial_D\mathcal{M}^h$} \} ,
\end{align*}
where $\partial_D\mathcal{M}^h \subset \partial\mathcal{M}^h$ 
is a suitable approximation of $\partial_D\mathcal{M}$.

In addition to ${\cal T}^h$, let
$0= t_0 < t_1 < \ldots < t_{N-1} < t_N = T$ be a
partitioning of $[0,T]$ into possibly variable time steps $\tau_n := t_n -
t_{n-1}$, $n=1,\ldots, N$. 

\subsection{The obstacle potential}

Let $\Phi^0 \in K^h$ be an approximation of $\varphi(0)$. 
Similarly, if $\vartheta>0$ let $W^0 \in S^h_D$ be an approximation of 
$w(0)$.
Then, for $n \geq 1$, find $(\Phi^n,W^n) \in K^h \times S^h_D$ such that
\begin{subequations} \label{eq:fea}
\begin{align}
& \vartheta \left\langle\dfrac{W^{n}-W^{n-1}}{\tau_n},
 \chi \right\rangle_{\mathcal{M}^h}^h + 
\tfrac12\lambda\left\langle\dfrac{\Phi^{n}-\Phi^{n-1}}{\tau_n},
 \chi \right\rangle_{\mathcal{M}^h}^h 
+  \left\langle\mathcal{K}\nabs W^{n},\nabs\chi \right\rangle_{\mathcal{M}^h}^h
 = 0 \quad \forall\ \chi \in S^h_0, \label{eq:feaU} \\
& 
\epsilon\frac\rho\alpha\left\langle\mu(\cdot, \nabs\Phi^{n-1})
\dfrac{\Phi^{n}-\Phi^{n-1}}{\tau_n}, \chi - \Phi^n
\right\rangle_{\mathcal{M}^h}^h +
\epsilon\left\langle B(\cdot,\nabs\Phi^{n-1})\nabs\Phi^{n}, \nabs [\chi -\Phi^n]
\right\rangle_{\mathcal{M}^h}^h
\nonumber \\ & \qquad\qquad
 \geq \left\langle\tfrac12 \cPsi \frac{a}\alpha W^{n} 
 + \epsilon^{-1}\Phi^{n-1},\chi - \Phi^n\right\rangle_{\mathcal{M}^h}^h
 \qquad \forall\ \chi \in K^h. \label{eq:feaW}
\end{align}
\end{subequations}

\begin{remark} \label{rem:fea}
It is possible to generalize the considered model in the following ways.
All of the results presented in this paper remain valid for these
generalizations, on using the techniques developed by the authors in 
\cite{vch}. 
\begin{itemize}
\item It is possible to consider a phase-dependent $\mathcal{K}$, e.g., by
defining 
$b(s) = \tfrac12(1+s)\mathcal{K}_+ + \tfrac12(1-s)\mathcal{K}_-$, 
and replacing $\mathcal{K}$ in \eqref{eq:weaka} by $b(\varphi)$.
\item It is possible to consider functions $\varrho$ that accelerate the
convergence of the phase field model to the sharp interface limit. To this end,
the two factors $\frac12$ in \eqref{eq:weak} are replaced by the coefficients
$\varrho(\varphi)$ in the two associated integrals, where, e.g.,
\begin{equation} \label{eq:varrho}
\text{(i)}\ 
\varrho(s) = \tfrac12\,,\qquad
\text{(ii)}\ 
\varrho(s) = \tfrac12\,(1 - s)\,,\qquad
\text{(iii)}\ 
\varrho(s) = \tfrac{15}{16}\,(s^2 - 1)^2\,.
\end{equation}
\item It is possible to extend the numerical analysis to a more general
family of anisotropies. In particular, in place of \eqref{eq:g1} the family of
densities
\[
\gamma(z, \vec{p}) = \left(\sum_{\ell=1}^L |\gamma_\ell(z,\vec p)|^r 
\right)^{\frac1r}, \qquad r \in [1,\infty),
\]
can be considered, so that \eqref{eq:g1} collapses to the case $r=1$.
In the case $r > 1$, the term $B(\cdot,\nabs\Phi^{n-1})$ in, e.g.,
\eqref{eq:feaW} needs to be changed to
$B_r(\cdot,\nabs\Phi^{n-1},\nabs\Phi^{n})$,
where the matrices $B_r$ are defined analogously to the flat and spatially
homogeneous case in \cite{vch}.
\end{itemize}
\end{remark}

Let 
\begin{equation*} 
\mathcal{E}^h(W, \Phi) = 
\frac\vartheta2\,|W - \uD|_h^2 +
\frac{\lambda\alpha}{a}\frac1\cPsi
\left[ \tfrac12\epsilon|\gamma(\cdot,\nabs\Phi)|_h^2 + 
\epsilon^{-1} \langle\Psi(\Phi), 1\rangle^h_{\mathcal{M}^h} \right] \,,
\end{equation*}
and define
\begin{equation*} 
\mathcal{F}^h(W, \Phi) = \mathcal{E}^h(W, \Phi)
- \lambda\uD \langle \tfrac12, 1+\Phi\rangle_{\mathcal{M}^h}^h
\end{equation*}
for all $W, \Phi \in S^h$.

\begin{theorem} \label{thm:stab2}
Let $\gamma$ be of the form \eqref{eq:g1}.
Then there exists a solution $(\Phi^n,W^n) \in K^h \times S^h_D$ to 
\eqref{eq:fea}, where $\Phi^n$ is unique and $W^n$ is unique up to 
an additive constant.
If $\vartheta > 0$ or $\partial_D\mathcal M^h\not=\emptyset$
then $W^n$ is unique.
If $\vartheta= 0$ and $\partial_D\mathcal M^h=\emptyset$, then 
$W^n$ is unique if there exists a $j\in J$ such that
$|\Phi^n(p_j)| < 1$.
In addition, it holds that any solution $(\Phi^n,W^n) \in K^h \times S^h_D$ to 
\eqref{eq:fea} satisfies the stability bound
\begin{align} \label{eq:stab3}
& \mathcal{F}^h(W^n, \Phi^n)  
+ \tau_n \langle \mathcal{K}\nabs W^{n} , \nabs W^n \rangle_{\mathcal{M}^h}^h 
+ \tau_n\,\frac{\lambda\,\rho}{a}\,\frac\epsilon\cPsi
\left| [\mu(\cdot, \nabs\Phi^{n-1})]^{\frac12}\, 
\dfrac {\Phi^{n}-\Phi^{n-1}}{\tau_n}\right|_h^2 
\nonumber \\ & \hspace{9cm}
\leq \mathcal{F}^h(W^{n-1}, \Phi^{n-1})
\,.
\end{align}
\end{theorem}
\begin{proof}
The existence and uniqueness results follow analogously to 
Lemma~3.1 in \cite{vch}, and so we omit their proof here. 

Choosing $\chi = W^n - \uD$ in (\ref{eq:feaU}) and 
$\chi = \Phi^{n-1}$ in (\ref{eq:feaW}) yields that
\begin{subequations} \label{eq:stab2}
\begin{align}
& \vartheta \langle W^n - W^{n-1}, W^n - \uD \rangle_{\mathcal{M}^h}^h 
+ \tfrac12\lambda \langle\Phi^{n}-\Phi^{n-1}, W^n - \uD 
\rangle_{\mathcal{M}^h}^h 
+ \tau_n \langle \mathcal{K}\nabs W^n, \nabs W^n \rangle_{\mathcal{M}^h}^h
\nonumber \\ & \hspace{5cm}
 = 0 \,, \label{eq:stab2u} \\
& 
\epsilon\frac\rho\alpha\tau_n^{-1}
\langle\mu(\cdot,\nabs\Phi^{n-1})[\Phi^{n}-\Phi^{n-1}], 
\Phi^{n-1} - \Phi^n\rangle_{\mathcal{M}^h}^h + \epsilon
 \langle B(\cdot,\nabs \Phi^{n-1})\nabs\Phi^{n}, \nabs [\Phi^{n-1} -\Phi^n]
\rangle_{\mathcal{M}^h}^h \nonumber \\ & \hspace{5cm}
 \geq 
 \langle\tfrac12\cPsi\frac{a}\alpha W^{n} 
 + \epsilon^{-1}\Phi^{n-1}, \Phi^{n-1} - \Phi^n \rangle_{\mathcal{M}^h}^h. 
\label{eq:stab2w}
\end{align}
\end{subequations}
It follows from \eqref{eq:stab2}, on noting the elementary identity 
\begin{equation*} 
2y(y-z) = y^2 - z^2 + (y-z)^2 \qquad \forall\ y,z \in \bR,
\end{equation*}
and on recalling \eqref{eq:Bmono}, that
\begin{align*}
& \tfrac12\epsilon|\gamma(\cdot,\nabs\Phi^n)|_0^2 
- \tfrac12\epsilon^{-1}|\Phi^n|_h^2 + 
\frac\vartheta2\frac{a}{\lambda\,a}\cPsi|W^n - \uD|_h^2 +
\tau_n\epsilon\frac\rho\alpha\left|
[\mu(\cdot,\nabs\Phi^{n-1})]^{\frac12}
\dfrac {\Phi^{n}-\Phi^{n-1}}{\tau_n}\right|_h^2 
\nonumber \\ & \hspace{2cm}
- \uD\frac{a}\alpha\cPsi\,
 \langle\tfrac12, \Phi^n - \Phi^{n-1} \rangle_{\mathcal{M}^h}^h
+ \tau_n\frac{a}{\lambda\,\alpha}\cPsi
\langle \mathcal{K}\nabs W^{n} , \nabs W^n \rangle_{\mathcal{M}^h}^h
\nonumber \\ & \hspace{5cm}
\leq \tfrac12\epsilon|\gamma(\cdot,\nabs\Phi^{n-1})|_0^2 
- \tfrac12 \epsilon^{-1} |\Phi^{n-1}|_h^2 +
\frac\vartheta2\frac{a}{\lambda\,a}\cPsi|W^{n-1} - \uD|_h^2 \,.
\end{align*}
This yields the desired result (\ref{eq:stab3}) on adding the constant
$\frac12\,\epsilon^{-1}\,\int_{\mathcal{M}^h} 1 \dH2$ on both sides, 
and then multiplying the inequality with 
$\frac{\lambda\,\alpha}{a\cPsi}$.
\end{proof}

\subsection{Smooth potentials} 
The unconditionally stable approximation \eqref{eq:fea} for the obstacle
potential (\ref{eq:obs}) can be easily adapted to the case of a smooth
potential such as (\ref{eq:quartic}). Such approximations rely on
a convex/concave splitting of $\phi = \Psi'$, i.e.,
$\phi = \phi^+ + \phi^-$ with
\begin{equation*} 
\pm (\phi^\pm)'(s) \geq 0 \quad \forall\ s \in \bR\,.
\end{equation*}
For the quartic potential (\ref{eq:quartic}) the natural choices are
\begin{equation} \label{eq:phi+-}
\phi^+(s) = s^3 \quad \text{and} \quad \phi^-(s) = -s\,, 
\end{equation}
and for simplicity we restrict our attention to that case. Details on how
to deal with more general potentials can be found in \cite{vch}. 

As before, given $\Phi^0 \in K^h$ and, if $\vartheta>0$, $W^0 \in S^h_D$, 
for $n \geq 1$, find $(\Phi^n,W^n) \in S^h \times S^h_D$ such that
\begin{subequations} \label{eq:qfea}
\begin{align}
& \vartheta \left\langle\dfrac{W^{n}-W^{n-1}}{\tau_n},
 \chi \right\rangle_{\mathcal{M}^h}^h + 
\tfrac12\lambda\left\langle\dfrac{\Phi^{n}-\Phi^{n-1}}{\tau_n},
 \chi \right\rangle_{\mathcal{M}^h}^h 
+  \left\langle\mathcal{K}\nabs W^{n},\nabs\chi \right\rangle_{\mathcal{M}^h}^h
 = 0 \quad \forall\ \chi \in S^h_0, \label{eq:qfeaU} \\
& 
\epsilon\frac\rho\alpha\left\langle\mu(\cdot, \nabs\Phi^{n-1})
\dfrac{\Phi^{n}-\Phi^{n-1}}{\tau_n}, \chi
\right\rangle_{\mathcal{M}^h}^h +
\epsilon\left\langle B(\cdot,\nabs\Phi^{n-1})\nabs\Phi^{n}, \nabs \chi
\right\rangle_{\mathcal{M}^h}^h
+ \epsilon^{-1} \left\langle \phi^+(\Phi^n),\chi \right\rangle_{\mathcal{M}^h}^h
\nonumber \\ & \qquad\qquad
 = \left\langle\tfrac12 \cPsi \frac{a}\alpha W^{n} 
 - \epsilon^{-1}\phi^-(\Phi^{n-1}),\chi \right\rangle_{\mathcal{M}^h}^h
 \qquad \forall\ \chi \in S^h. \label{eq:qfeaW}
\end{align}
\end{subequations}

\begin{theorem} \label{thm:qstab}
Let $\gamma$ be of the form \eqref{eq:g1}.
Let $\Psi$ be given by \eqref{eq:quartic} and let \eqref{eq:phi+-} hold.
Then there exists a unique solution
$(\Phi^n,W^n) \in S^h \times S^h_D$ to \eqref{eq:qfea}. Moreover, the solution
satisfies the stability bound \eqref{eq:stab3}. 
\end{theorem}
\begin{proof}
The existence and uniqueness results follow analogously to
Theorem~3.8 in \cite{vch}. The proof of the stability bound \eqref{eq:stab3} 
is similar to the proof of Theorem~\ref{thm:stab2}, making use of the
convex/concave splitting $\Psi' = \phi = \phi^+ + \phi^-$.
\end{proof}

\setcounter{equation}{0}
\section{Numerical results} \label{sec:nr}

We implemented the scheme \eqref{eq:fea} with the help of the finite element 
toolbox ALBERTA, see \cite{Alberta}. 
To increase computational efficiency, we employ adaptive meshes, which have a 
finer mesh size $h_{f} \approx \frac{1}{N_f}$ within the diffuse interfacial 
regions and a coarser mesh size $h_{c} \approx \frac{1}{N_c}$ away from them, 
with $N_f, N_c \in \bN$, see \cite{voids3d,voids} for a more detailed 
description in the planar case. For the solution of the nonlinear systems of
equations we employ the iterative solvers discussed in \cite{vch}.

Unless otherwise stated, for the physical parameters in \eqref{eq:general} 
we choose
$\vartheta = \rho = 0$, $\lambda = \mathcal K = a = \alpha = 1$, $\mu=1$ 
and \eqref{eq:obs} for the potential $\Psi$.
In addition, we employ
uniform time steps $\tau_n = \tau$, $n = 1,\ldots, N$.

\subsection{Spatially inhomogeneous anisotropies in 2d} 

In this section we consider the special case
$\mathcal M \subset \bR^2 \times \{0\} \subset \bR^3$, i.e., we reformulate
evolutions in $\bR^2$ within the framework of this paper.

\begin{figure}
\center
\includegraphics[angle=-0,width=0.3\textwidth]{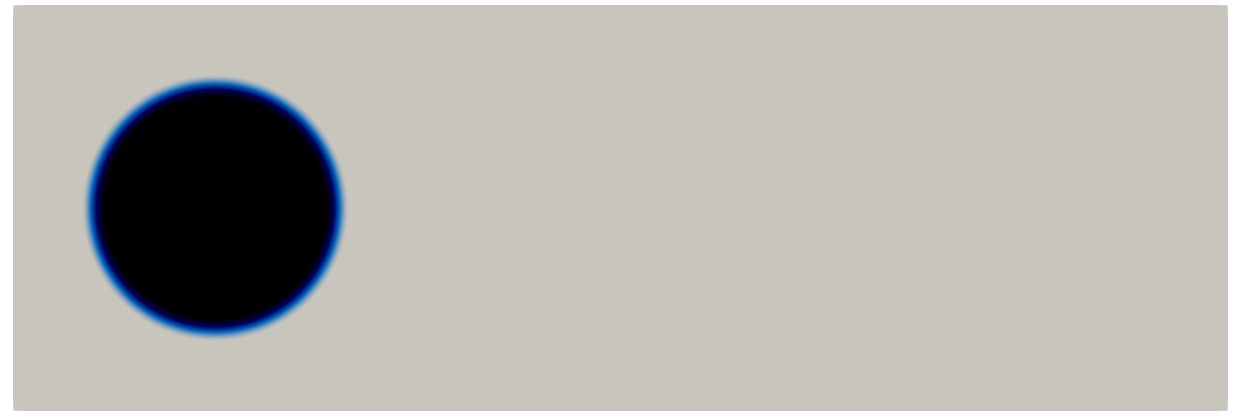}
\includegraphics[angle=-0,width=0.3\textwidth]{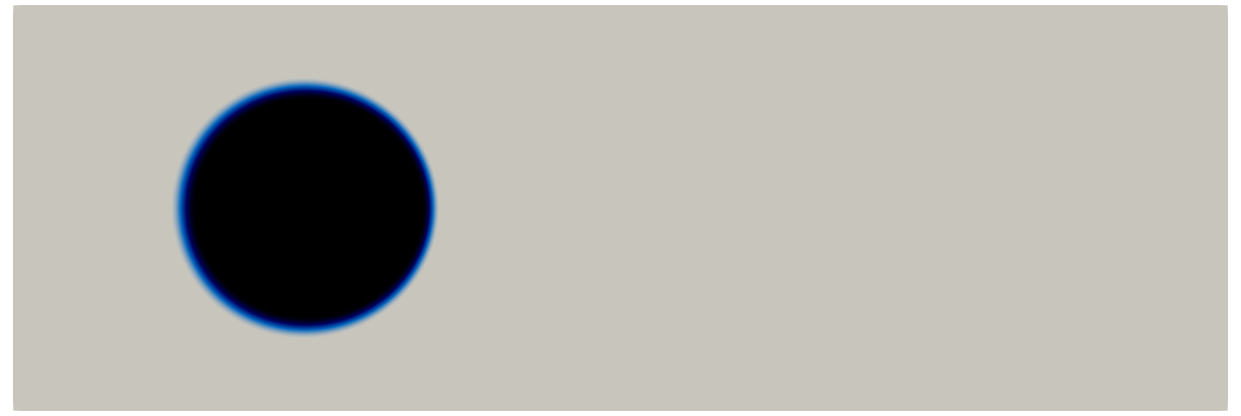}
\includegraphics[angle=-0,width=0.3\textwidth]{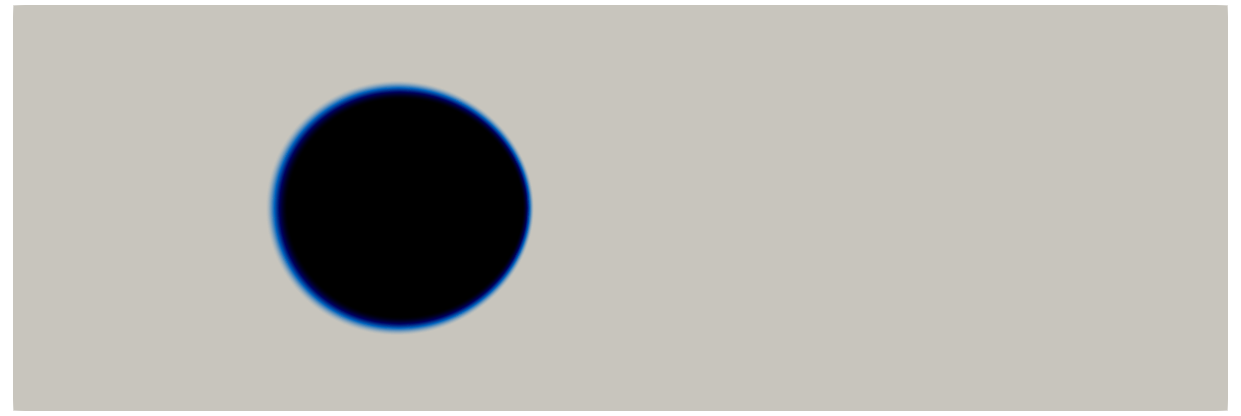}
\includegraphics[angle=-0,width=0.3\textwidth]{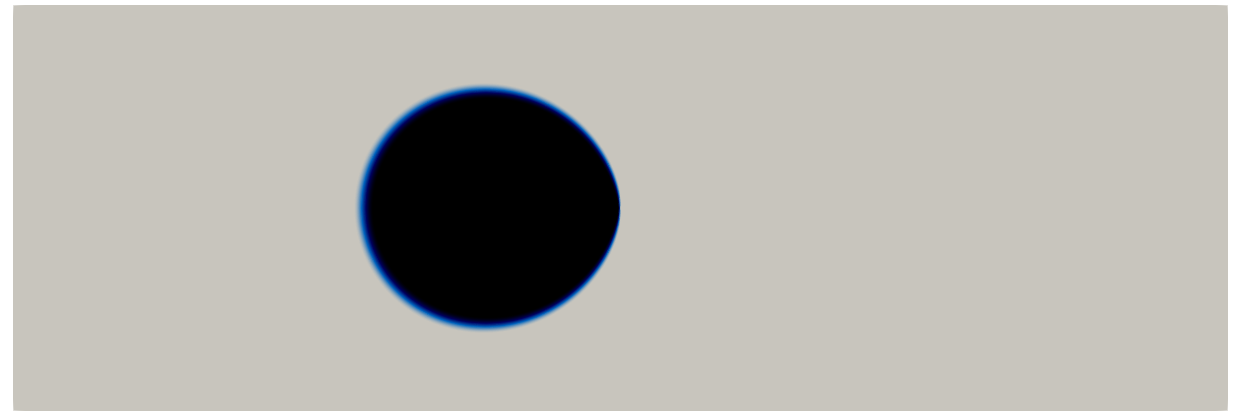}
\includegraphics[angle=-0,width=0.3\textwidth]{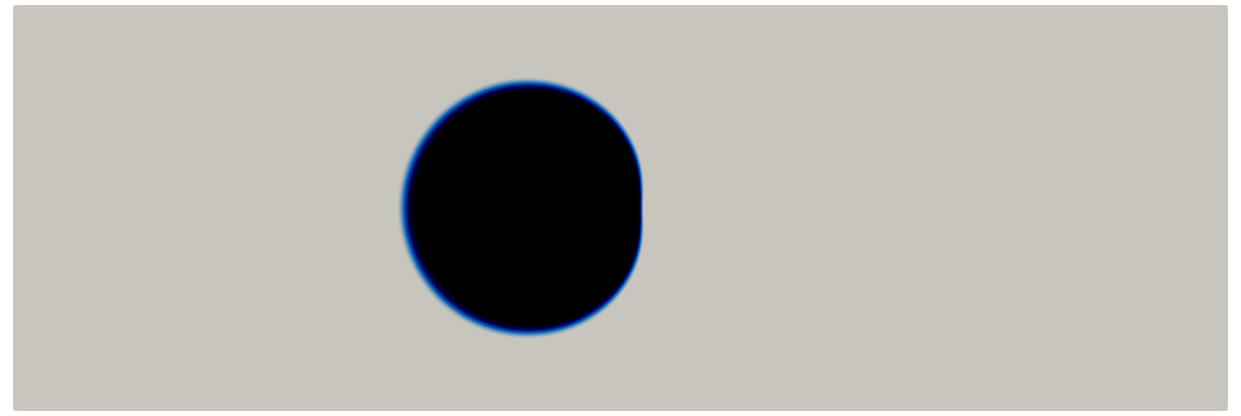}
\includegraphics[angle=-0,width=0.3\textwidth]{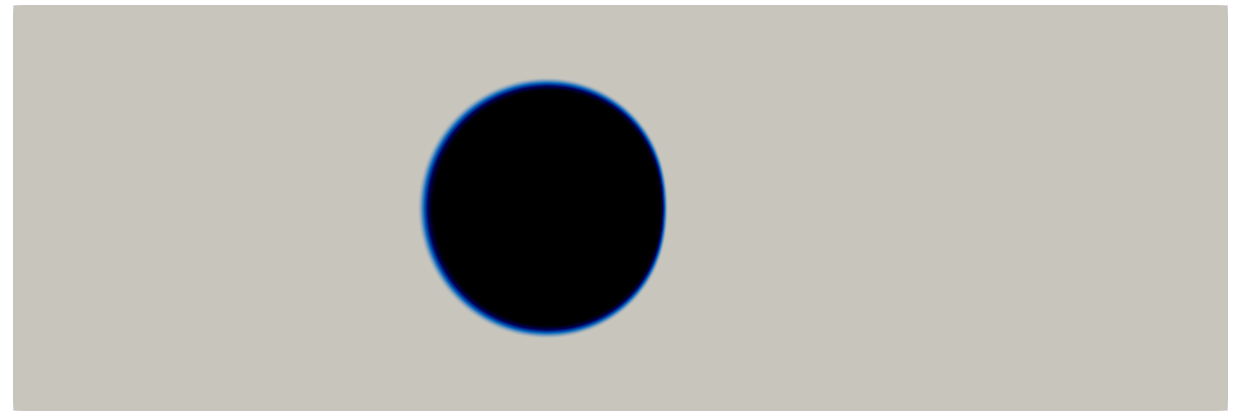}
\includegraphics[angle=-0,width=0.3\textwidth]{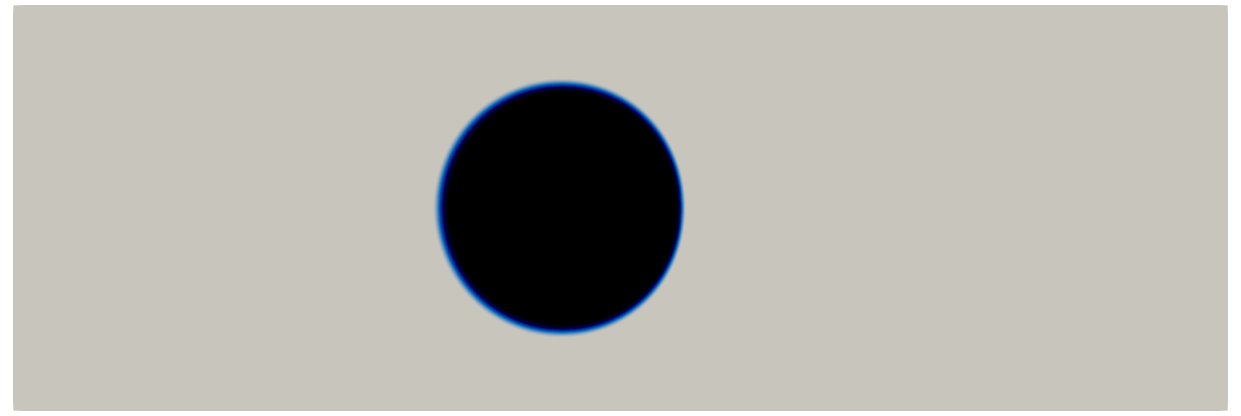}
\includegraphics[angle=-0,width=0.3\textwidth]{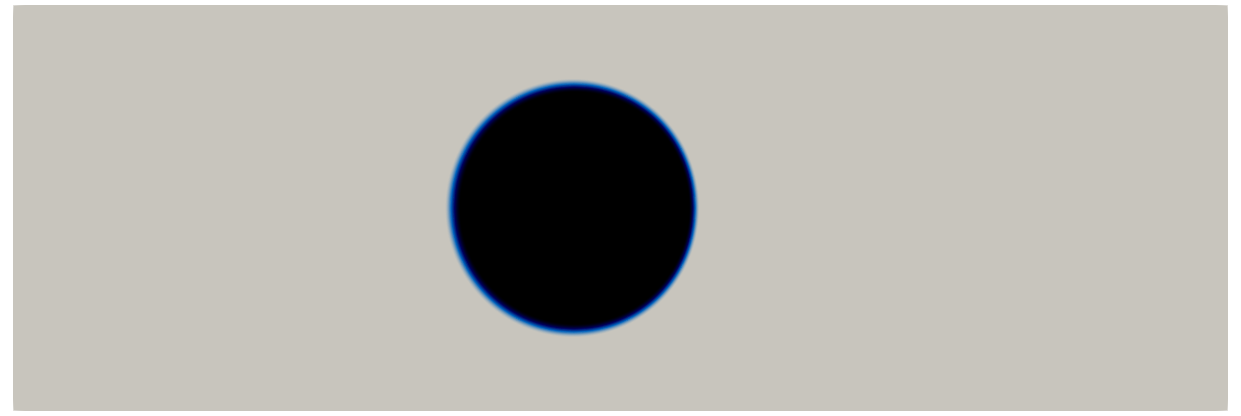}
\includegraphics[angle=-0,width=0.3\textwidth]{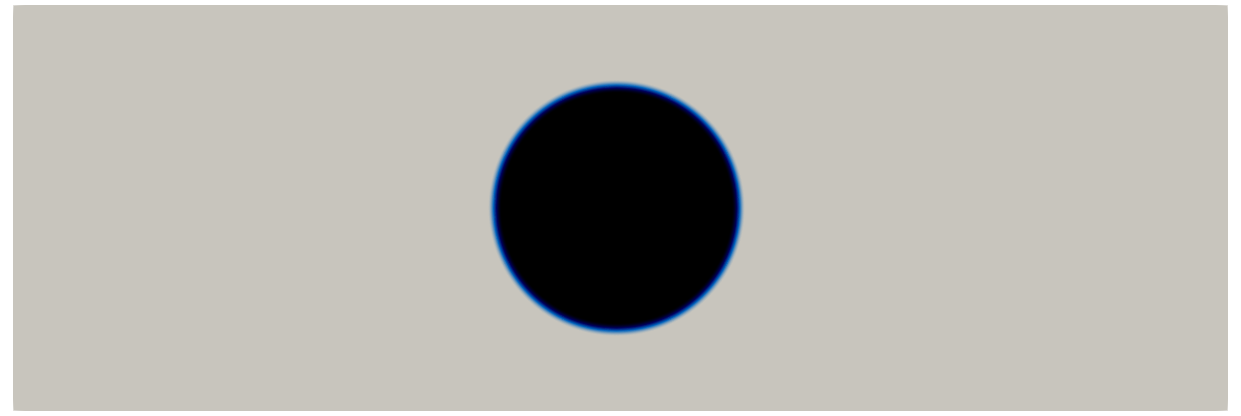}
\caption{($\epsilon=(16\pi)^{-1}$) 
Solutions for the Cahn--Hilliard equation in the case of a spatially inhomogeneous 
energy in the plane.
Snapshots of the evolutions at times $t=0, 0.1,0.2,0.3,0.4,0.5,0.6,0.7, 2$.
}
\label{fig:sink}
\end{figure}

Then, on 
$\mathcal M = (-\frac32,\frac32) \times (-\frac12,\frac12) \times \{0\}$ 
we use the spatially weighted isotropic surface energy 
\begin{equation*} 
\gamma(z, p) = (0.01 + |z|) |p|.
\end{equation*}
As the initial interface we choose a circle of radius $0.3$ centred at 
$(-1,0)^T$. Clearly, the interface can reduce its energy by moving towards a
circle centred at the origin, which we conjecture to be the global minimizer
for this setting.
The evolution for the solutions from our scheme \eqref{eq:fea} 
can be seen in Figure~\ref{fig:sink}. Here we let
$\epsilon=(16\pi)^{-1}$, $N_c = 32$, $N_f = 512$ and $\tau = 10^{-4}$.

\subsection{Spatially homogeneous anisotropies in 3d}

In this section we consider anisotropies of the form \eqref{eq:gamma0}.

\subsubsection{Simulations for Lemma~\ref{lem:greatcircle}}
\begin{figure}
\center
\includegraphics[angle=-0,width=0.3\textwidth]{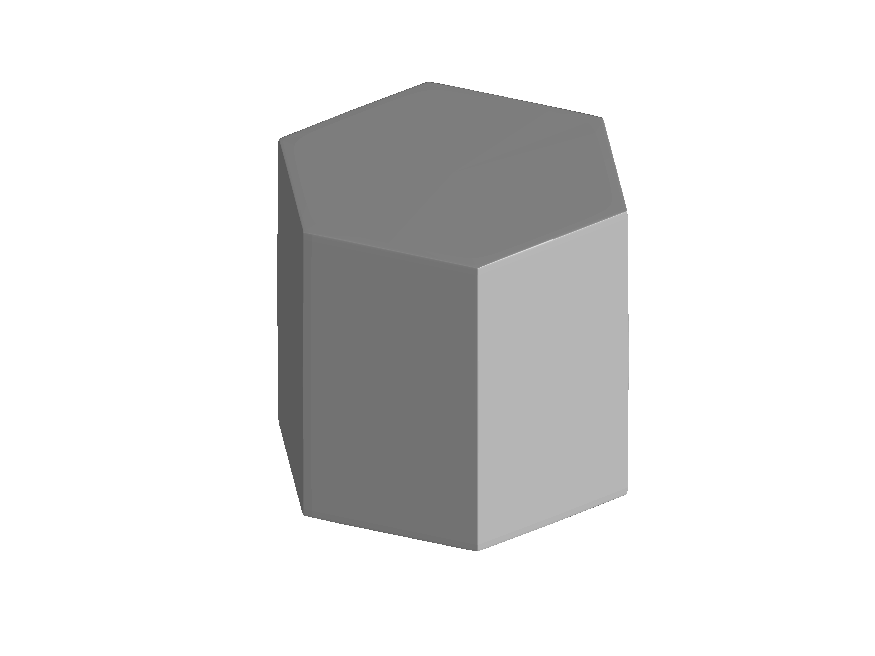}
\includegraphics[angle=-0,width=0.3\textwidth]{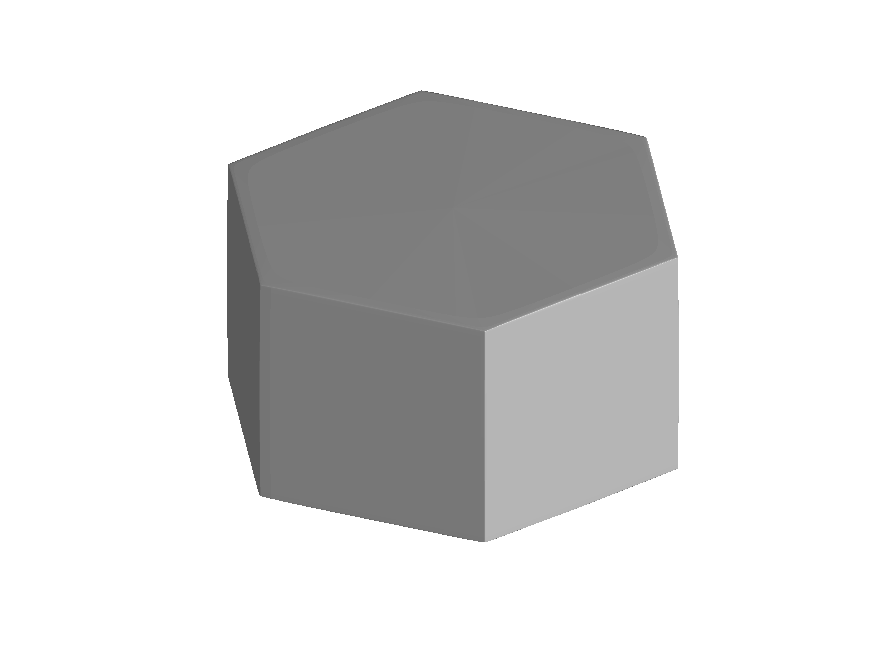}
\includegraphics[angle=-0,width=0.3\textwidth]{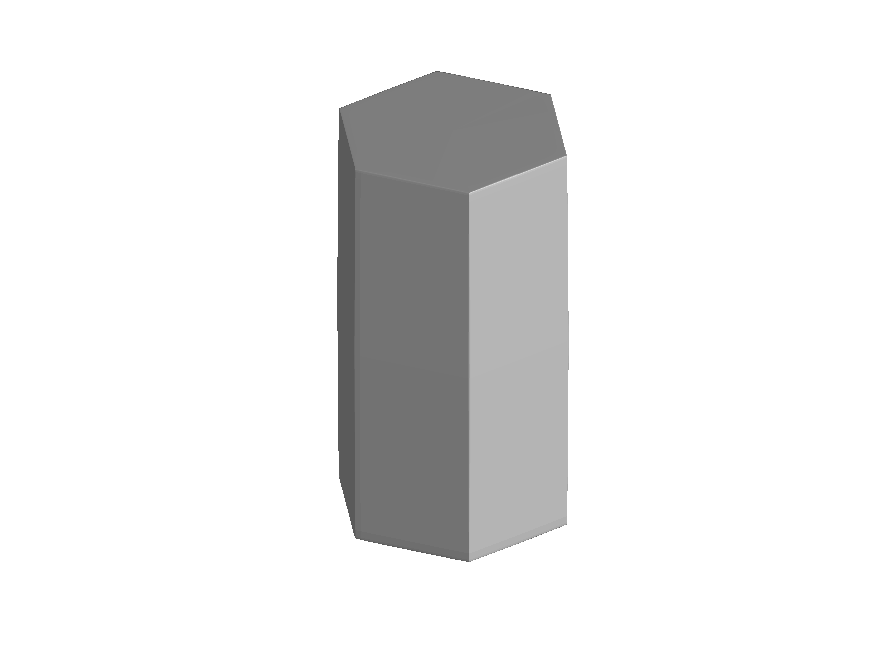}
\caption{The Wulff shapes for \eqref{eq:L4abc}(a), \eqref{eq:L4abc}(b) and
\eqref{eq:L4abc}(c).}
\label{fig:wulff3d}
\end{figure}

We choose as anisotropy the density
\begin{equation} \label{eq:L4abc}
\gamma(\cdot,p) = 
\gamma_0(\vec{p}) = 
l_\delta(R_2(\tfrac{\pi}2)\,\vec{p}) + w
\sum_{\ell = 1}^3
l_\delta(R_1(\tfrac{\ell\,\pi}3)\,\vec{p}),\ \delta=0.01,\
\begin{cases} 
w = \tfrac{1}{\sqrt{3}} & \text{(a)} \\
w = 1 & \text{(b)} \\
w = \tfrac{1}{2\sqrt{3}} & \text{(c)}
\end{cases},
\end{equation}
where $R_{1}(\theta)=\left(\!\!\!\scriptsize
\begin{array}{rrr} \cos\theta & \sin\theta&0 \\
-\sin\theta & \cos\theta & 0 \\ 0 & 0 & 1 \end{array}\!\! \right)$, 
$R_{2}(\theta)=\left(\!\!\!\scriptsize
\begin{array}{rrr} \cos\theta & 0 & \sin\theta \\
0 & 1 & 0 \\ -\sin\theta & 0 & \cos\theta \end{array}\!\! \right)$
and $l_\delta(\vec{p}) = \left[ \delta^2\,|\vec{p}|^2 
+ p_1^2\,(1-\delta^2) \right]^{\frac12}$.
The Wulff shapes for \eqref{eq:L4abc} are shown in Figure~\ref{fig:wulff3d}.
We recall from \cite{FonsecaM91} that for a given anisotropy $\gamma_0$
the boundary of its Wulff shape is the solution of the isoperimetric
problem for $\mathcal E_0(\Gamma) = \int_\Gamma \gamma_0(\nu) \dH2$.
Hence they are often used to visualize the properties of $\gamma_0$.
We note that for $\delta=0$, the densities \eqref{eq:L4abc}
represent crystalline anisotropies, whose Wulff shapes are given by hexagonal
prisms, cf.\ \cite[Fig.\ 3]{crystal}.
The coefficients in \eqref{eq:L4abc}(a) are chosen such that, in the limit
$\delta=0$, the eight admissible normal directions all have the same energy, 
while in \eqref{eq:L4abc}(b) and \eqref{eq:L4abc}(c) the two normals $\pm e_3$ 
have a smaller/larger energy than the other six directions.

\begin{figure}
\center
\includegraphics[angle=-0,width=0.3\textwidth]{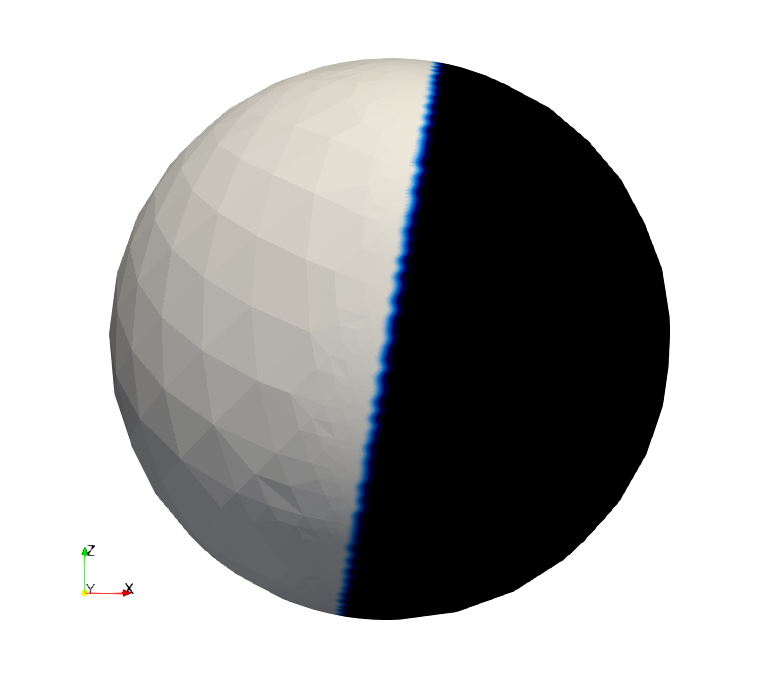}
\includegraphics[angle=-0,width=0.3\textwidth]{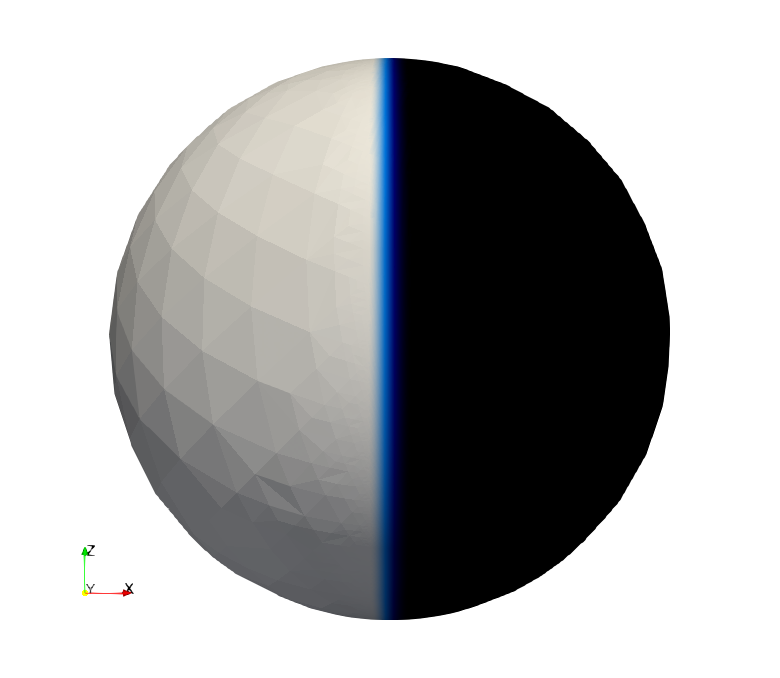}
\includegraphics[angle=-0,width=0.3\textwidth]{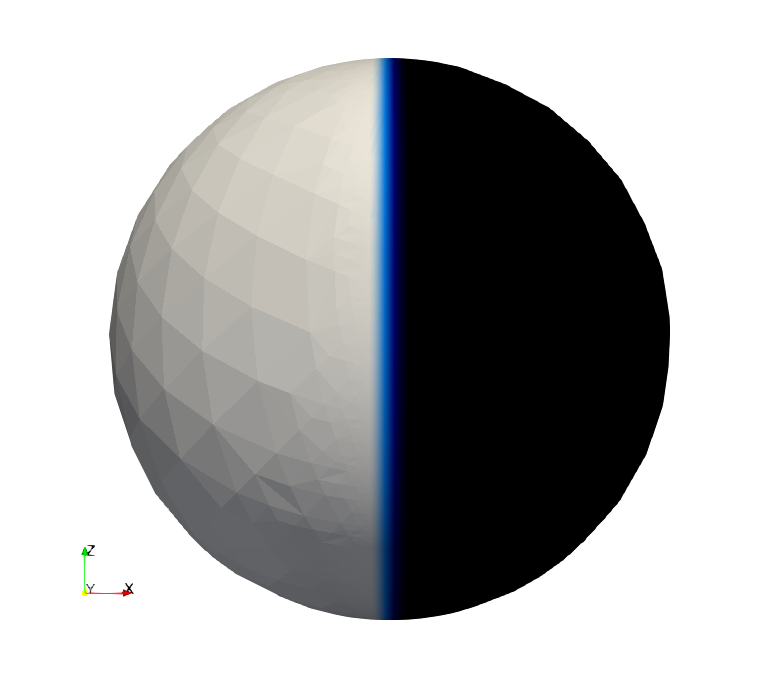}
\includegraphics[angle=-0,width=0.3\textwidth]{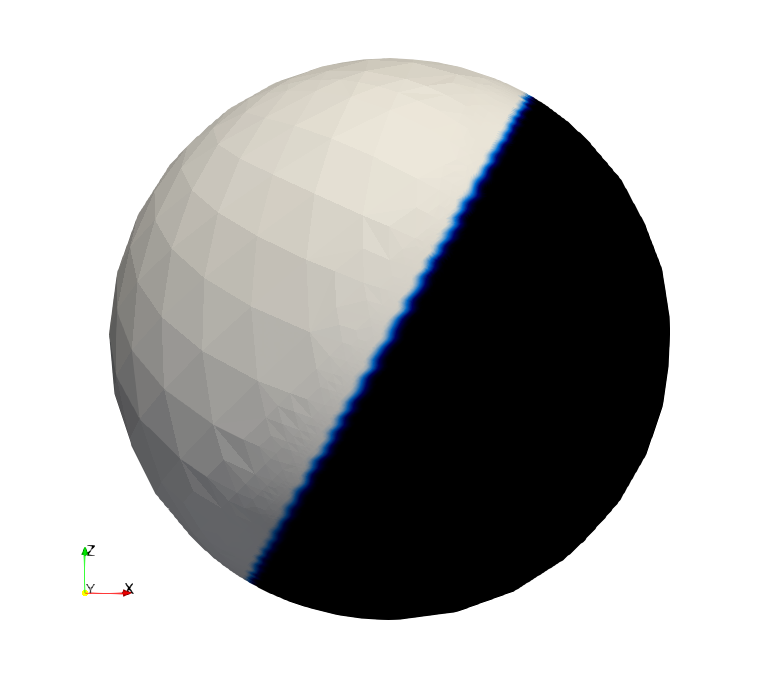}
\includegraphics[angle=-0,width=0.3\textwidth]{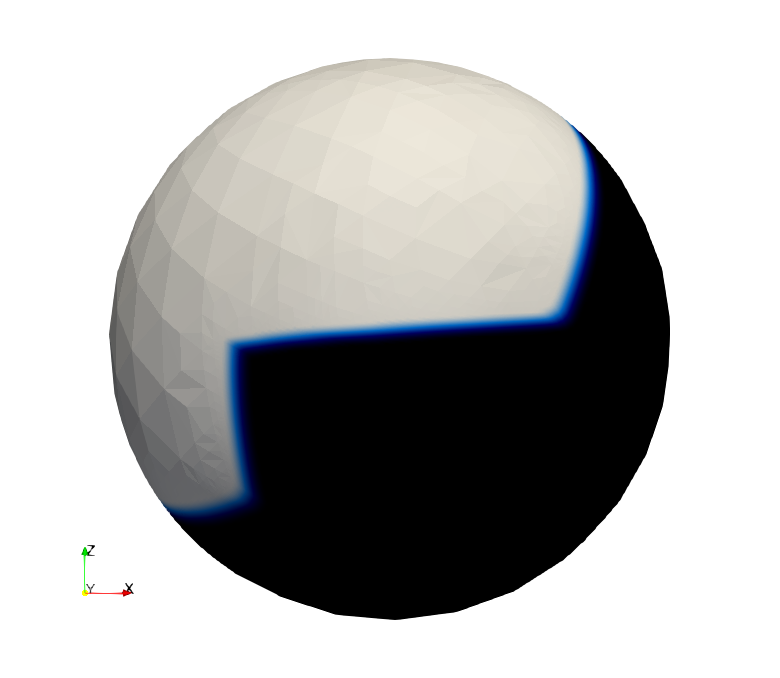}
\includegraphics[angle=-0,width=0.3\textwidth]{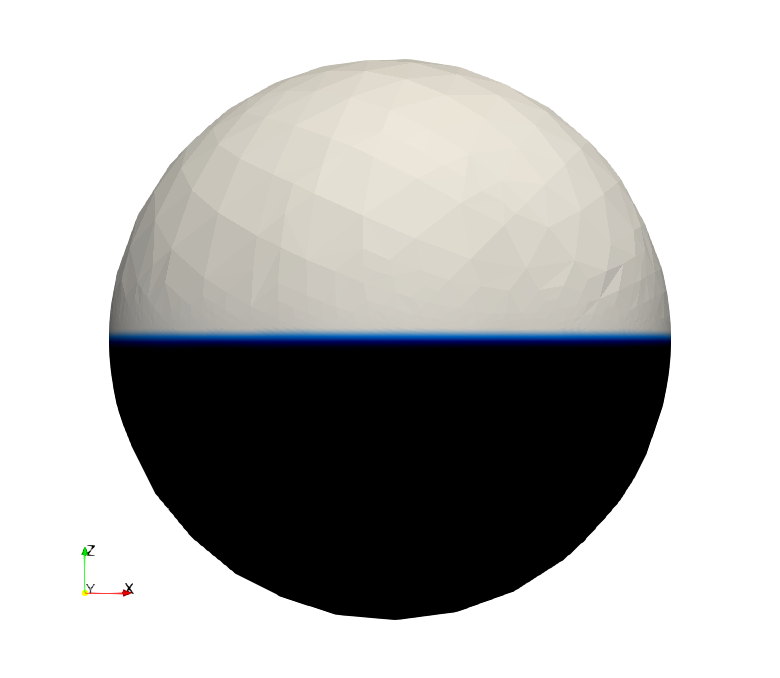}
\includegraphics[angle=-0,width=0.3\textwidth]{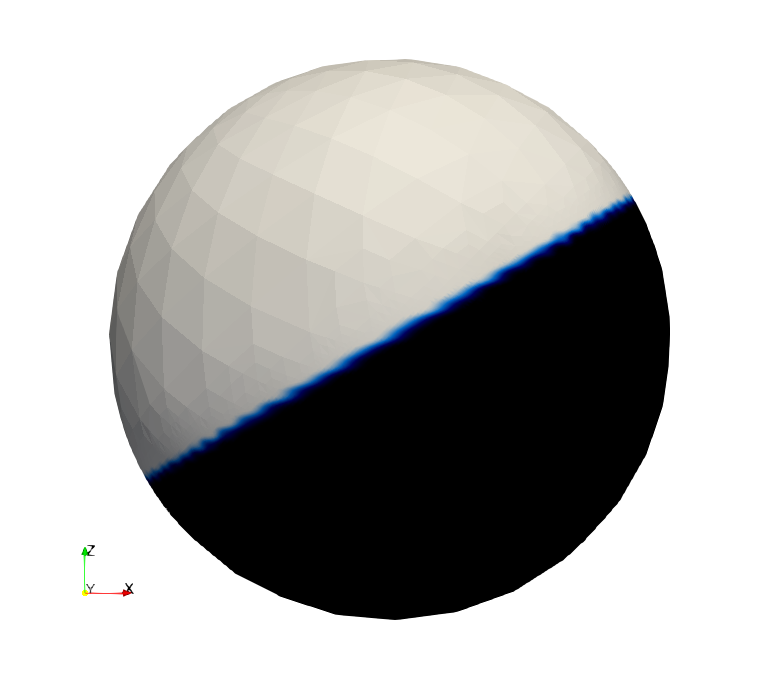}
\includegraphics[angle=-0,width=0.3\textwidth]{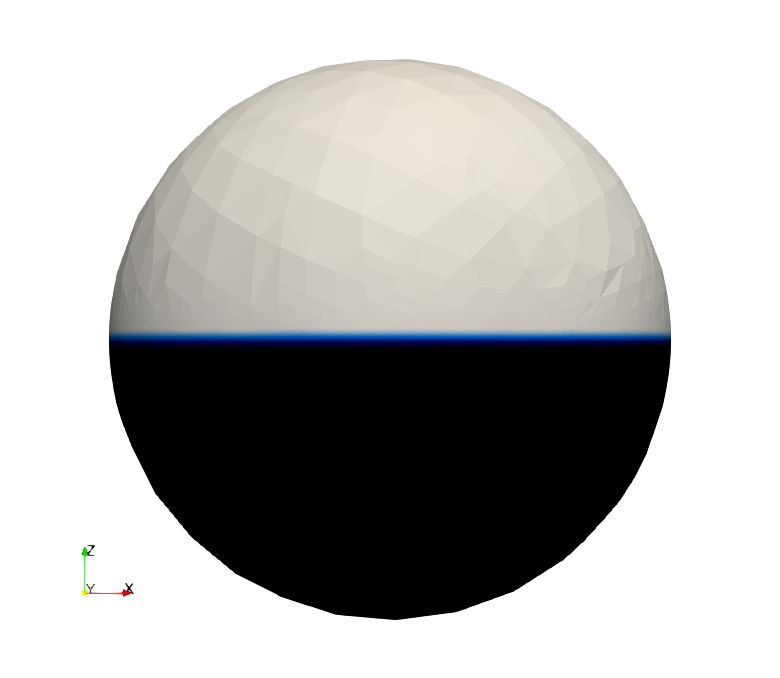}
\includegraphics[angle=-0,width=0.3\textwidth]{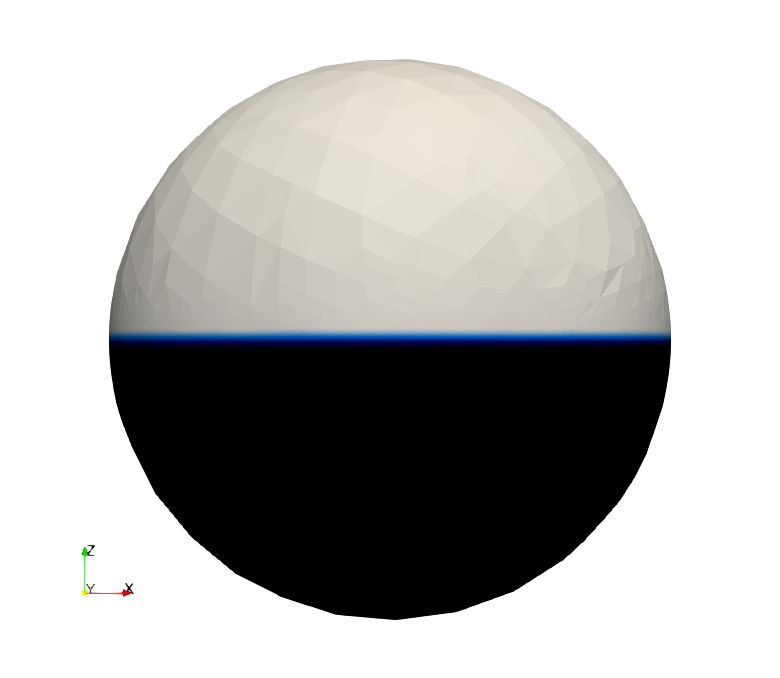}
\caption{($\epsilon=(16\pi)^{-1}$) 
Anisotropic Cahn--Hilliard equation on a sphere.
Snapshots of the evolutions at times $t=0, 1, 5$ (from left to right).
The initial data approximates a great circle within a hyperplane that
is tilted by $10^\circ, 30^\circ, 60^\circ$ (from top to bottom).
}
\label{fig:lem21tilts}

\end{figure}

We use the anisotropy \eqref{eq:L4abc}(b) with $\delta=10^{-2}$ to numerically
investigate the statement from Lemma~\ref{lem:greatcircle}. To this end, we
choose as initial interface a perturbed great circle on the unit sphere, that
lies in a hyperplane that makes a certain fixed angle $\theta$ with the 
$x-z$-plane. In particular, in Figure~\ref{fig:lem21tilts} we show simulations
for $\theta = 10^\circ, 30^\circ, 60^\circ$. 
Here we let 
$\epsilon=(16\pi)^{-1}$, $N_c = 1$, $N_f = 128$ and $\tau = 10^{-4}$.
We note that each simulation settles on a great circle on the sphere. Here we
note that as $\nu_{\min} = e_3$, the equator has the lowest energy of all the
great circles, and so if the initial curve on $\mathcal M$ is sufficiently
inclined, that global minimizer is indeed reached by the evolution. For only
small deviations from a north-south great circle, the evolutions 
settle on that local minimizer instead.
We remark that the discrete steady state solutions exhibit final
discrete energies $\mathcal{F}^h(W^n, \Phi^n)$ of 12.7, 6.7 and 6.7,
respectively. This once again confirms that the equator has the lowest energy.

\subsubsection{Convergence experiment on the unit sphere}

We use the rotationally symmetric solution for the Mullins--Sekerka problem on
the unit sphere from \cite{Ratz16} for a convergence experiment for
our approximations as $\epsilon\to0$, see also \cite{dendritic}. 

We consider an annulus domain on the unit sphere, which is bounded by two
circles on the lower half with radii $1 > R_1(t) > R_2(t) > 0$. This uniquely
defines angles $\theta_i(t) \in (\frac\pi2, \pi)$ such that 
$R_i(t) = \sin\theta_i(t)$ and $h_i(t) = \cos\theta_i(t) <0 $ denotes the
heights of the two circles.
That means at any given time, the annulus will enclose a surface area of
$2\pi(h_1(t) - h_2(t)) = 2\pi(\cos\theta_1(t) - \cos\theta_2(t))$.
Hence $a_0 := \cos\theta_1(0) - \cos\theta_2(0) = 
\cos\theta_1(t) - \cos\theta_2(t)$ is an invariant, and so
\begin{equation} \label{eq:theta2}
\theta_2(t) = \arccos(\cos\theta_1(t) - a_0),
\end{equation}
which means that the ODE system for $(\theta_1,\theta_2)$ from \cite{Ratz16} 
can be reduced to a scalar differential equation. The former is given by
\begin{subequations} \label{eq:ODE}
\begin{equation} \label{eq:ODEa}
\lambda \theta_i'(t) = - \frac{c_+(t)}{\sin\theta_i(t)},\ i=1,2,
\end{equation}
where
\begin{equation} \label{eq:ODEb}
c_+(t) = - \alpha \frac{\cot\theta_1(t) + \cot\theta_2(t)}{
\ln ( \tan\frac{\theta_1(t)}2) - \ln(\tan\frac{\theta_2(t)}2)}.
\end{equation}
\end{subequations}
Combining \eqref{eq:ODE} with \eqref{eq:theta2} yields a scalar differential 
equation for $\theta_1$.

\begin{figure}
\center
\includegraphics[angle=-90,width=0.45\textwidth]{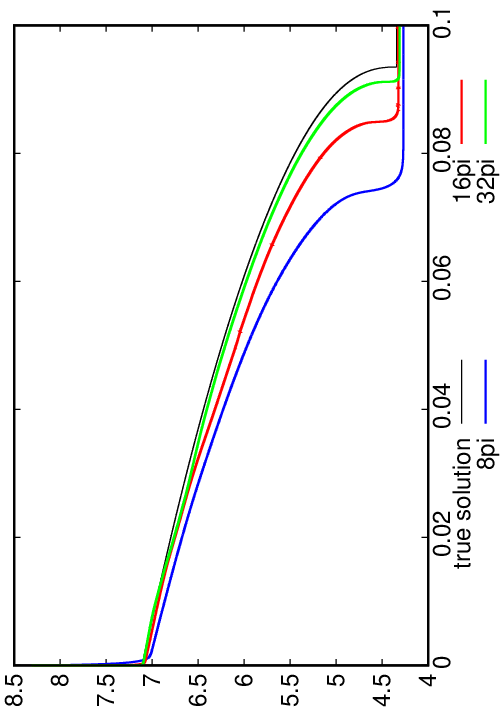}
\includegraphics[angle=-90,width=0.45\textwidth]{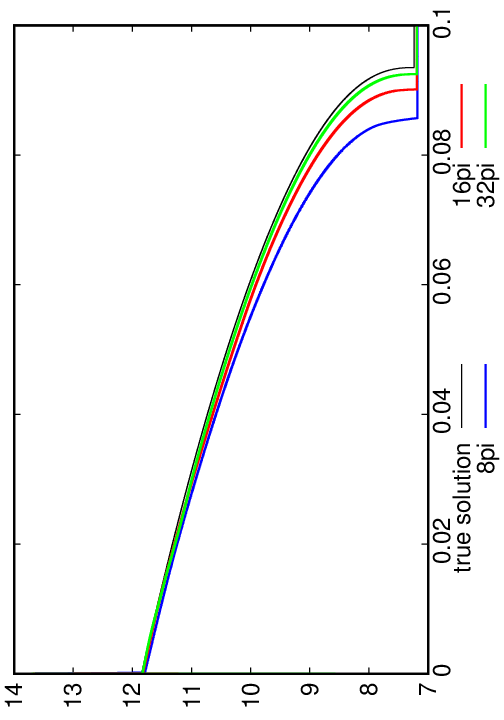}
\caption{($\alpha = \frac{\sqrt{2}}3$, $\lambda=2$)
Energy plots for \eqref{eq:qfea} (left) and \eqref{eq:fea} (right) 
for $\epsilon=(8\pi)^{-1}$,
$\epsilon=(16\pi)^{-1}$ and $\epsilon=(32\pi)^{-1}$, compared to the sharp
interface solution.
}
\label{fig:annulus_ratz3}
\end{figure}

For the initial data we choose $R_1(0) = 2 R_2(0) = 0.8$, so that 
$h_1(0) = -0.6$ and $h_2(0) = -\sqrt{0.84}$, respectively.
In order to get the same time scale as in \cite{Ratz16}, 
we use $\alpha = \frac{\sqrt{2}}3$ and $\lambda=2$ for these experiments,
and we visualize the phase field energies
\begin{equation*} 
\tfrac12\epsilon|\nabs\Phi^m|_h^2 + 
\epsilon^{-1} \langle\Psi(\Phi^m), 1\rangle^h_{\mathcal{M}^h} \,,
\end{equation*}
for deceasing values of $\epsilon$, 
compared to the energy $2\pi\cPsi (\sin\theta_1(t_m) +
\sin\theta_2(t_m))$ of the sharp interface solution. See
Figure~\ref{fig:annulus_ratz3}, where our numerical results confirm the
asymptotic convergence as $\epsilon\to0$ to the sharp interface limit.

\subsubsection{Spinodal decomposition on closed surfaces}

\begin{figure}
\center
\mbox{
\includegraphics[angle=-0,width=0.15\textwidth]{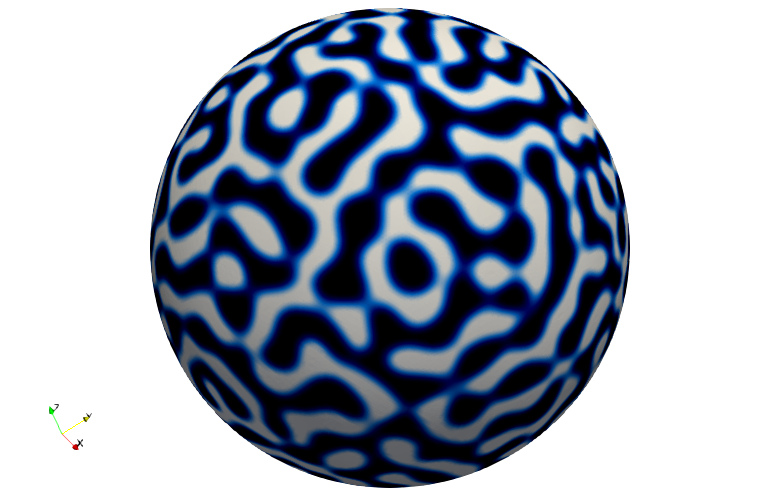}
\includegraphics[angle=-0,width=0.15\textwidth]{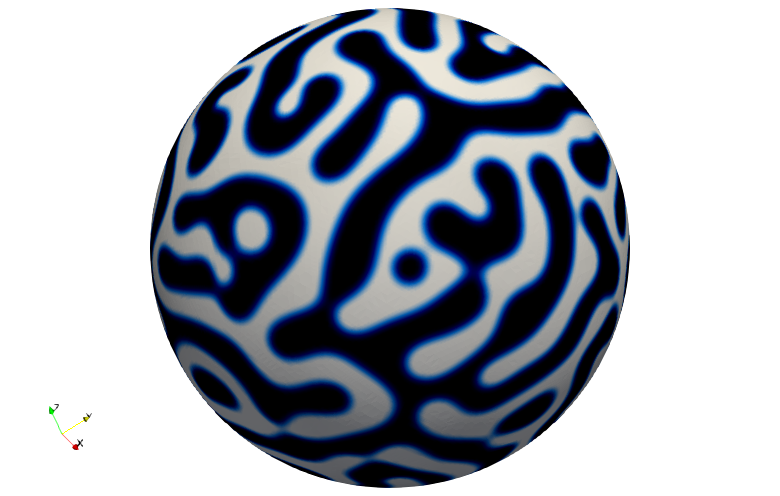}
\includegraphics[angle=-0,width=0.15\textwidth]{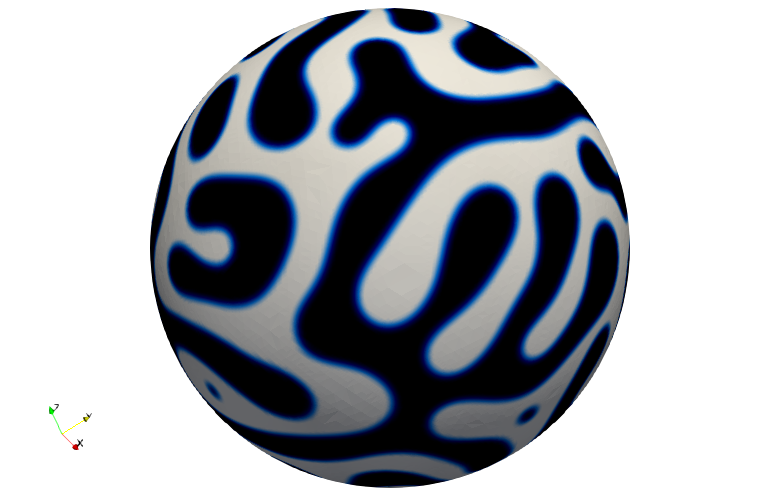}
\includegraphics[angle=-0,width=0.15\textwidth]{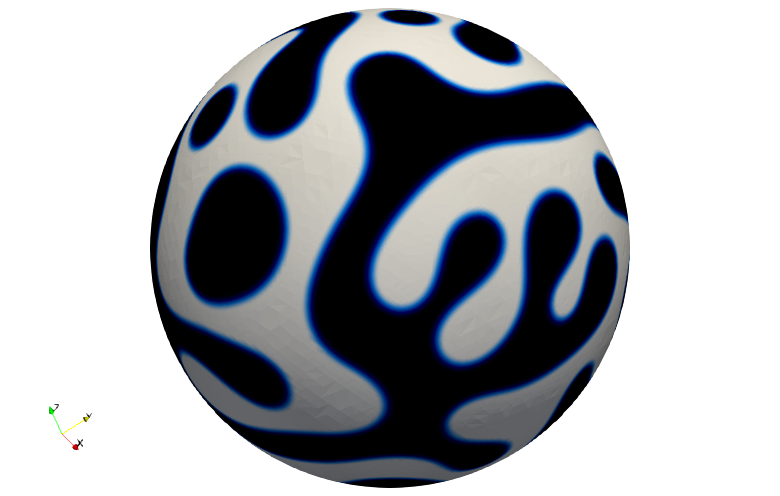}
\includegraphics[angle=-0,width=0.15\textwidth]{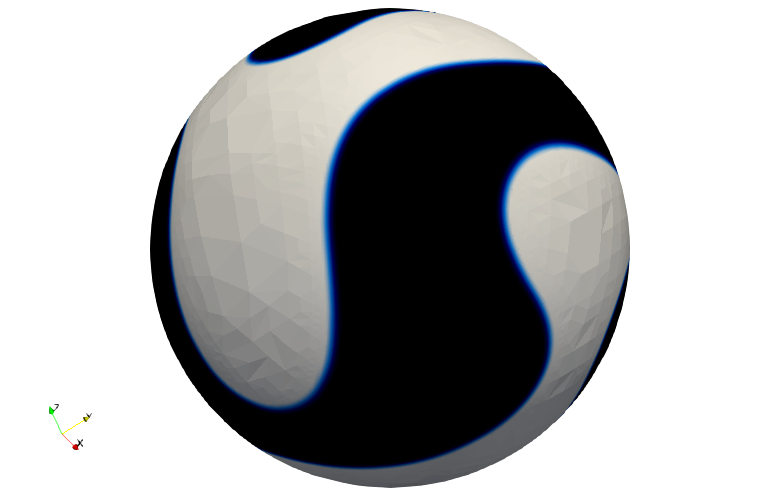}
\includegraphics[angle=-0,width=0.15\textwidth]{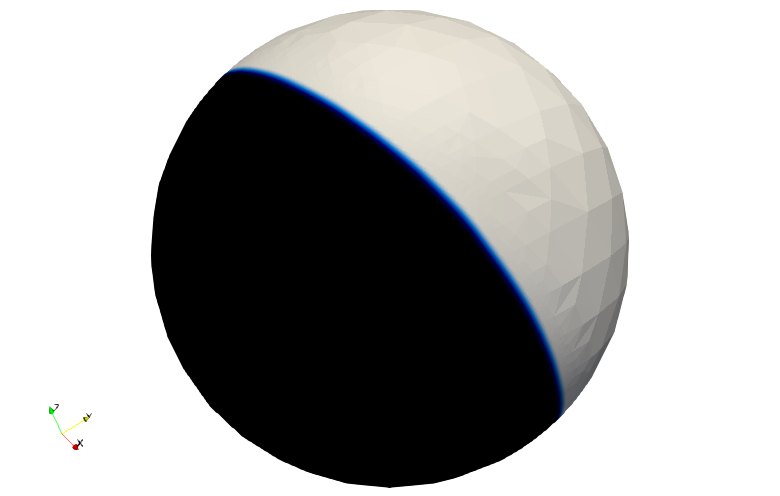}
}
\caption{($\epsilon=(16\pi)^{-1}$) 
Spinodal decomposition on a sphere, for $\gamma(p) = |p|$.
Snapshots of the evolutions at times $t=10^{-4}, 2\times 10^{-4}, 5\times10^{-4}, 10^{-3}, 0.01, 1$.
}
\label{fig:spinodal_iso}
\end{figure}

\begin{figure}
\center
\mbox{
\includegraphics[angle=-0,width=0.15\textwidth]{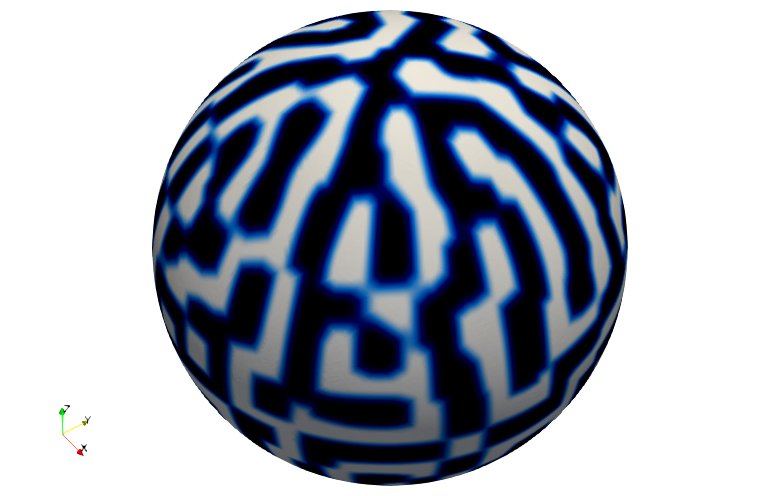}
\includegraphics[angle=-0,width=0.15\textwidth]{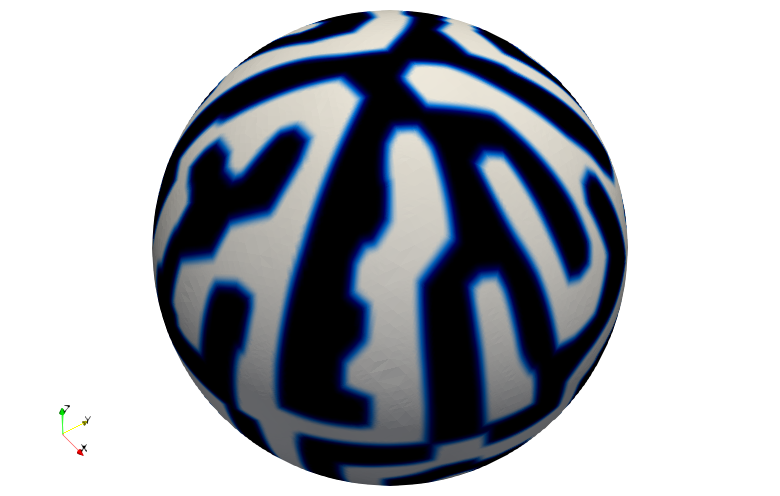}
\includegraphics[angle=-0,width=0.15\textwidth]{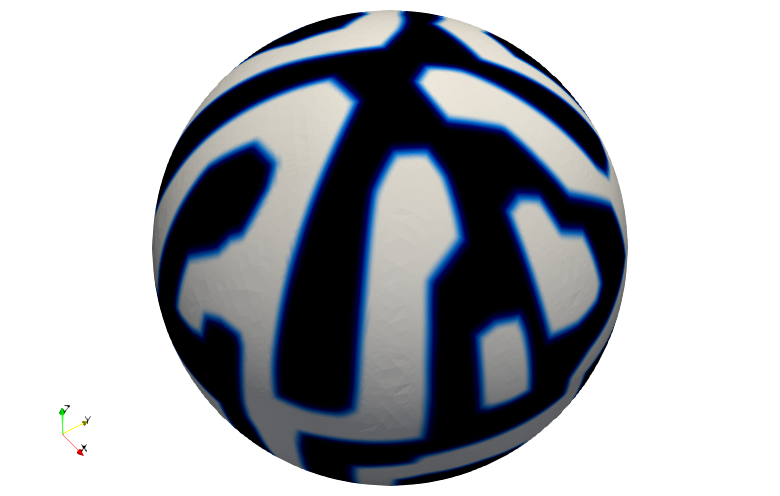}
\includegraphics[angle=-0,width=0.15\textwidth]{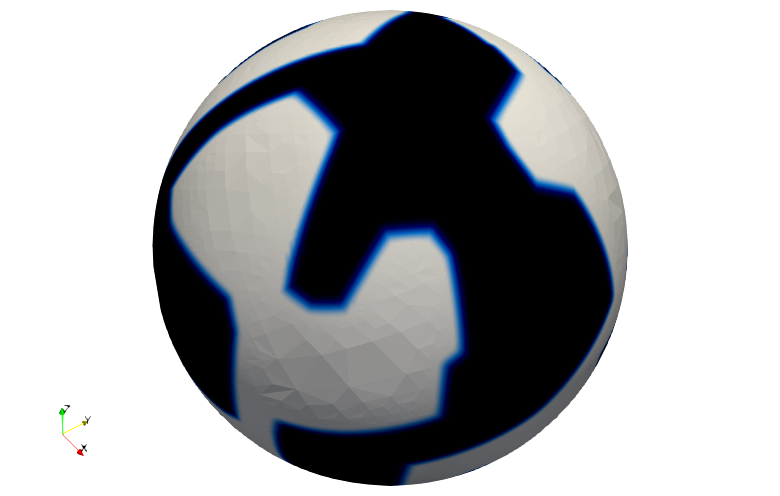}
\includegraphics[angle=-0,width=0.15\textwidth]{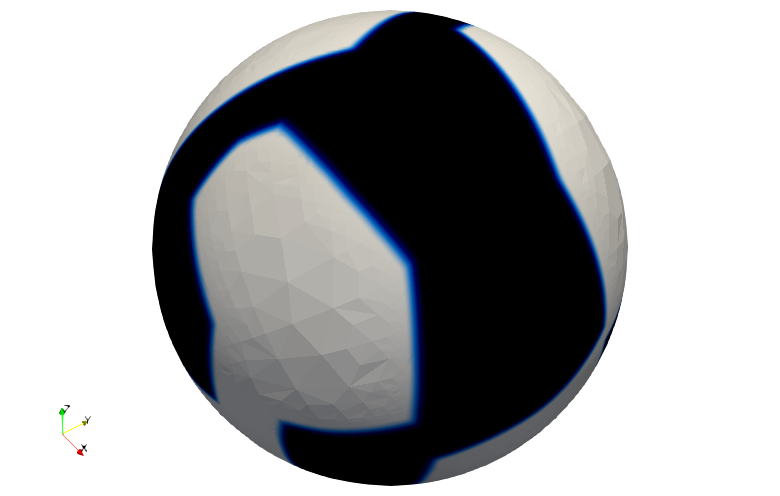}
\includegraphics[angle=-0,width=0.15\textwidth]{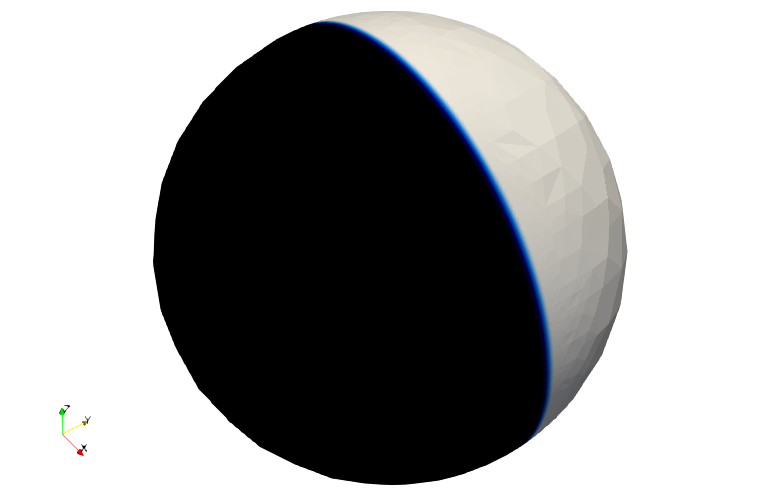}
}
\caption{($\epsilon=(16\pi)^{-1}$) 
Spinodal decomposition on a sphere, for \eqref{eq:L4abc}(a).
Snapshots of the evolutions at times $t=2\times 10^{-4}, 5\times10^{-4}, 
10^{-3}, 5\times10^{-3}, 0.01, 2$.
}
\label{fig:spinodal_L44}
\end{figure}

\begin{figure}
\center
\mbox{
\includegraphics[angle=-0,width=0.15\textwidth]{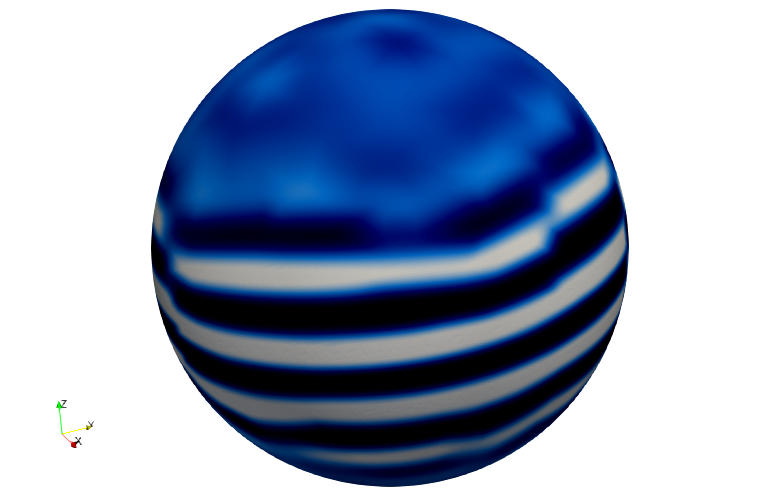}
\includegraphics[angle=-0,width=0.15\textwidth]{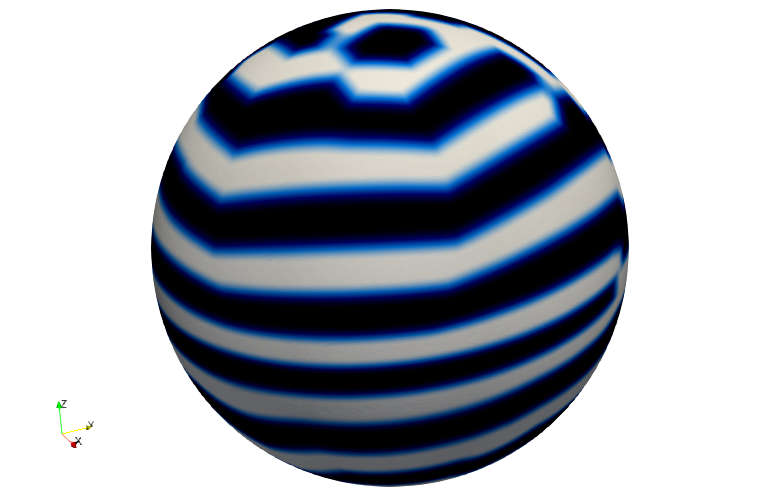}
\includegraphics[angle=-0,width=0.15\textwidth]{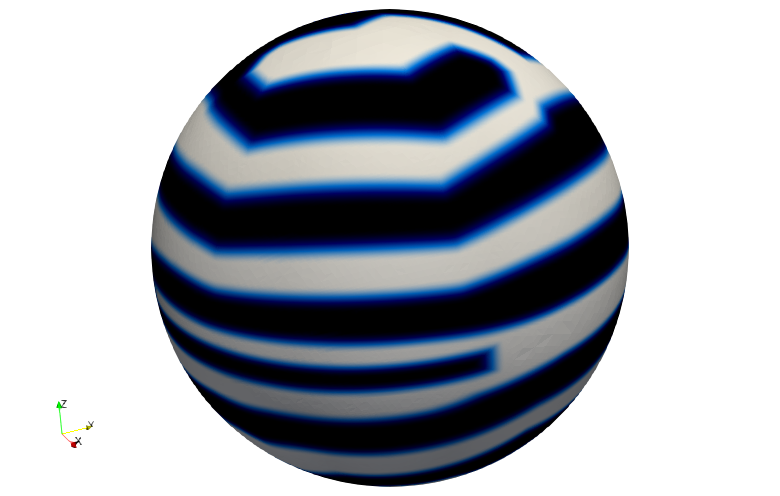}
\includegraphics[angle=-0,width=0.15\textwidth]{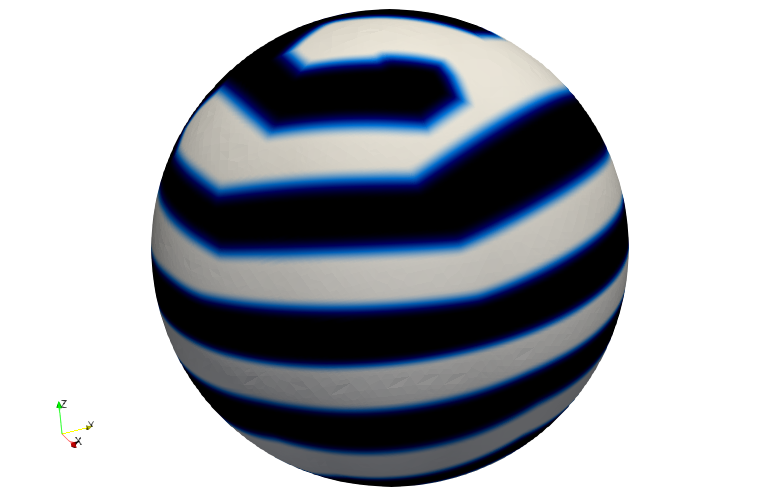}
\includegraphics[angle=-0,width=0.15\textwidth]{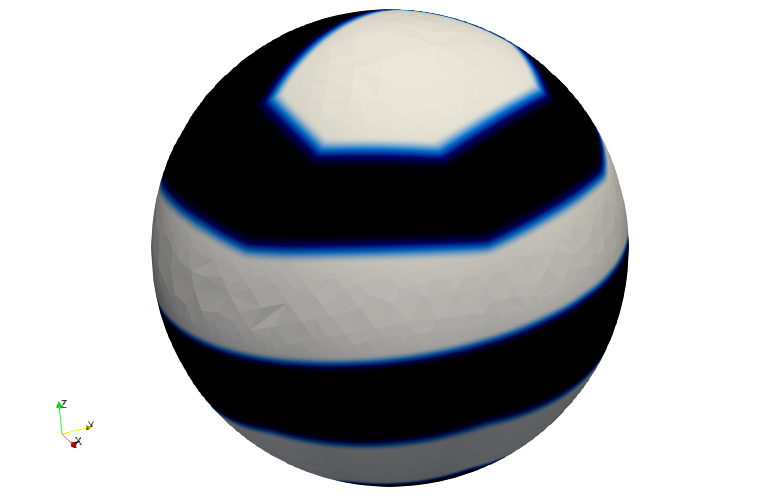}
\includegraphics[angle=-0,width=0.15\textwidth]{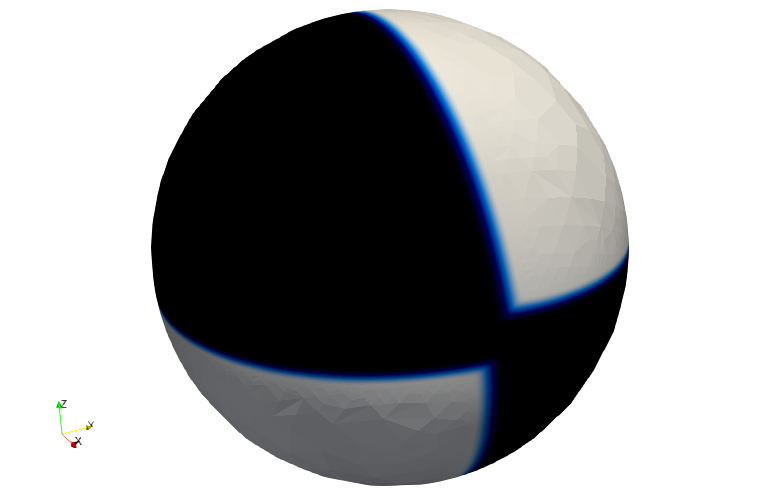}
}
\caption{($\epsilon=(16\pi)^{-1}$) 
Spinodal decomposition on a sphere, for \eqref{eq:L4abc}(b).
Snapshots of the evolutions at times $t=2\times 10^{-4}, 5\times10^{-4}, 
10^{-3}, 2\times10^{-3}, 0.01, 1$.
}
\label{fig:spinodal_L4}
\end{figure}

\begin{figure}
\center
\mbox{
\includegraphics[angle=-0,width=0.15\textwidth]{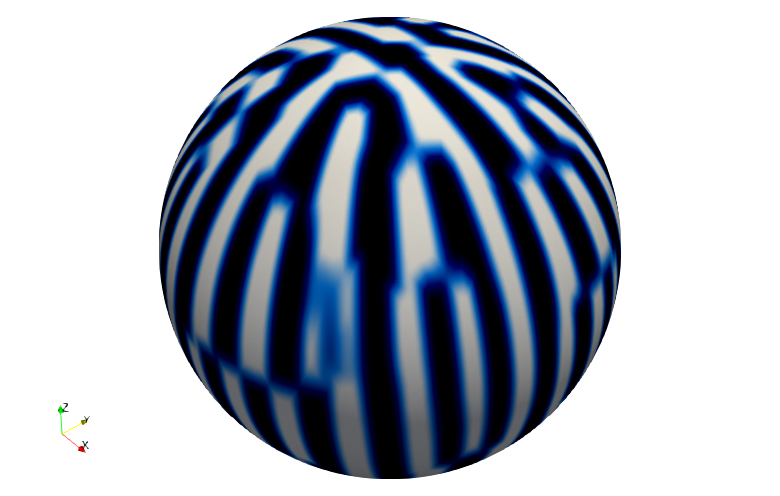}
\includegraphics[angle=-0,width=0.15\textwidth]{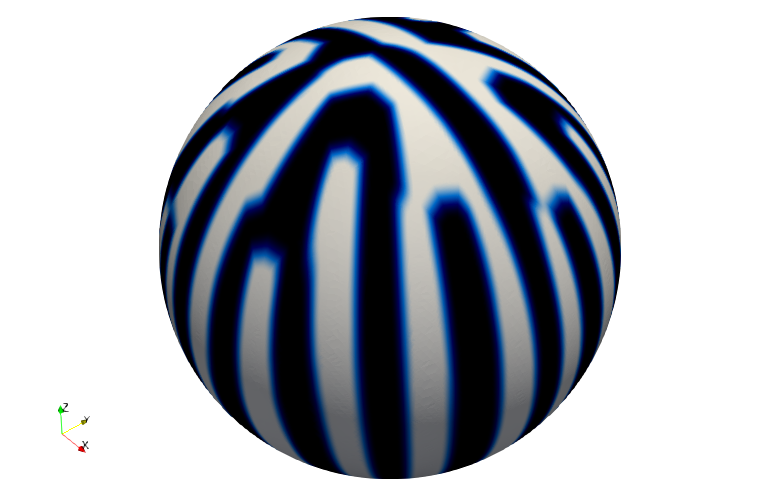}
\includegraphics[angle=-0,width=0.15\textwidth]{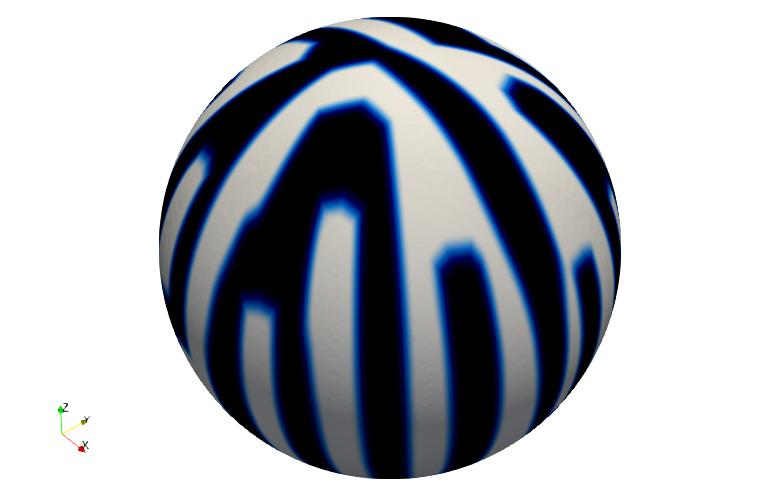}
\includegraphics[angle=-0,width=0.15\textwidth]{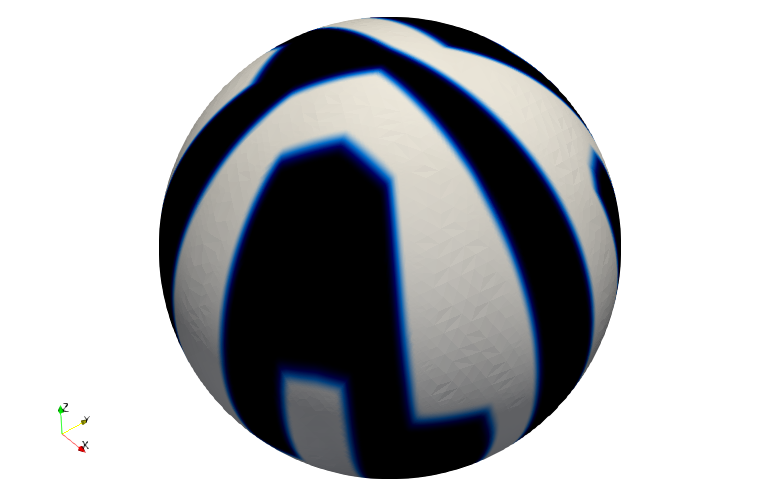}
\includegraphics[angle=-0,width=0.15\textwidth]{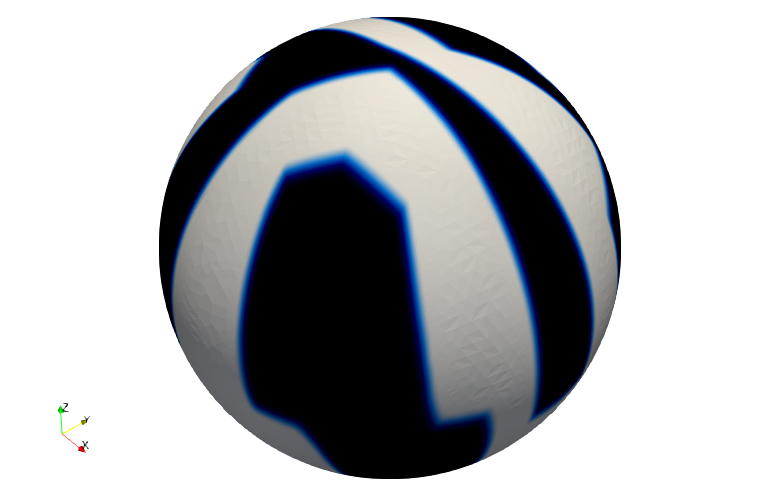}
\includegraphics[angle=-0,width=0.15\textwidth]{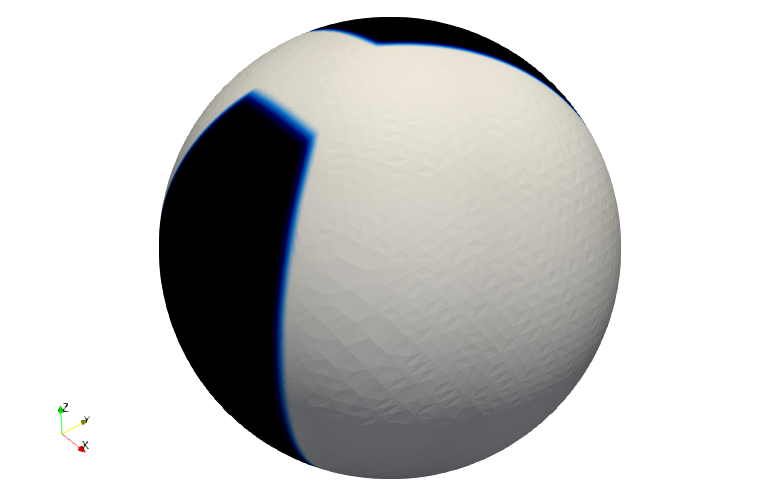}
}
\caption{($\epsilon=(16\pi)^{-1}$) 
Spinodal decomposition on a sphere, for \eqref{eq:L4abc}(c).
Snapshots of the evolutions at times $t=2\times 10^{-4}, 5\times10^{-4}, 
10^{-3}, 5\times10^{-3}, 0.01, 1$.
}
\label{fig:spinodal_L44fac}
\end{figure}

In this section we consider 
$\epsilon=(16\pi)^{-1}$, as well as $N_c = 1$, $N_f = 128$, $\tau=10^{-6}$.
We use different surface energy densities for evolutions that start from a
random mixture with mean zero and values in $[-0.1,0.1]$.
At first we choose the unit sphere $\mathcal M = \bS^2$.
In this isotropic case, we observe the well known spinodal decomposition 
patterns, see Figure~\ref{fig:spinodal_iso}.
When we use the anisotropy \eqref{eq:L4abc}(a), for which in the limit 
$\delta=0$ the normals $\nu=\pm e_3$ have the same energy density as the 
other six main facet normals of the crystalline Wulff shape,
during the spinodal decomposition we can observe corners, see 
Figure~\ref{fig:spinodal_L44}. Note that the spatially homogeneous anisotropy
leads to hexagonal symmetries only at the two poles, whereas the patterns near
the equator resemble squares.
The corresponding evolution for the anisotropy \eqref{eq:L4abc}(b)
is shown in Figure~\ref{fig:spinodal_L4}. Here it can be seen 
that $\nu=\pm e_3$ is now the preferred normal direction, and so the 
interfaces very quickly align with it.
As a final experiment, we also use the anisotropic density
\eqref{eq:L4abc}(c), see Figure~\ref{fig:spinodal_L44fac} for the numerical 
results. Now the opposite effect can be observed: since the normals 
$\nu=\pm e_3$ are relatively expensive, they are avoided by the developing
interfaces.
We also observe that out of the four spinodal decomposition experiments, only
two settle on a global energy minimizer as described in 
Lemma~\ref{lem:greatcircle}. We conjecture that the final shape for the other
two simulations is made up of ``straight'' segments that are aligned with the
Wulff shape, in the sense that their normal vectors correspond to directions
with minimal energy density. The anisotropic nature of the surface energy then
seems to make it impossible to go from this local minimizer to one of the
global ones.

\begin{figure}
\center
\mbox{
\includegraphics[angle=-0,width=0.15\textwidth]{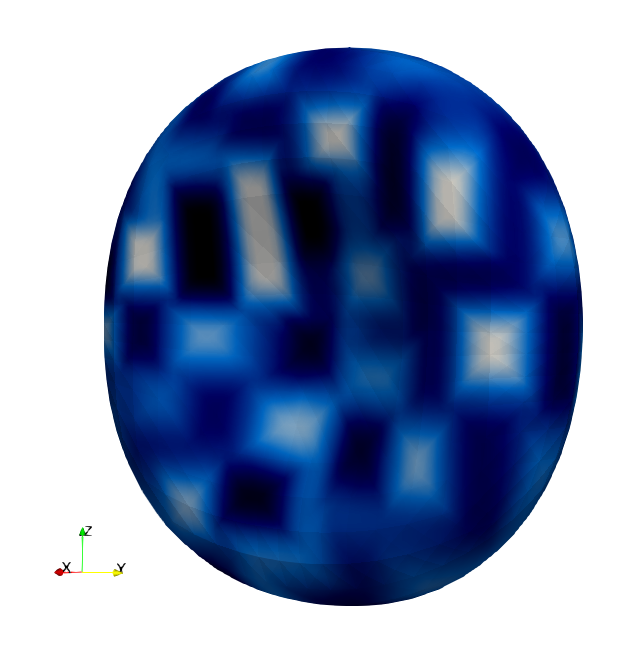}
\includegraphics[angle=-0,width=0.15\textwidth]{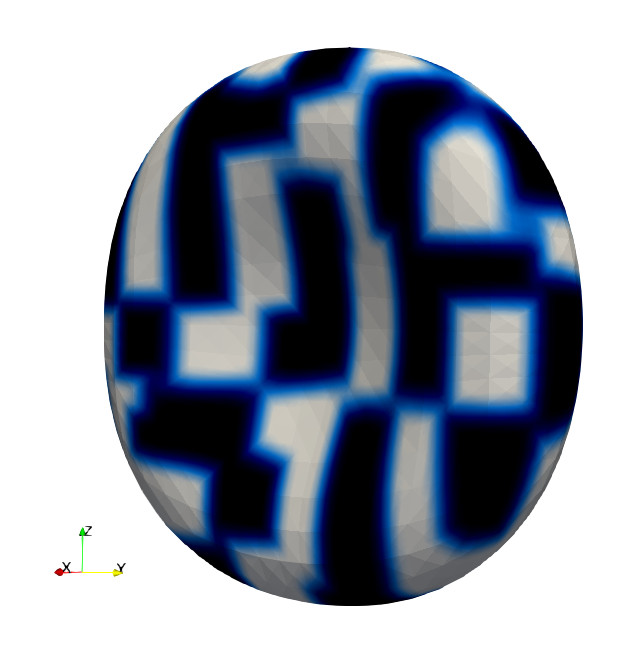}
\includegraphics[angle=-0,width=0.15\textwidth]{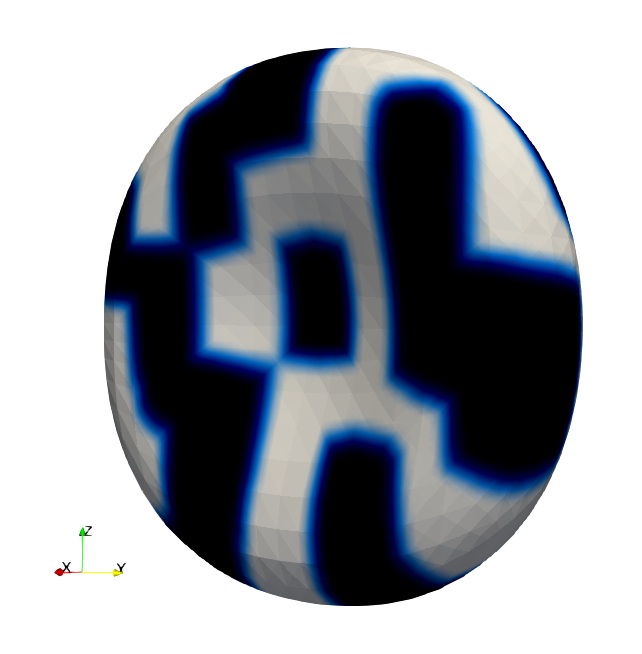}
\includegraphics[angle=-0,width=0.15\textwidth]{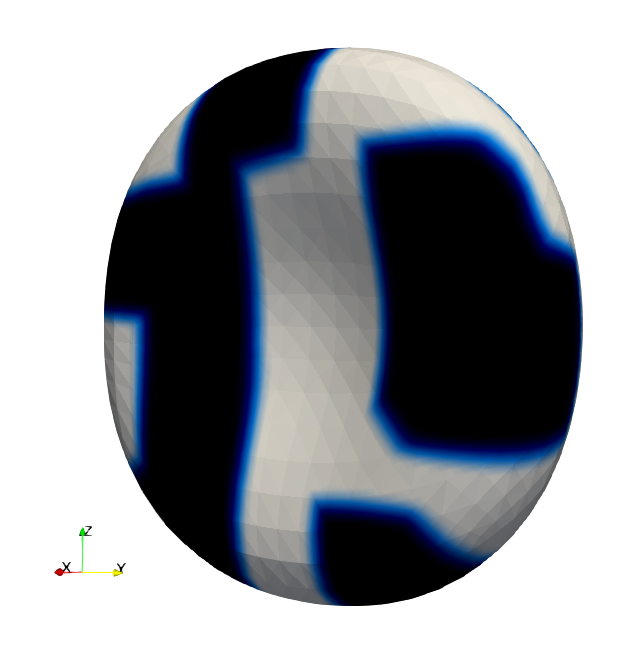}
\includegraphics[angle=-0,width=0.15\textwidth]{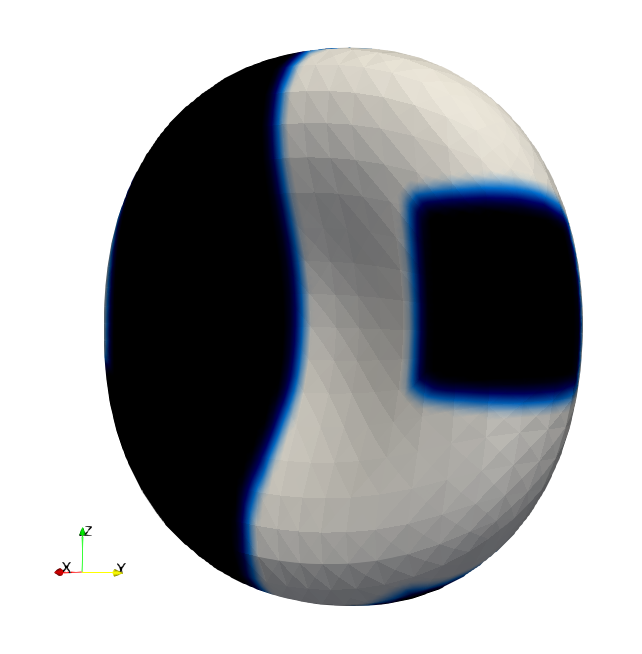}
\includegraphics[angle=-0,width=0.15\textwidth]{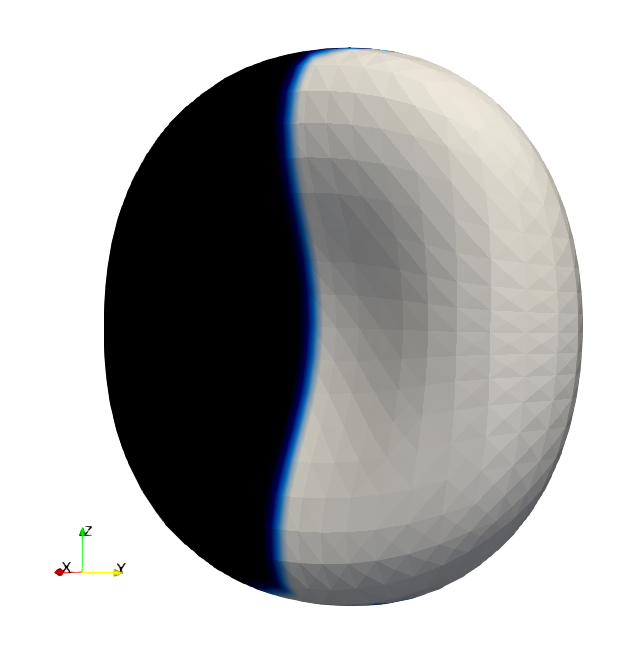}
}
\caption{($\epsilon=(16\pi)^{-1}$) 
Spinodal decomposition on a biconcave disk, for \eqref{eq:L4abc}(a).
Snapshots of the evolutions at times 
$t=10^{-4},2\times 10^{-4}, 5\times10^{-4}, 10^{-3}, 5\times10^{-3}, 0.01$.
}
\label{fig:bloodcell_L44}
\end{figure}

We end this subsection with a simulation on a surface that is the 
boundary of a nonconvex domain. In particular, we base $\mathcal{M}^h$ on a
rescaled version of the biconcave disk obtained in \cite[Fig.\ 8]{nsns}. 
The total
dimensions of $\mathcal{M}^h$ are about $0.28 \times 1 \times 1$. On this
surface we repeat the simulation from Figure~\ref{fig:spinodal_L44}. The results
of this new numerical experiment are shown in Figure~\ref{fig:bloodcell_L44}, 
where we once again observe the squared patterns due to the alignment of the 
surface with respect to the Wulff shape of the anisotropy \eqref{eq:L4abc}(a).

\subsubsection{Crystal growth on a sphere cap} \label{sec:cap}

In this section we use the hexagonal, spatially homogeneous anisotropic density
\eqref{eq:L4abc}(a) for computations of crystal growth on
a sphere cap. Here the sphere cap $\mathcal M \subset \bS^2$ is given by
\begin{align*}
& \mathcal M = 
\left\{ \begin{pmatrix} 0 \\ 0 \\ 1 \end{pmatrix} \right\} \cup
\left\{ \begin{pmatrix} \tfrac{z_1}{|z|}\,\cos\theta(z) \\
\tfrac{z_2}{|z|}\,\cos\theta(z) \\ 
\sin\theta(z)
\end{pmatrix}:
z = \begin{pmatrix} z_1 \\ z_2 \\ 0 \end{pmatrix} \in \bR^3, 0 < |z| \leq 1
\right\}, \\ & \text{where} \
\theta(z)=(\tfrac\pi2+\tfrac{\pi}{18})(1-|z|) - \tfrac{\pi}{18}.
\end{align*}
For the physical parameters, similarly to \cite[Fig.\ 9]{vch},
we choose $\alpha=0.03$, $\rho = 0.01$ and $\uD = -8$.

\begin{figure}
\center
\includegraphics[angle=-0,width=0.22\textwidth]{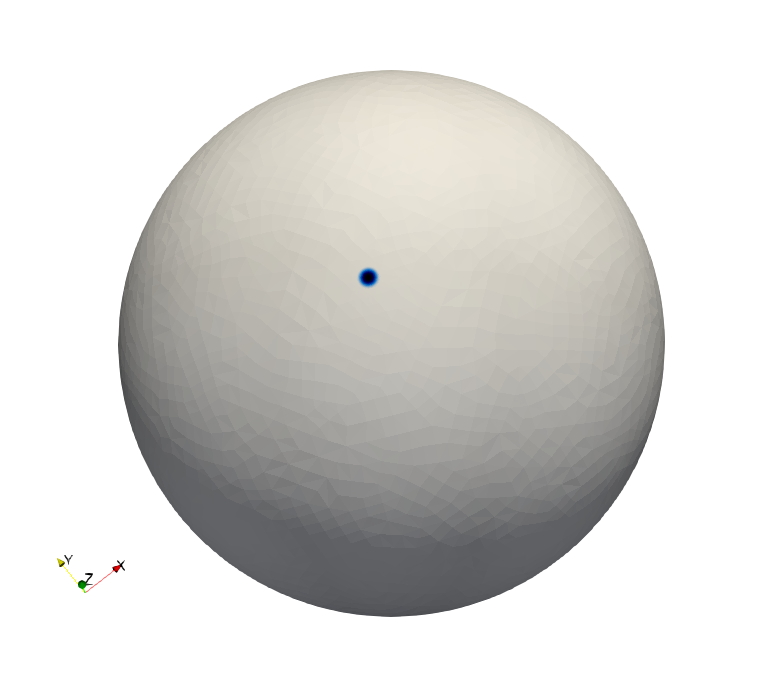}
\includegraphics[angle=-0,width=0.22\textwidth]{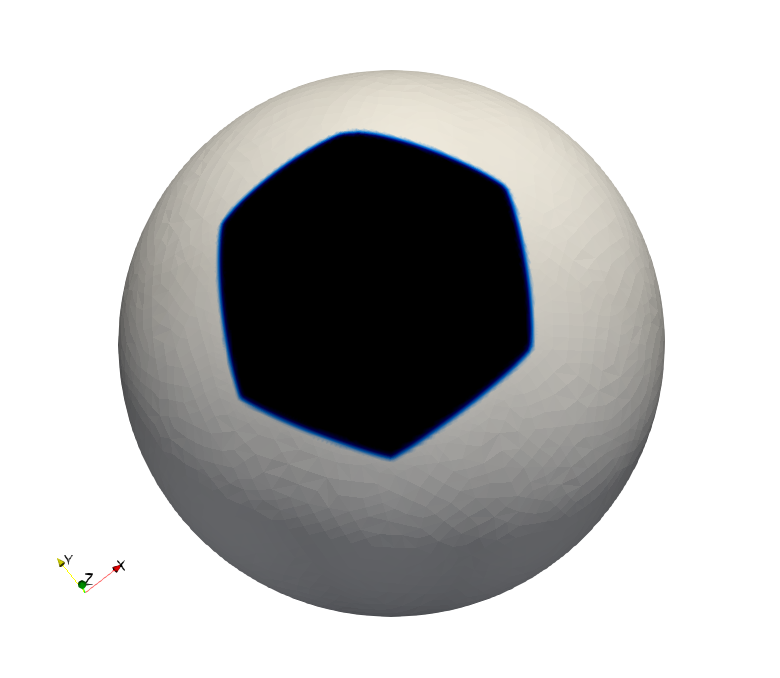}
\includegraphics[angle=-0,width=0.22\textwidth]{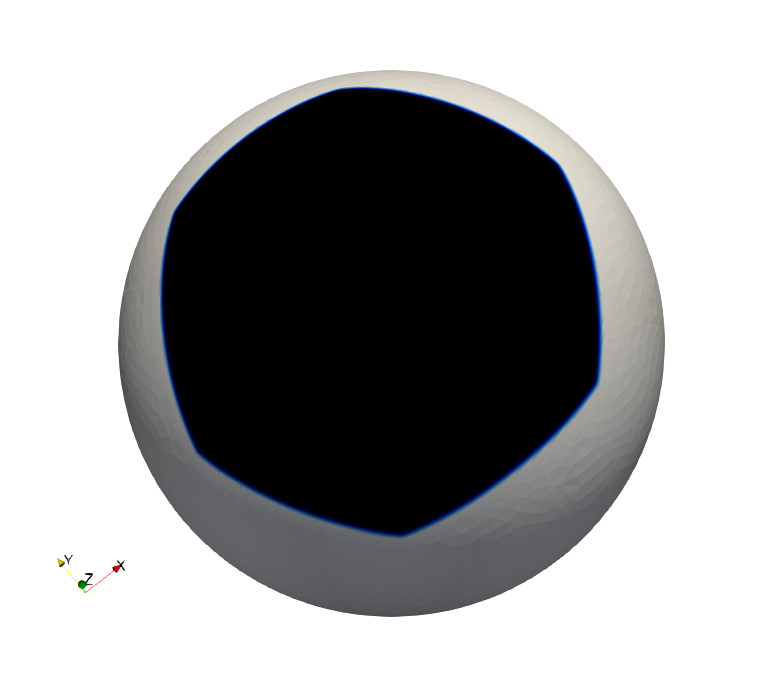}
\\
\includegraphics[angle=-0,width=0.22\textwidth]{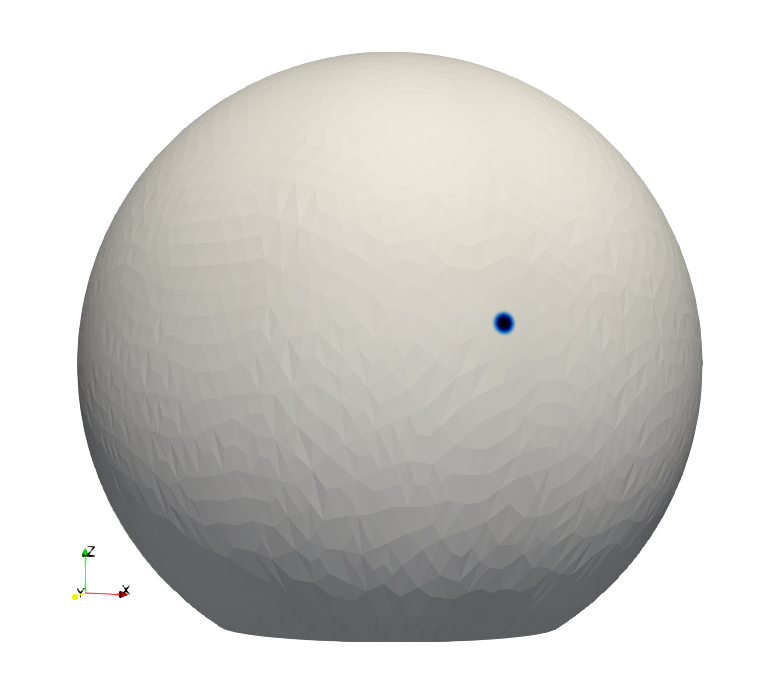}
\includegraphics[angle=-0,width=0.22\textwidth]{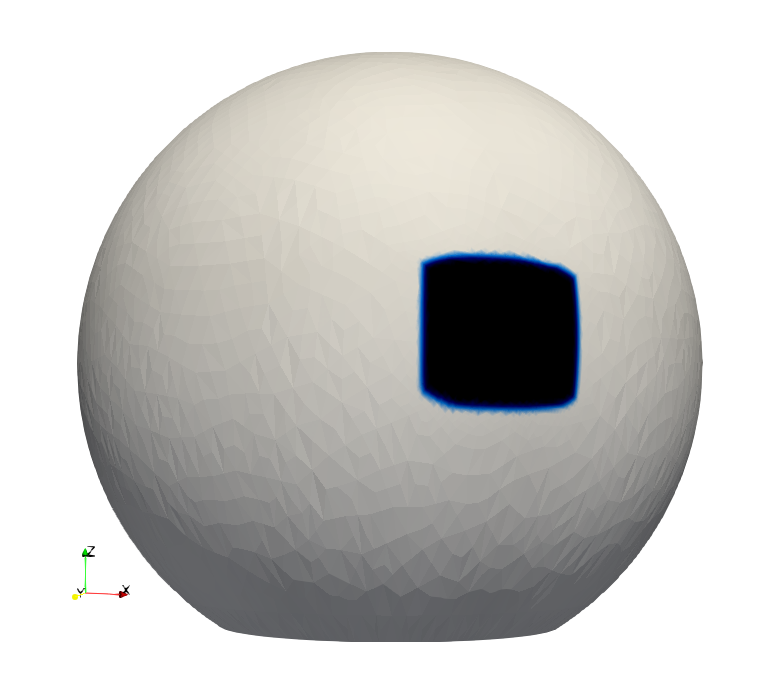}
\includegraphics[angle=-0,width=0.22\textwidth]{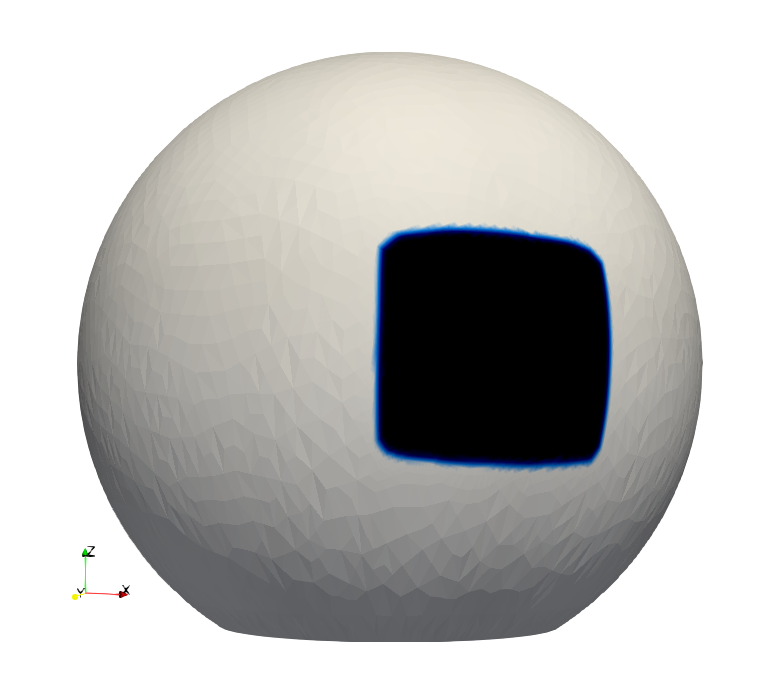}
\\
\includegraphics[angle=-0,width=0.22\textwidth]{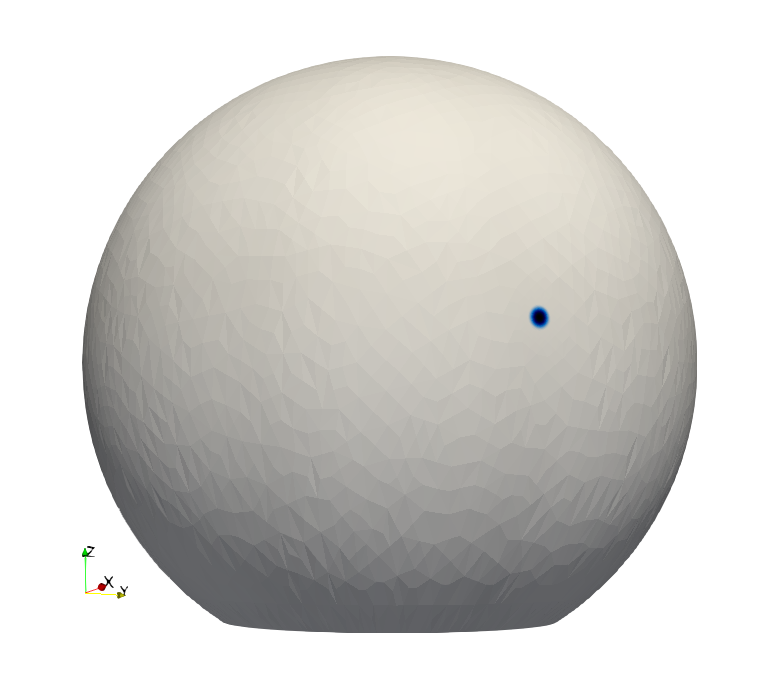}
\includegraphics[angle=-0,width=0.22\textwidth]{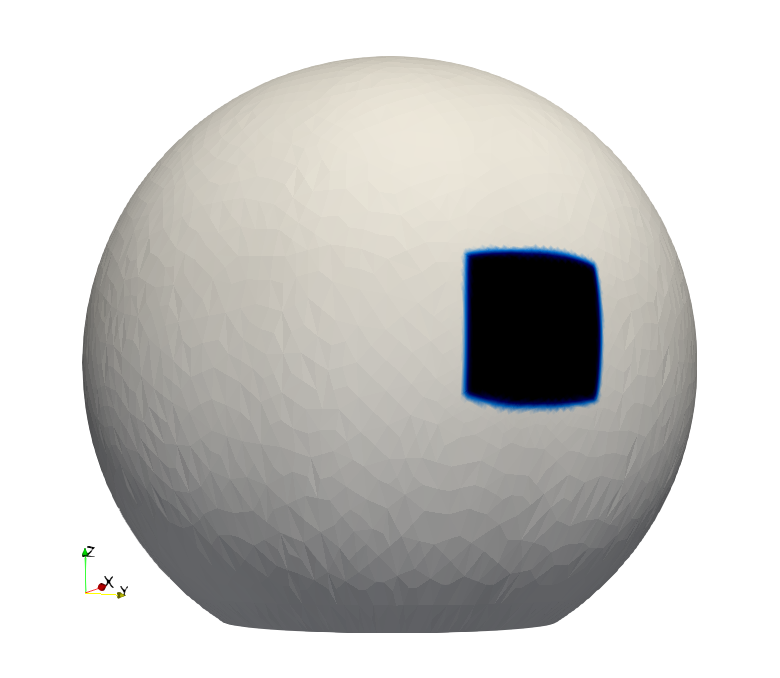}
\includegraphics[angle=-0,width=0.22\textwidth]{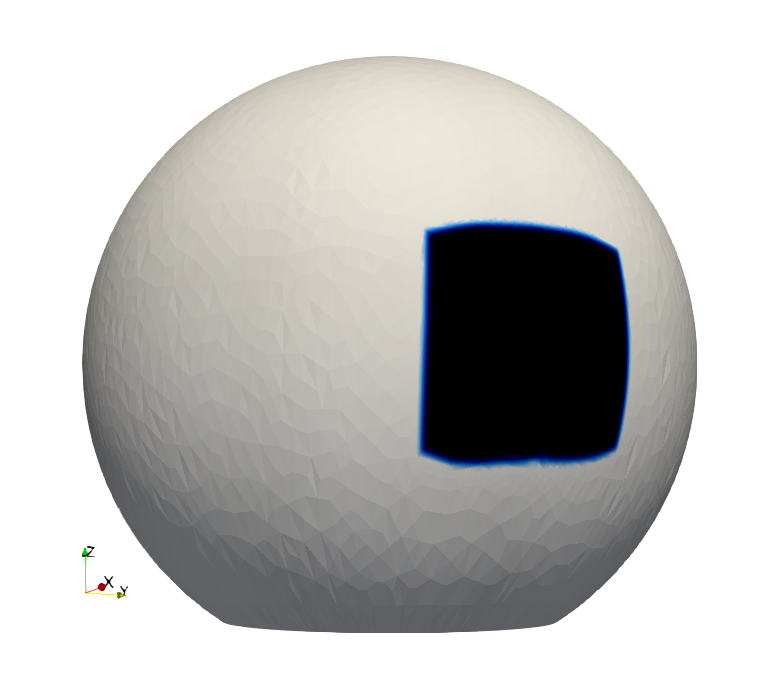}
\caption{($\epsilon=(32\pi)^{-1}$) 
Parameters as in \cite[Fig.\ 9]{vch}, but here $\epsilon=(32\pi)^{-1}$
and $\uD = -8$. Starting seed on top (top), at the front (middle) and on the
right (bottom). Displayed times are $t = 0, 0.05, 0.1$ (top) and
$t = 0, 0.01, 0.02$ (middle and bottom).
}
\label{fig:spherecap-10}
\end{figure}

Some first computations are shown in Figure~\ref{fig:spherecap-10},
where for $\epsilon=(32\pi)^{-1}$ we select the discretization parameters 
$N_c = 16$, $N_f = 128$ and $\tau = 10^{-5}$.
We start with an initial seed of radius $r_0=0.02$: either at the north pole,
or on the equator on the $x$-axis, or on the equator on the $y$-axis.
The evolutions in Figure~\ref{fig:spherecap-10} make it clear that the
hexagonal aspect of the anisotropy \eqref{eq:L4abc}(a) only comes to the fore 
at
the north pole. At the equator, on the other hand, we see interfaces with a 
four-fold symmetry, consistent with the Wulff shape displayed in 
Figure~\ref{fig:wulff3d}.
We remark that these computations are for \eqref{eq:varrho}(ii), compare with 
\cite[(2.11)(ii)]{vch}. Repeating the same simulations for
\eqref{eq:varrho}(i) leads to the creation of boundary layers, which is why
we prefer the choice \eqref{eq:varrho}(ii) here and from now. 
We refer to \cite{vch} for a more detailed discussion of this aspect.

\begin{figure}
\center
\includegraphics[angle=-0,width=0.22\textwidth]{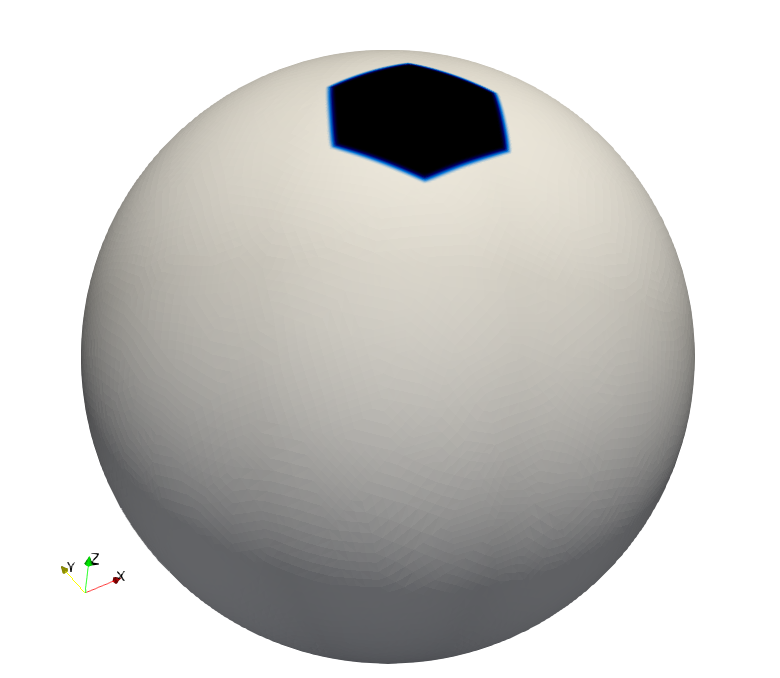}
\includegraphics[angle=-0,width=0.22\textwidth]{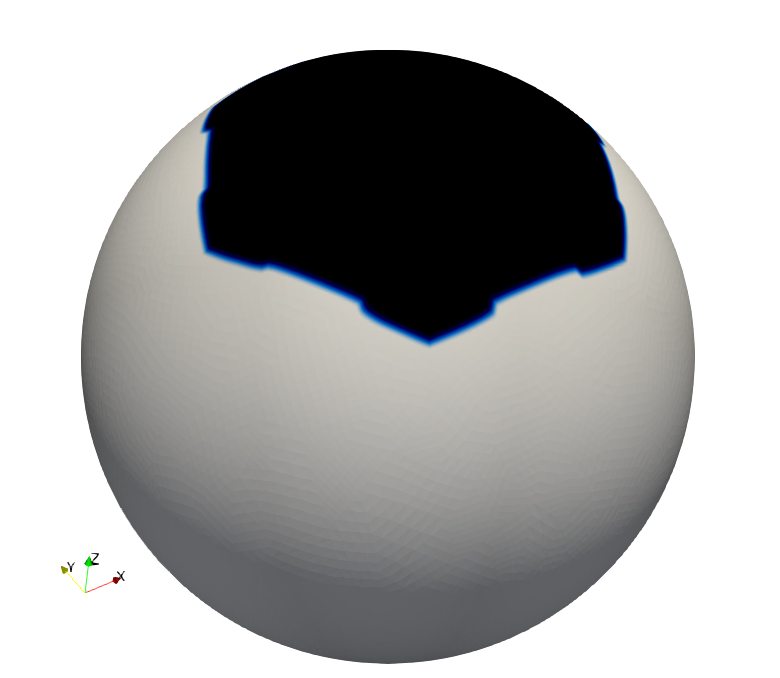}
\includegraphics[angle=-0,width=0.22\textwidth]{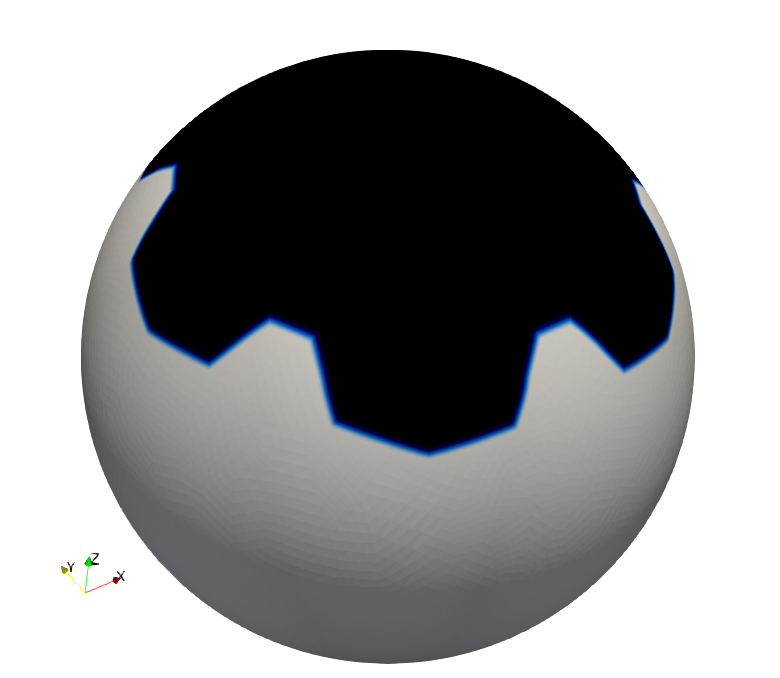}
\includegraphics[angle=-0,width=0.22\textwidth]{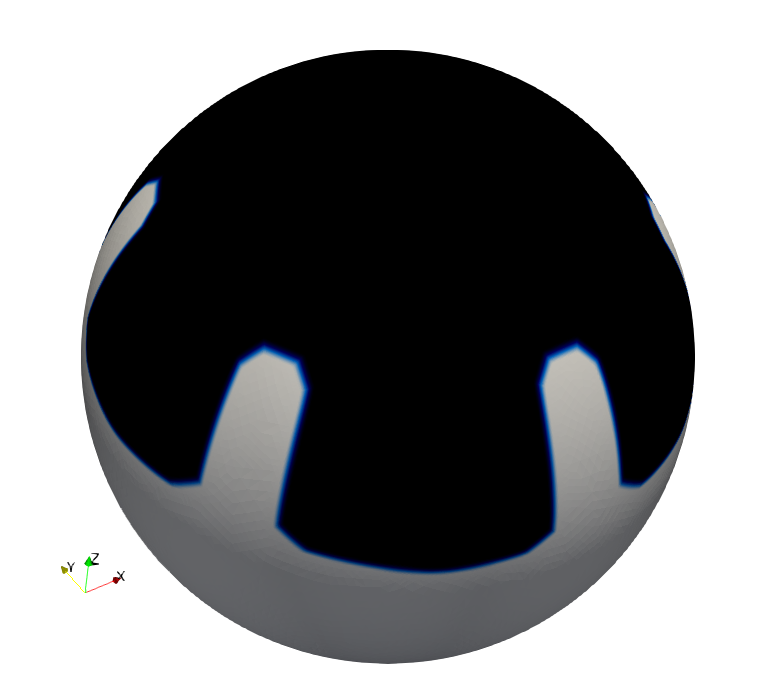} 
\includegraphics[angle=-0,width=0.22\textwidth]{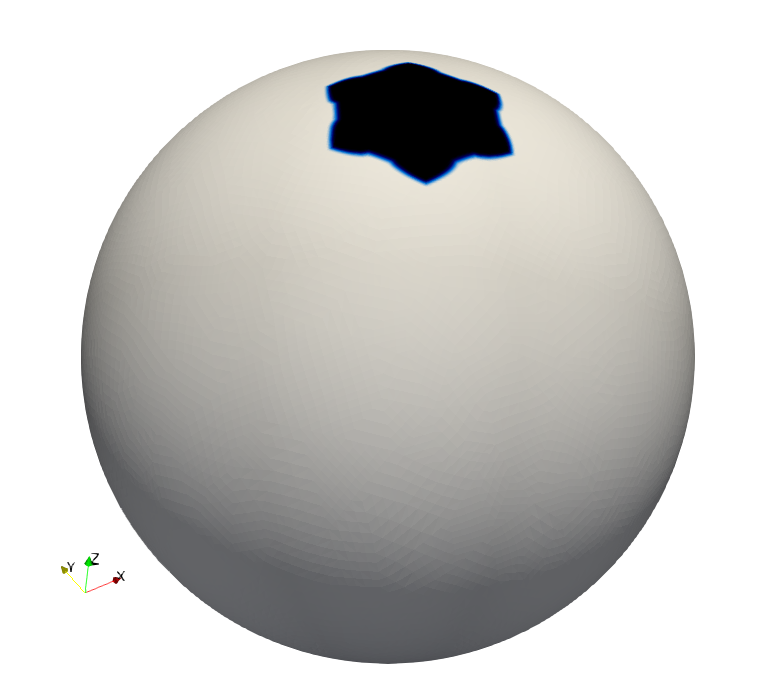}
\includegraphics[angle=-0,width=0.22\textwidth]{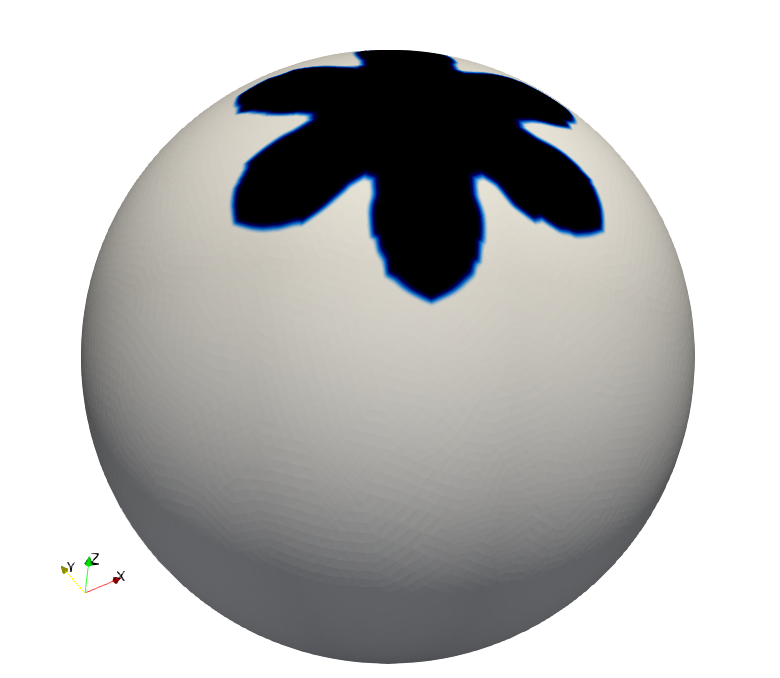}
\includegraphics[angle=-0,width=0.22\textwidth]{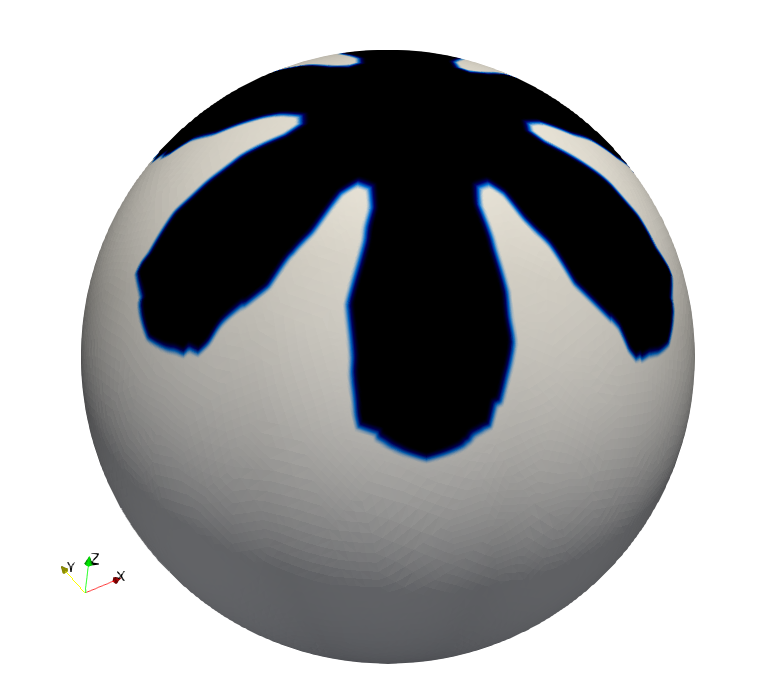}
\includegraphics[angle=-0,width=0.22\textwidth]{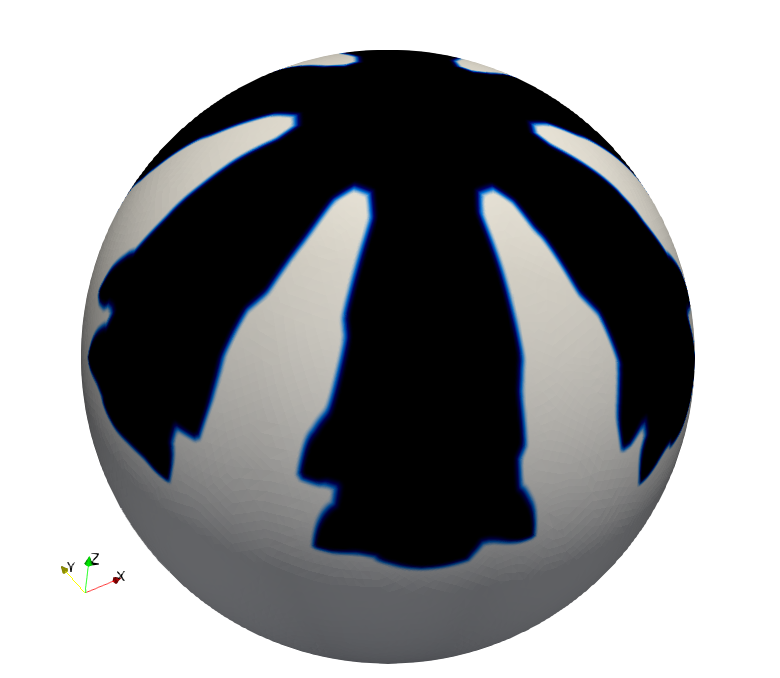} 
\includegraphics[angle=-0,width=0.22\textwidth]{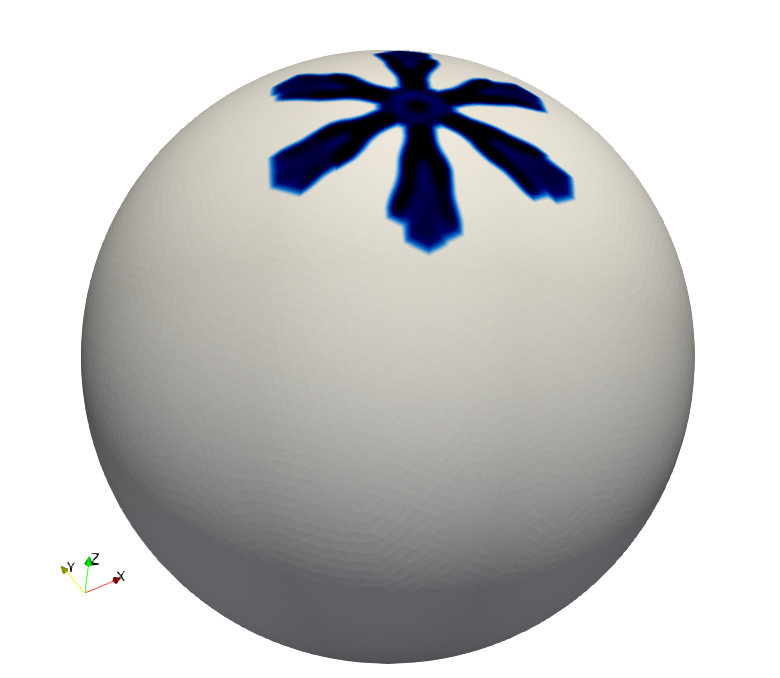}
\includegraphics[angle=-0,width=0.22\textwidth]{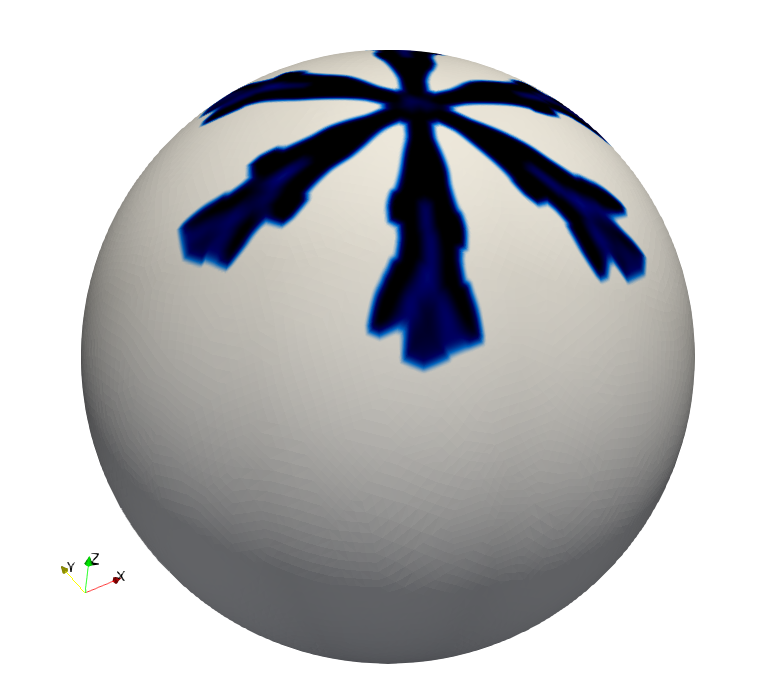}
\includegraphics[angle=-0,width=0.22\textwidth]{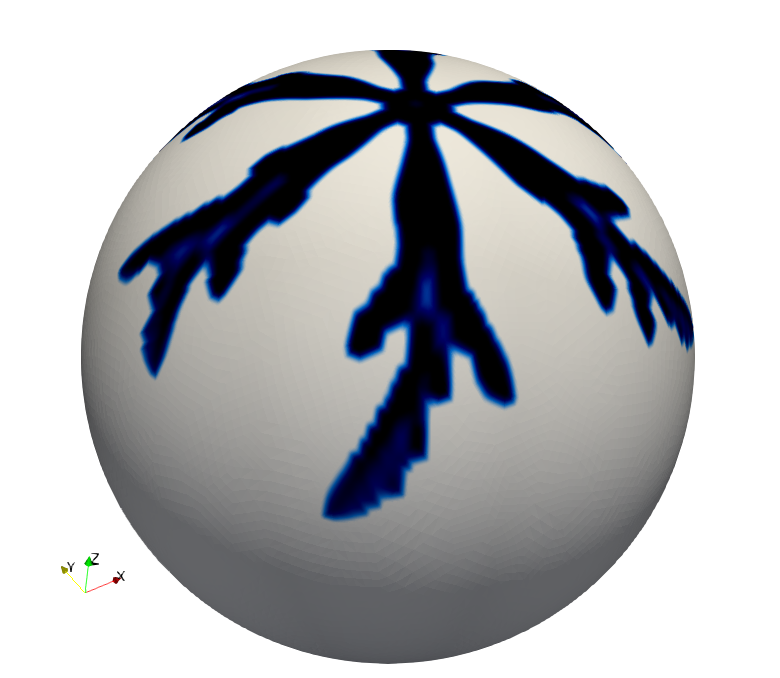}
\includegraphics[angle=-0,width=0.22\textwidth]{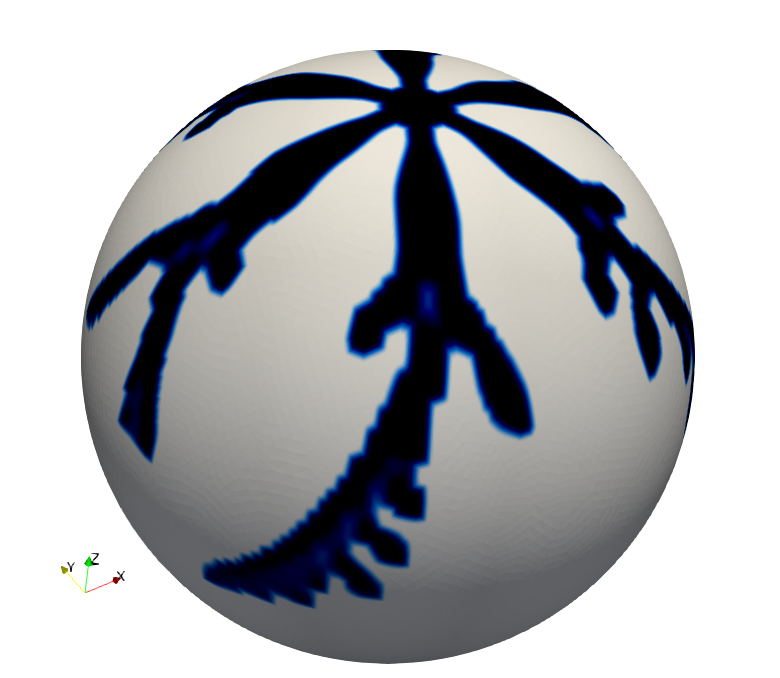}
\caption{($\epsilon=(32\pi)^{-1}$) 
Parameters as in Figure~\ref{fig:spherecap-10} but $\uD = -12$ and:
$\rho=\alpha = 0.01$ (top), $\rho = 10\alpha = 0.001$ (middle) and
$\rho = \alpha = 0.001$ (bottom).
Displayed times are $t=0.01, 0.05, 0.08, 0.12$ (top), 
$t=0.01, 0.03, 0.06, 0.09$ (middle)
and $t=0.01, 0.02, 0.03, 0.04$ (bottom).
}
\label{fig:spherecap-10more12f}
\end{figure}

\begin{figure}
\center
\includegraphics[angle=-0,width=0.22\textwidth]{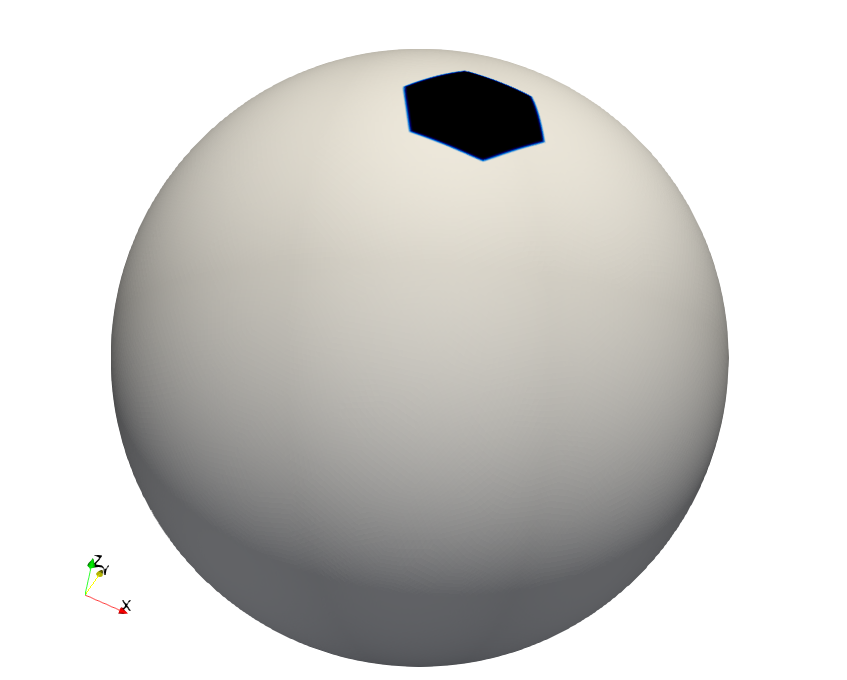}
\includegraphics[angle=-0,width=0.22\textwidth]{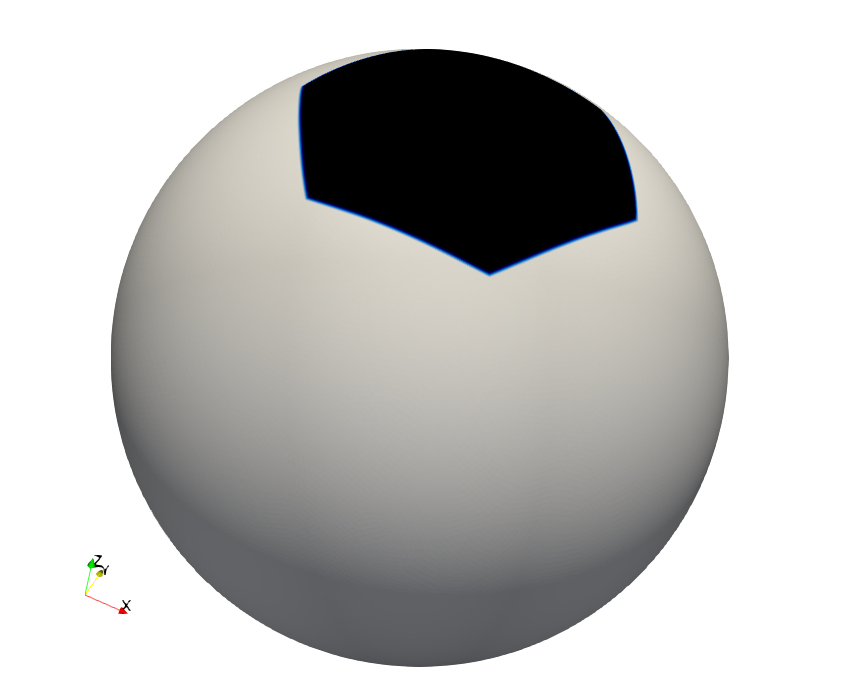}
\includegraphics[angle=-0,width=0.22\textwidth]{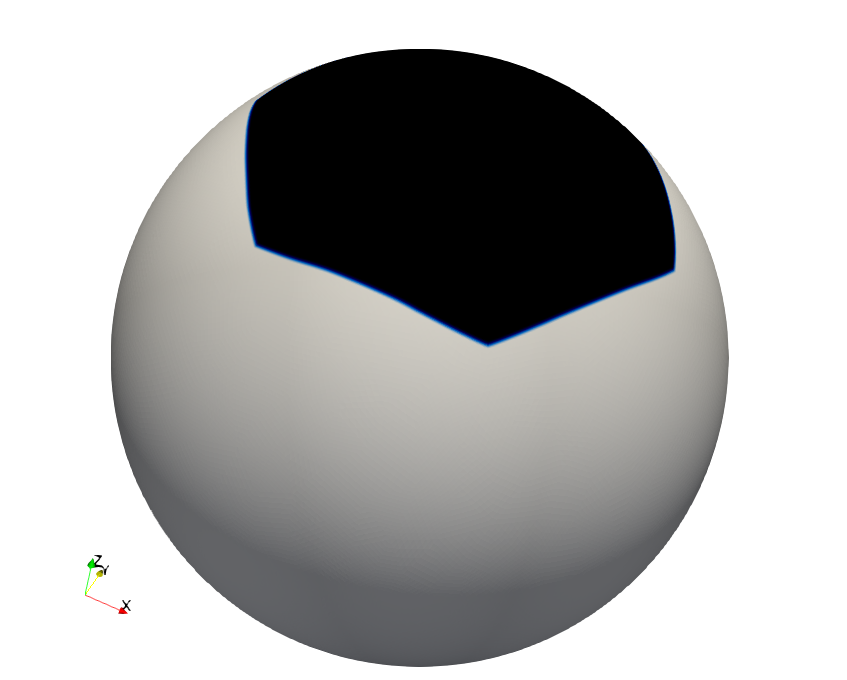}
\includegraphics[angle=-0,width=0.22\textwidth]{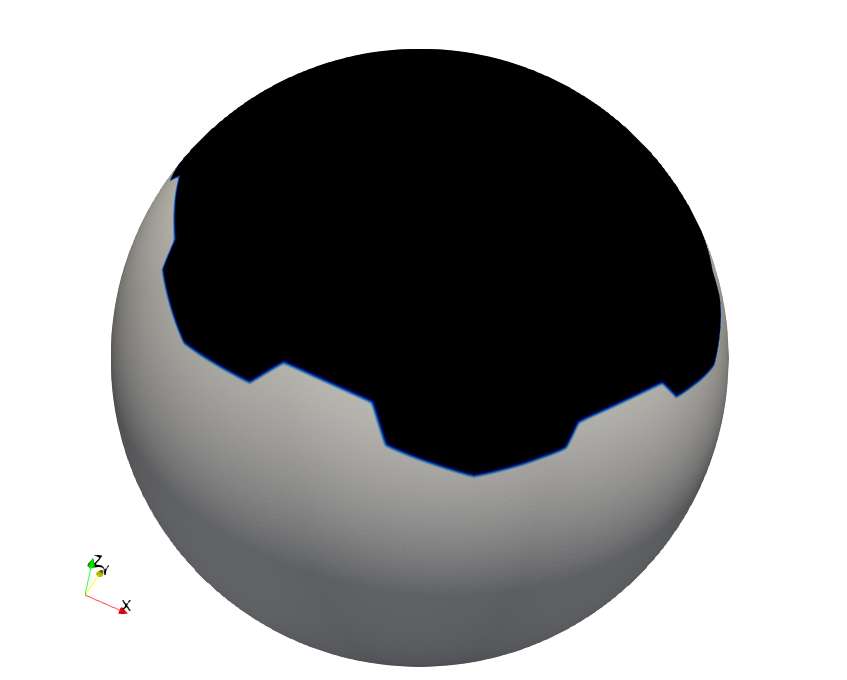} 
\includegraphics[angle=-0,width=0.22\textwidth]{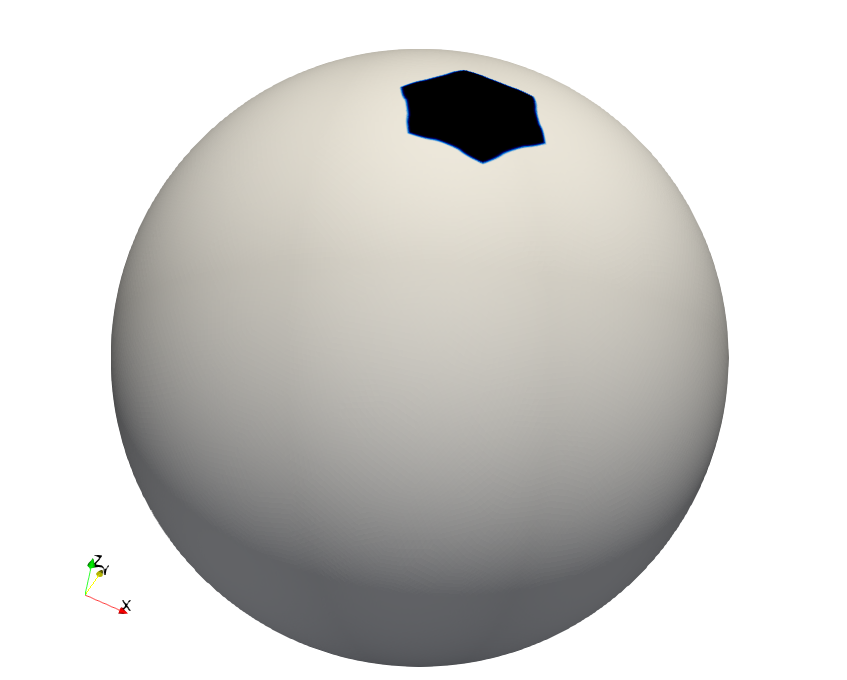}
\includegraphics[angle=-0,width=0.22\textwidth]{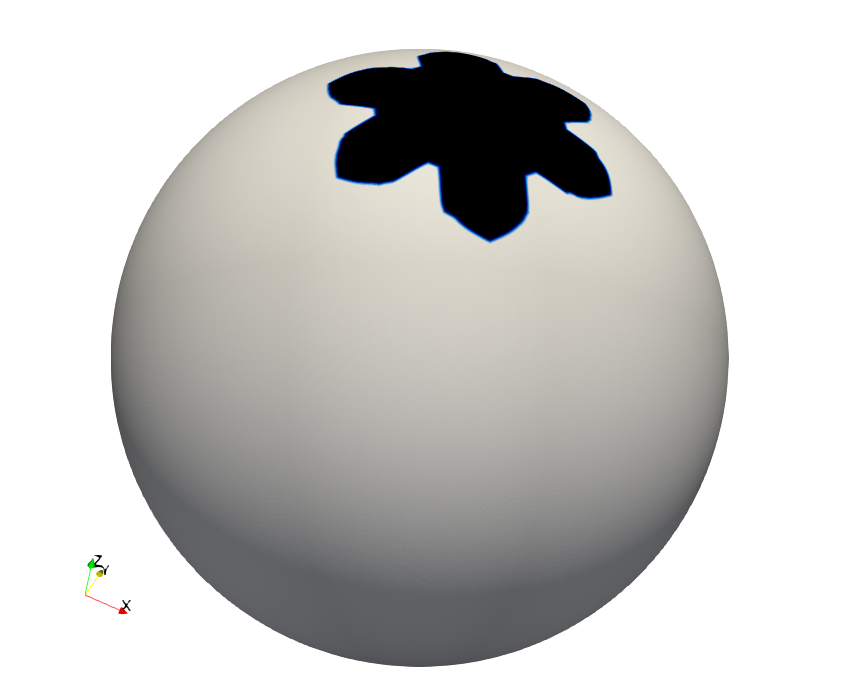}
\includegraphics[angle=-0,width=0.22\textwidth]{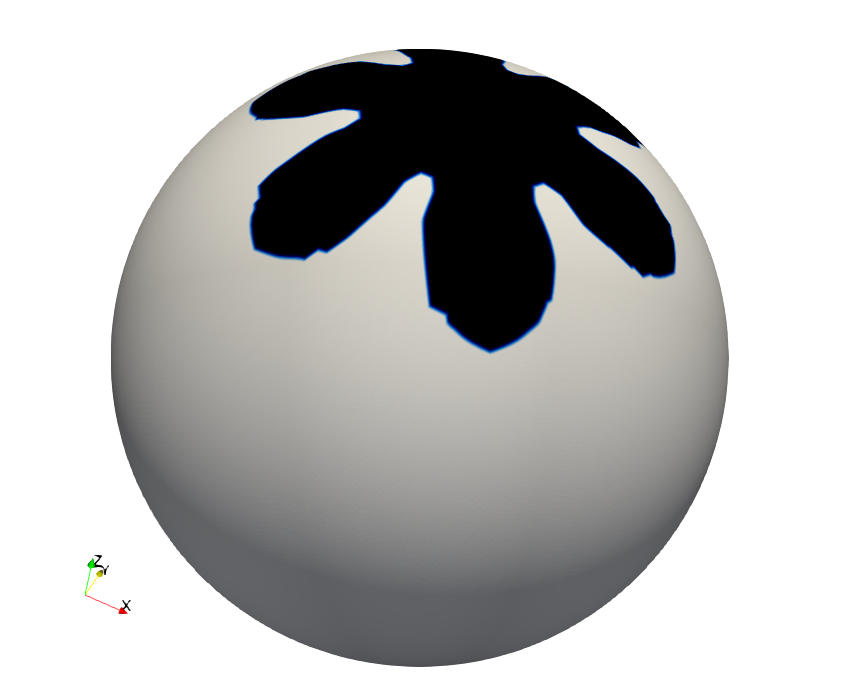}
\includegraphics[angle=-0,width=0.22\textwidth]{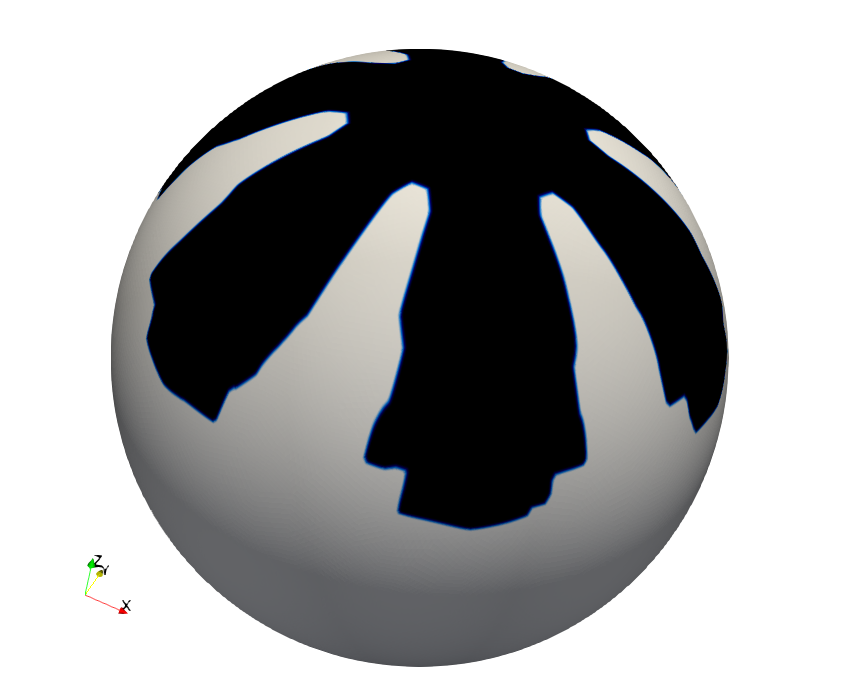} 
\includegraphics[angle=-0,width=0.22\textwidth]{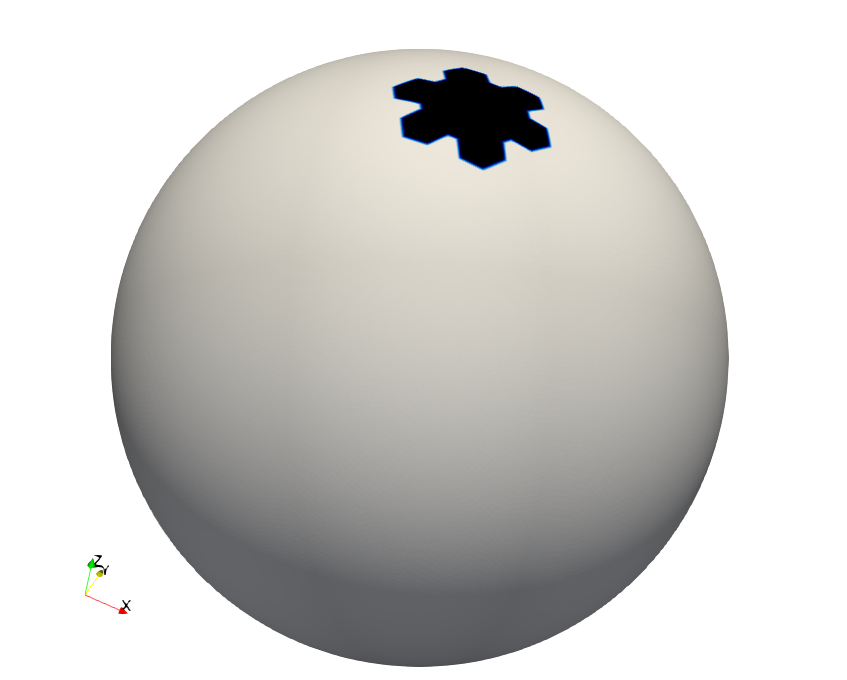}
\includegraphics[angle=-0,width=0.22\textwidth]{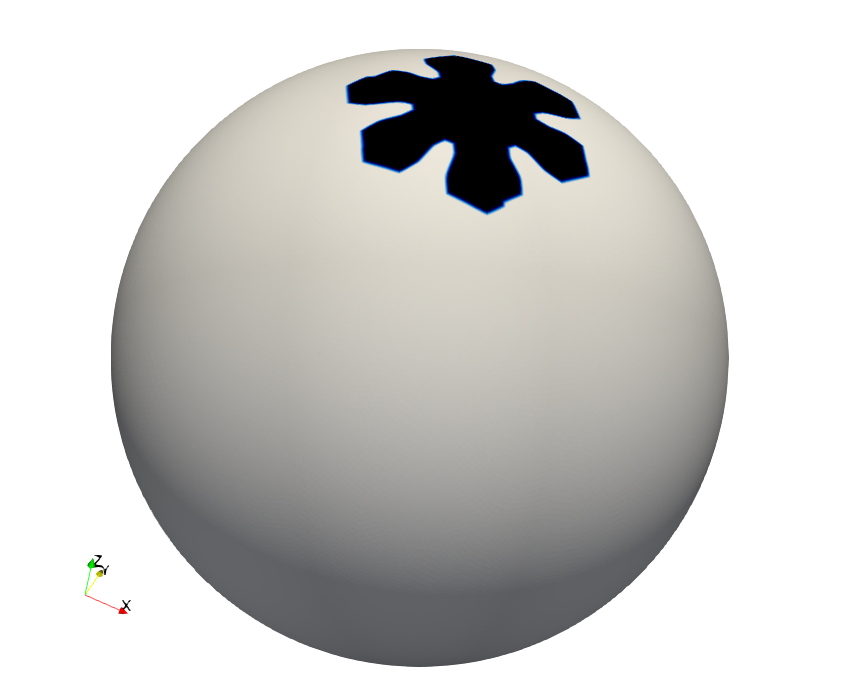}
\includegraphics[angle=-0,width=0.22\textwidth]{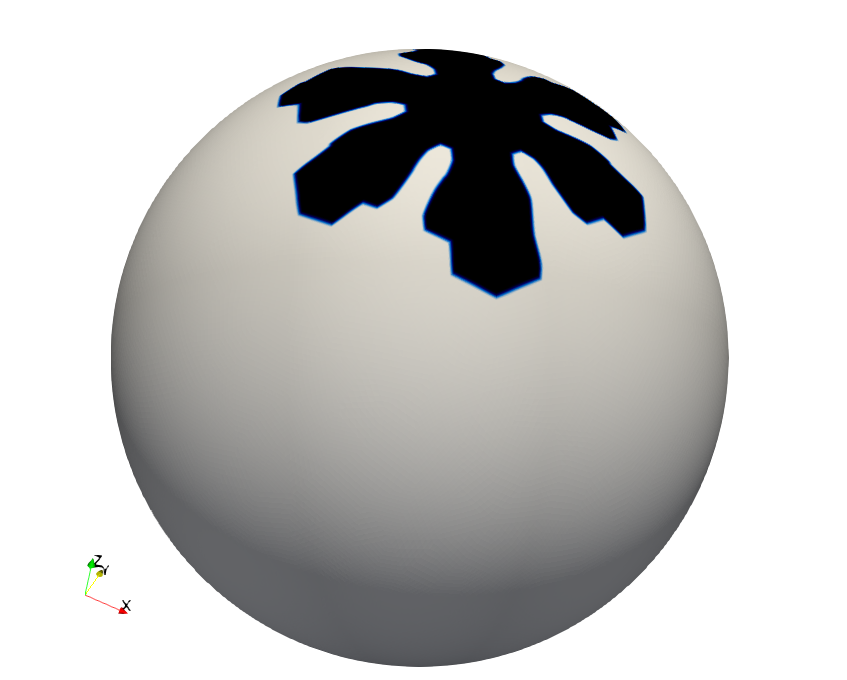}
\includegraphics[angle=-0,width=0.22\textwidth]{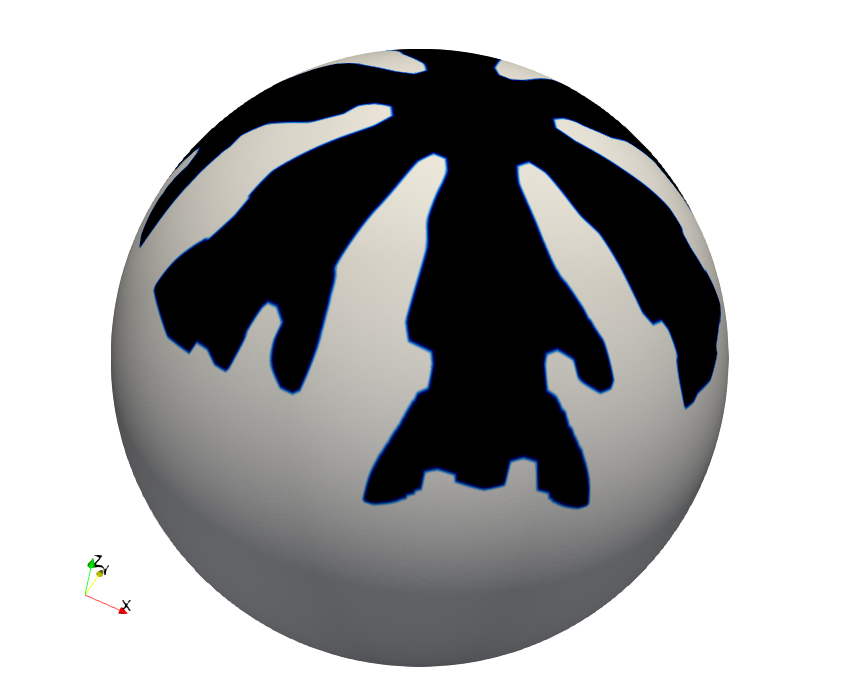}
\caption{($\epsilon=(64\pi)^{-1}$) 
Parameters as in Figure~\ref{fig:spherecap-10more12f} but with $\uD=-8$ and
starting from a smaller initial seed, with radius $r_0=0.005$.
Displayed times are $t=0.01, 0.05, 0.08, 0.14$ (top), 
$t=0.01, 0.03, 0.06, 0.012$ (middle)
and $t=0.01, 0.02, 0.04, 0.1$ (bottom).
}
\label{fig:spherecap-10more8r0}
\end{figure}

Having established that for the homogeneous hexagonal anisotropies in this
section it is beneficial for six-fold structures to start with a seed 
at the north pole, we present a few more computations with a seed of radius 
$r_0=0.02$ at the north pole. Then for different physical parameters, we obtain
slightly different evolutions. But note that what they all have in common is
that the crystal will initially grow sixth arms, but these will then often
display non-hexagonal sidearms and structures.
The simulations shown in Figure~\ref{fig:spherecap-10more12f} 
are for $\uD = -12$, with either
$\rho = \alpha = 0.01$, $\rho = 10 \alpha = 0.01$ or
$\rho = \alpha = 0.001$, respectively.
Here we let $\epsilon=(32\pi)^{-1}$ and 
$N_c = 32$, $N_f = 256$ and $\tau = 10^{-5}$.
For the next simulations we left everything unchanged, apart from
$\uD=-8$ and starting from a smaller seed, with radius $r_0 = 0.005$. 
As a consequence
we also change the values of $\epsilon = (64\pi)^{-1}$, $N_c = 32$,
$N_f=512$ and $\tau = 2.5 \times 10^{-6}$.
The new results are shown in Figure~\ref{fig:spherecap-10more8r0}.
We stress once more that if the seed would not be placed at the north
pole, or if the ambient anisotropy would be rotated so that the tangent
space at the north pole is not aligned with the six-fold symmetry, then the
resulting growth would not show the desired six-fold symmetric patterns.

\subsection{Consistent 2d anisotropies on the unit sphere}

\begin{figure}
\center
\includegraphics[angle=-0,width=0.3\textwidth]{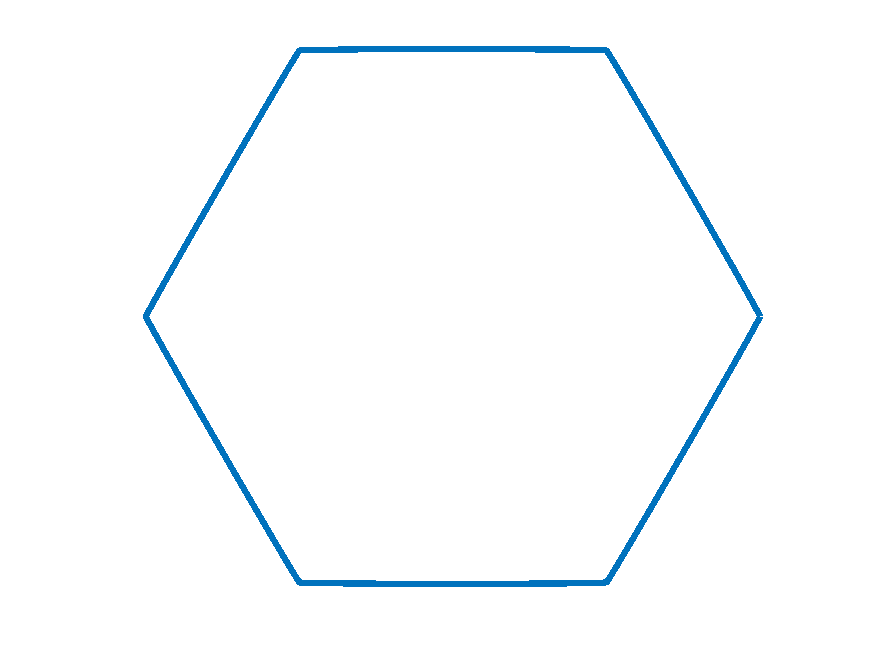}
\caption{The Wulff shape for \eqref{eq:L3}.}
\label{fig:wulff2d}
\end{figure}

In this subsection we present some numerical evidence in support of the
construction of consistent anisotropies along a surface as proposed in
\S\ref{sec:cons2d}.
We first repeat the simulations in Figure~\ref{fig:spherecap-10},
but now for the anisotropy \eqref{eq:gammasphere} with
\eqref{eq:hatg1}, for $L=3$, the usual 2d hexagonal anisotropy from, e.g.,
\cite{triplejANI,dendritic,eck,vch}, i.e.,
\begin{equation} \label{eq:L3}
\hat\gamma(\vec{p}) = 
\sum_{\ell = 1}^3
l_\delta(\hat R_1(\tfrac{\ell\,\pi}3)\,\vec{p}),
\end{equation}
with 
$\hat R_{1}(\theta):=\left(\!\!\!\scriptsize
\begin{array}{rr} \cos\theta & \sin\theta \\
-\sin\theta & \cos\theta \end{array}\!\! \right)$, 
similarly to \eqref{eq:L4abc}. We visualize the Wulff shape of \eqref{eq:L3} 
in Figure~\ref{fig:wulff2d}.

\begin{figure}
\center
\includegraphics[angle=-0,width=0.22\textwidth]{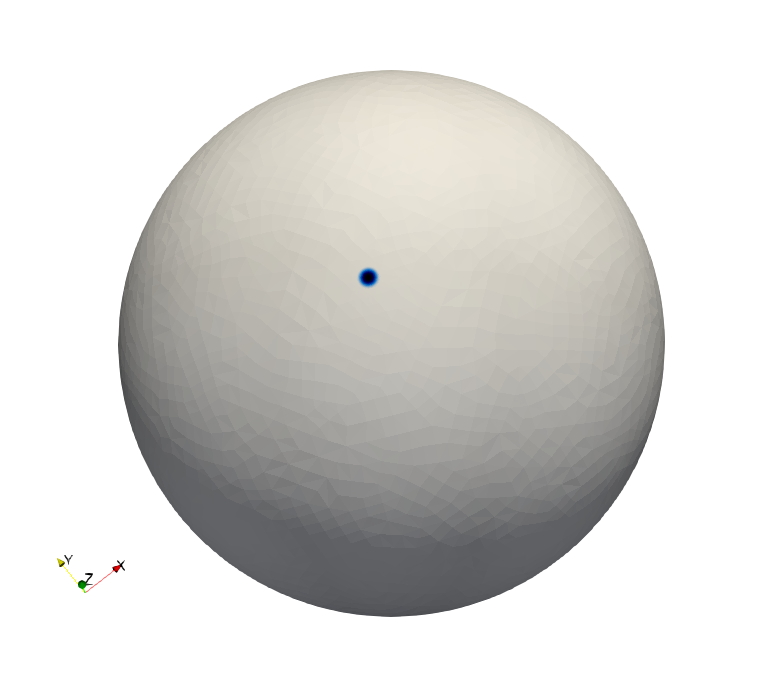}
\includegraphics[angle=-0,width=0.22\textwidth]{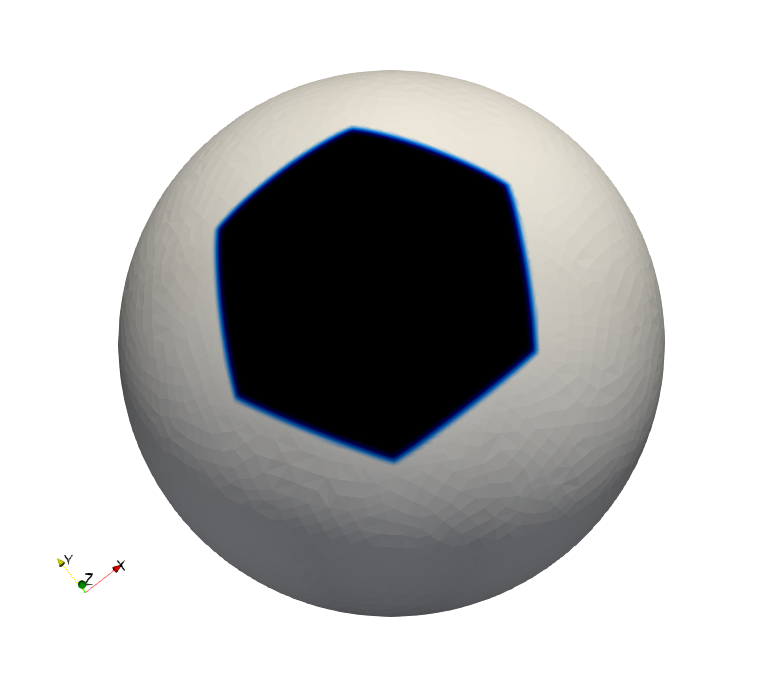}
\includegraphics[angle=-0,width=0.22\textwidth]{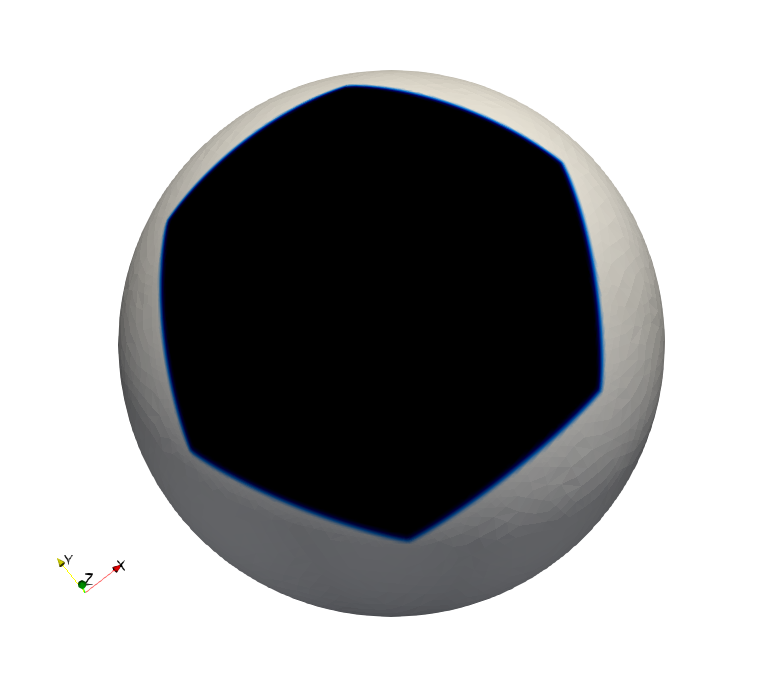}
\\
\includegraphics[angle=-0,width=0.22\textwidth]{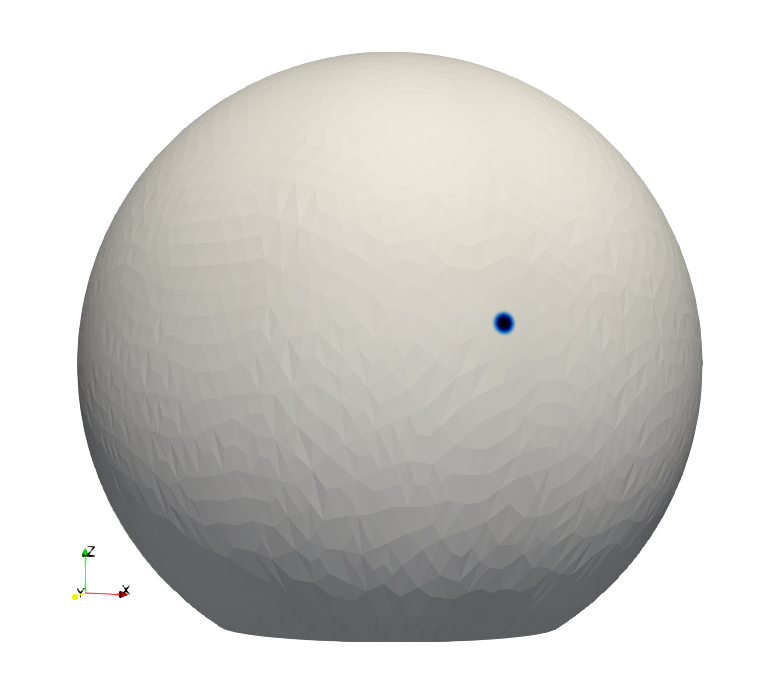}
\includegraphics[angle=-0,width=0.22\textwidth]{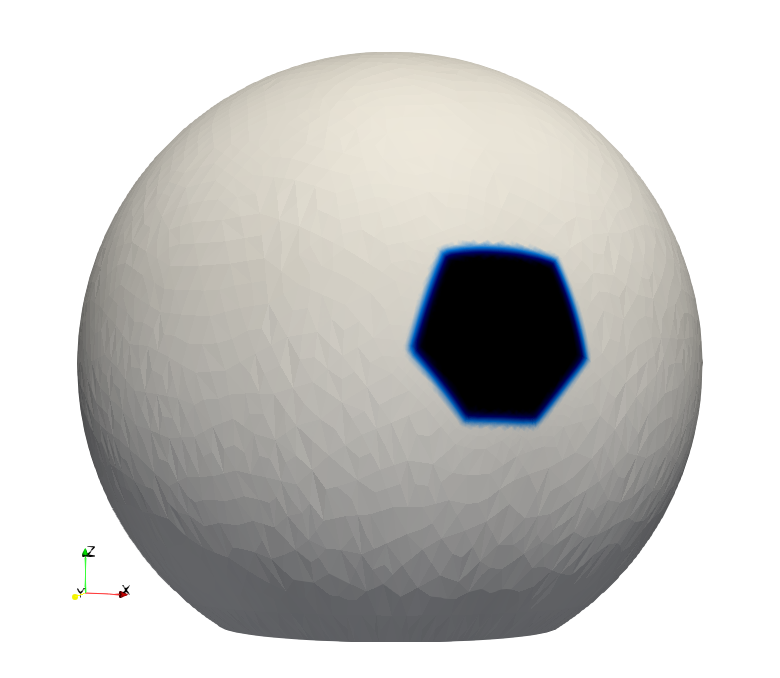}
\includegraphics[angle=-0,width=0.22\textwidth]{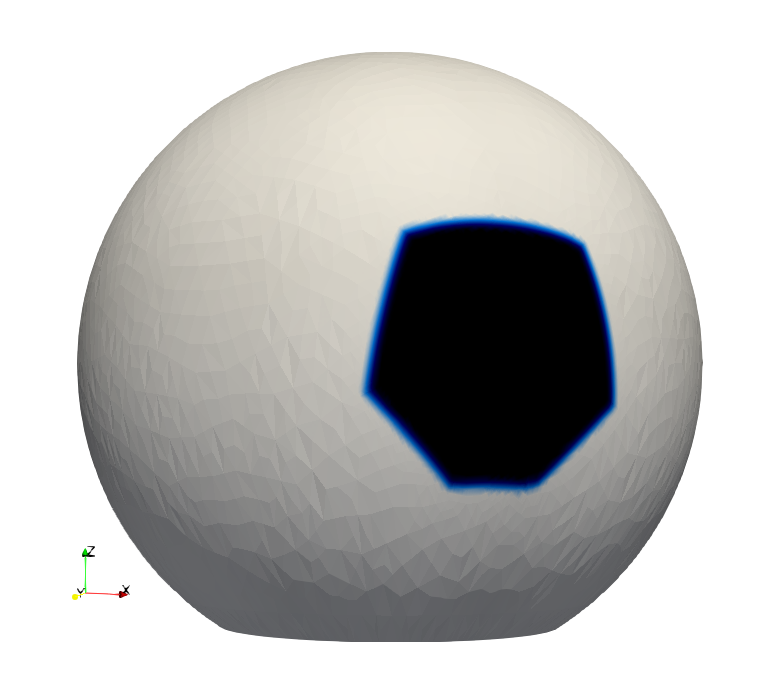}
\\
\includegraphics[angle=-0,width=0.22\textwidth]{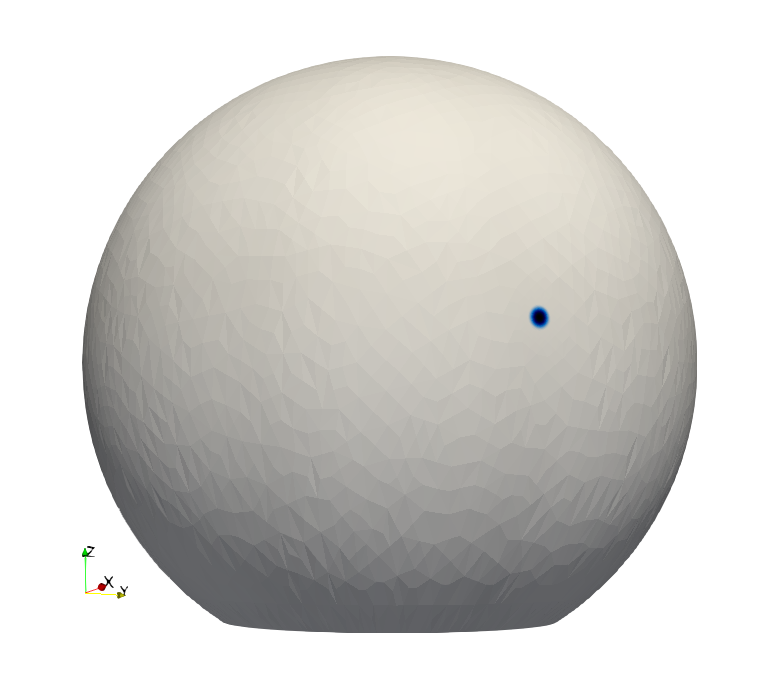}
\includegraphics[angle=-0,width=0.22\textwidth]{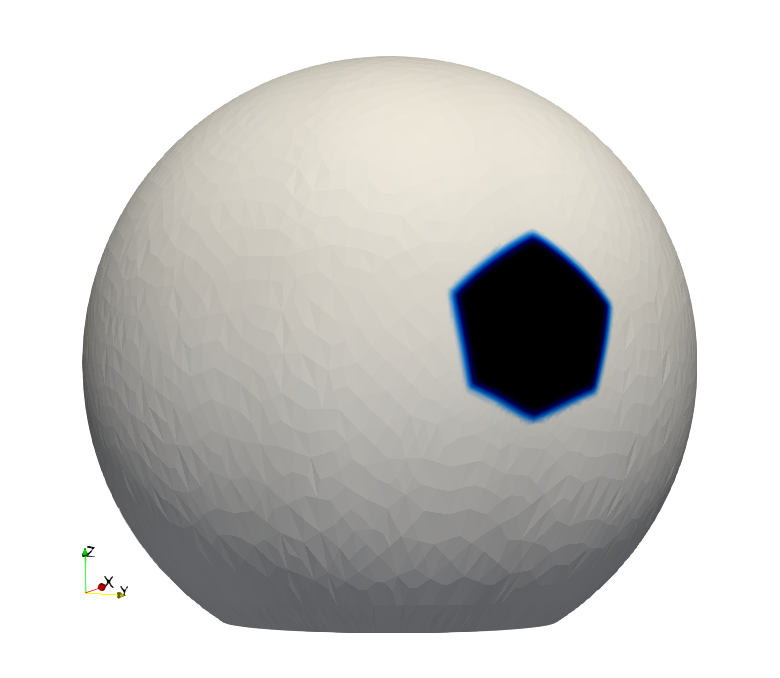}
\includegraphics[angle=-0,width=0.22\textwidth]{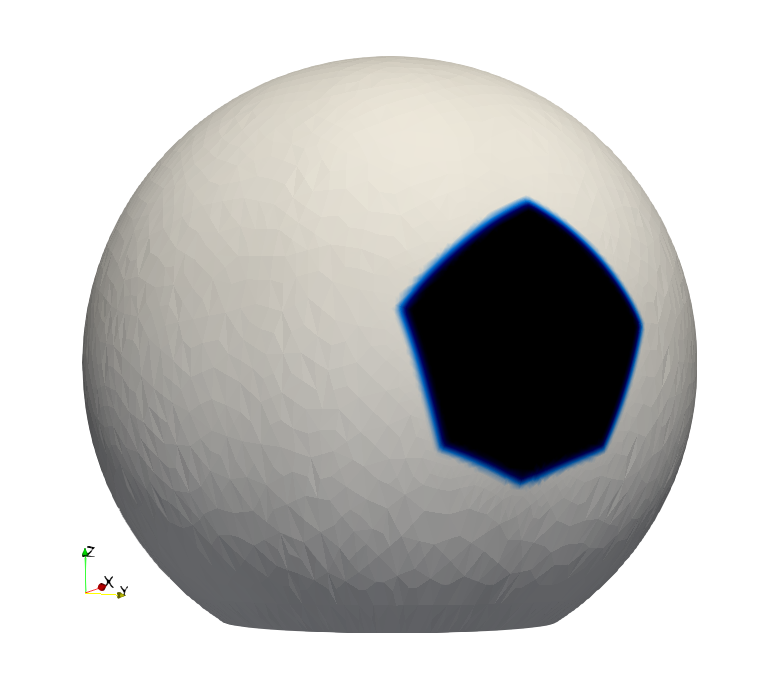}
\caption{($\epsilon=(32\pi)^{-1}$) 
Parameters as in Figure~\ref{fig:spherecap-10}, but for the
anisotropy \eqref{eq:gammasphere} with \eqref{eq:L3}.
Starting seed on top (top), at the front (middle) and on the
right (bottom). Displayed times are $t = 0, 0.05, 0.1, 0.4$ (top) and
$t = 0, 0.01, 0.02, 0.05$ (middle and bottom).
}
\label{fig:spherecap-10L103}
\end{figure}

The numerical results are shown in 
Figure~\ref{fig:spherecap-10L103}, and the difference to 
Figure~\ref{fig:spherecap-10} is obvious. While the two simulations for the 
initial seed at the north pole are qualitatively very close, the seeds that
start on the equator have very different evolutions. 
In Figure~\ref{fig:spherecap-10}, the seeds grow into a quadrilateral 
interface, while in Figure~\ref{fig:spherecap-10L103} these also grow
hexagonally. The explanation can be found in the Wulff shape of
\eqref{eq:L4abc}, which for tangent spaces on the 
equator of the sphere will lead to four-fold structures, 
rather than six-fold structures at the two poles.

\begin{figure}
\center
\includegraphics[angle=-0,width=0.22\textwidth]{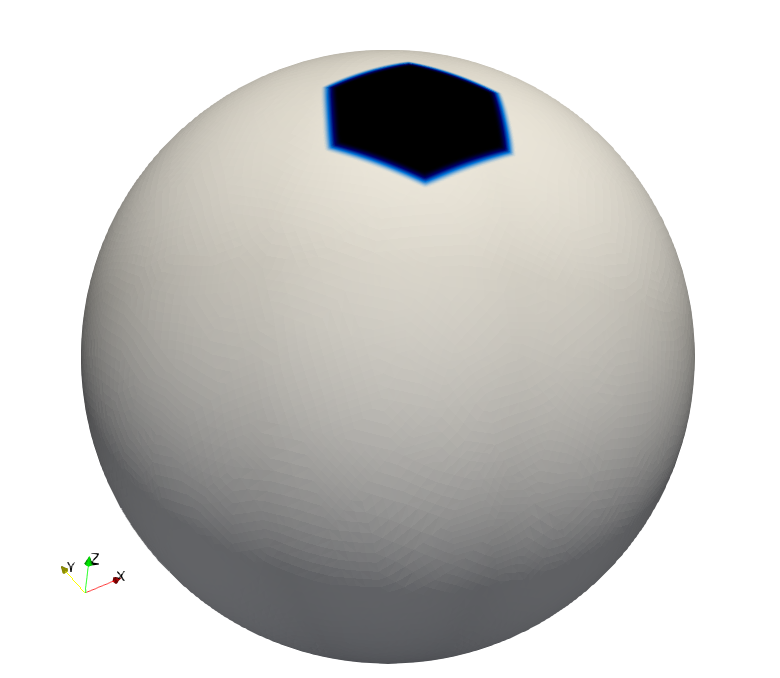}
\includegraphics[angle=-0,width=0.22\textwidth]{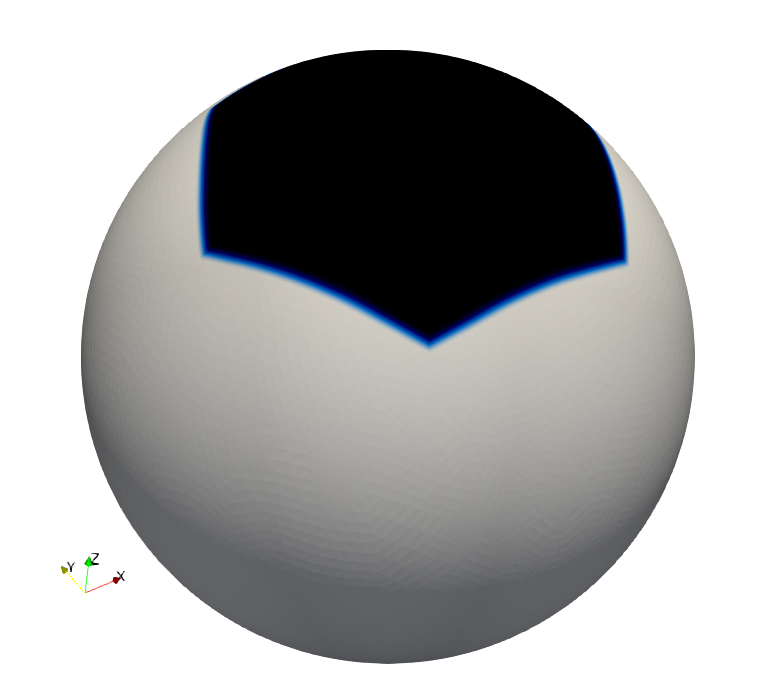}
\includegraphics[angle=-0,width=0.22\textwidth]{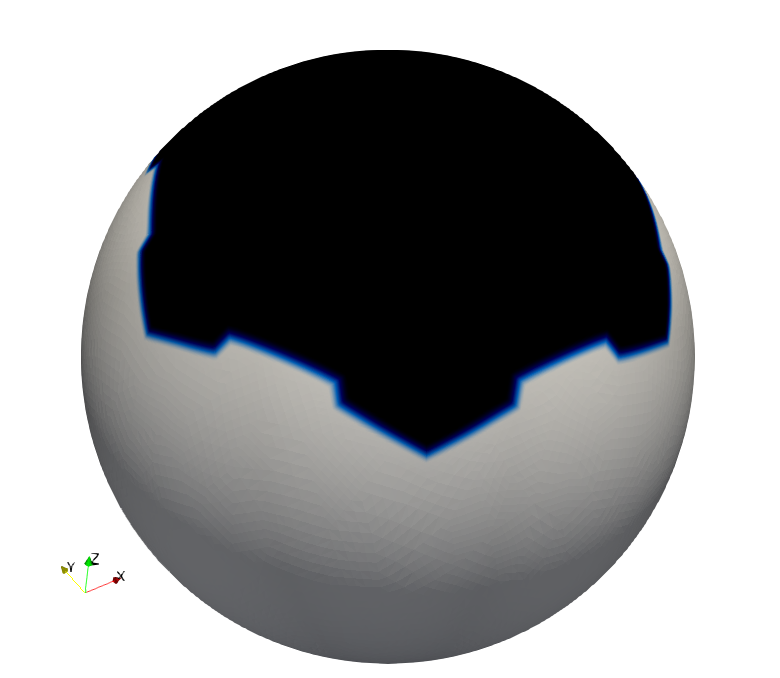}
\includegraphics[angle=-0,width=0.22\textwidth]{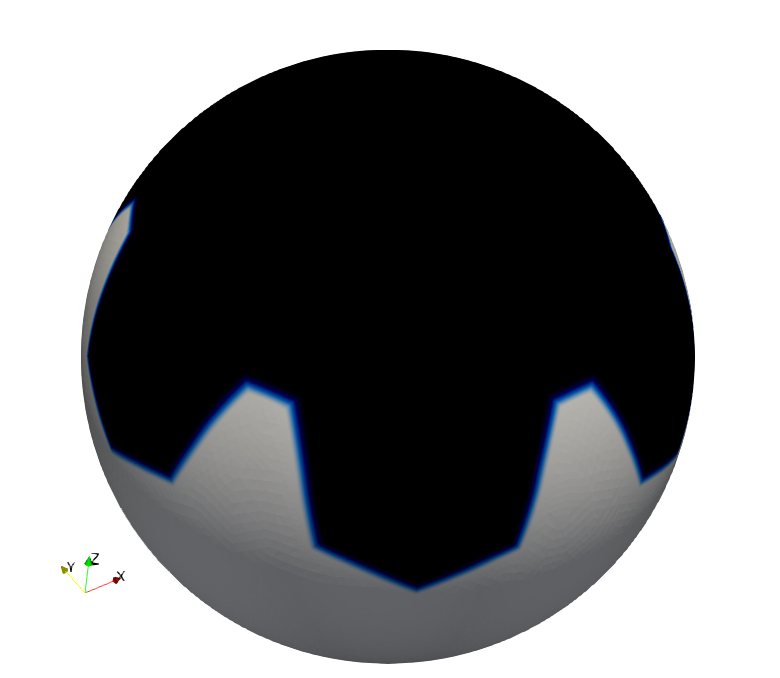} 
\includegraphics[angle=-0,width=0.22\textwidth]{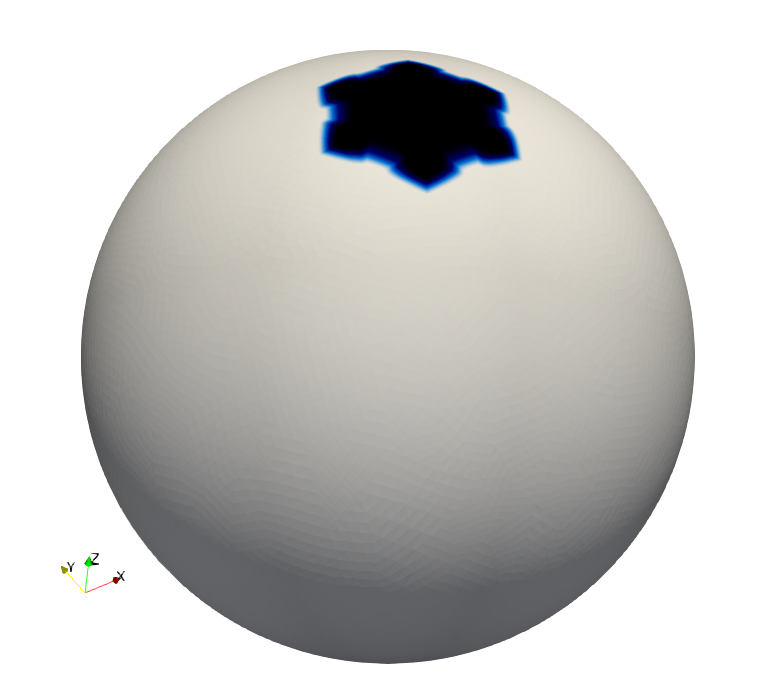}
\includegraphics[angle=-0,width=0.22\textwidth]{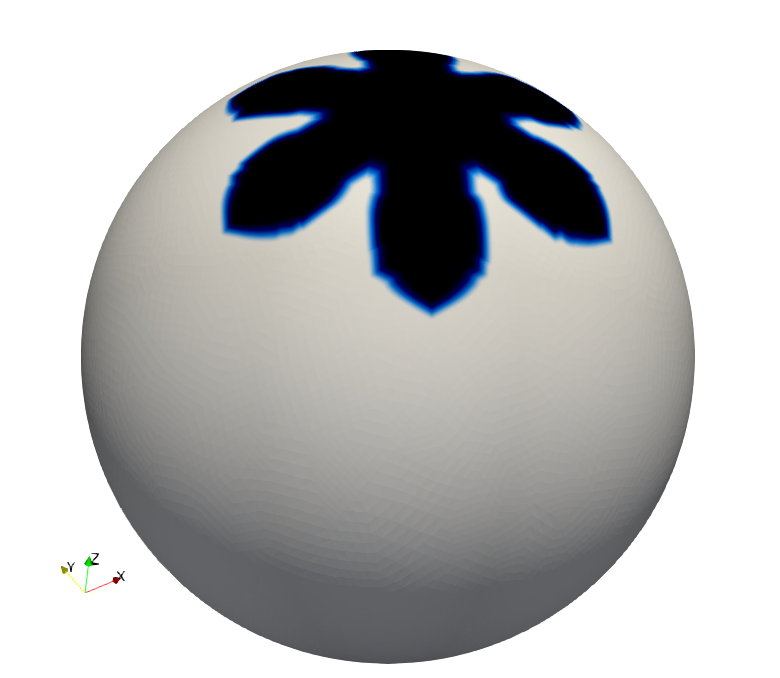}
\includegraphics[angle=-0,width=0.22\textwidth]{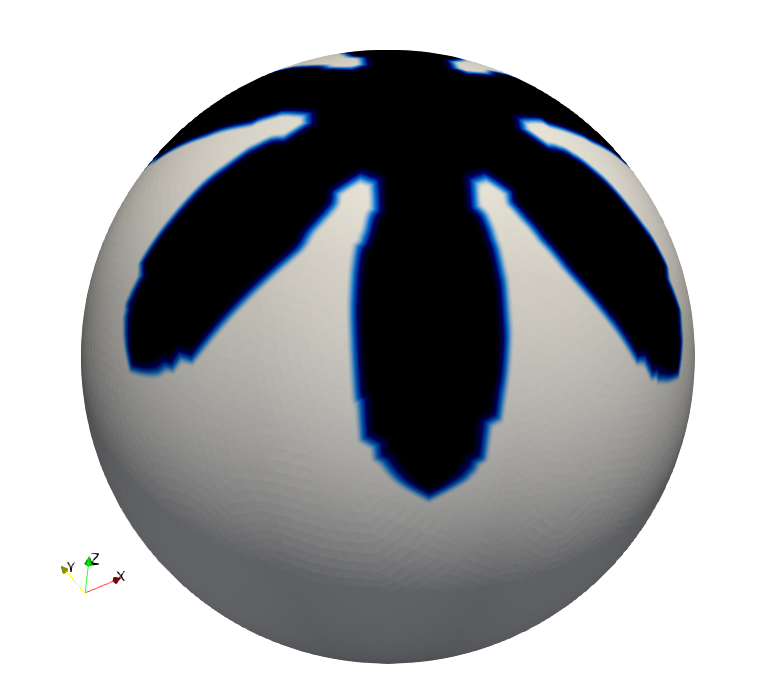}
\includegraphics[angle=-0,width=0.22\textwidth]{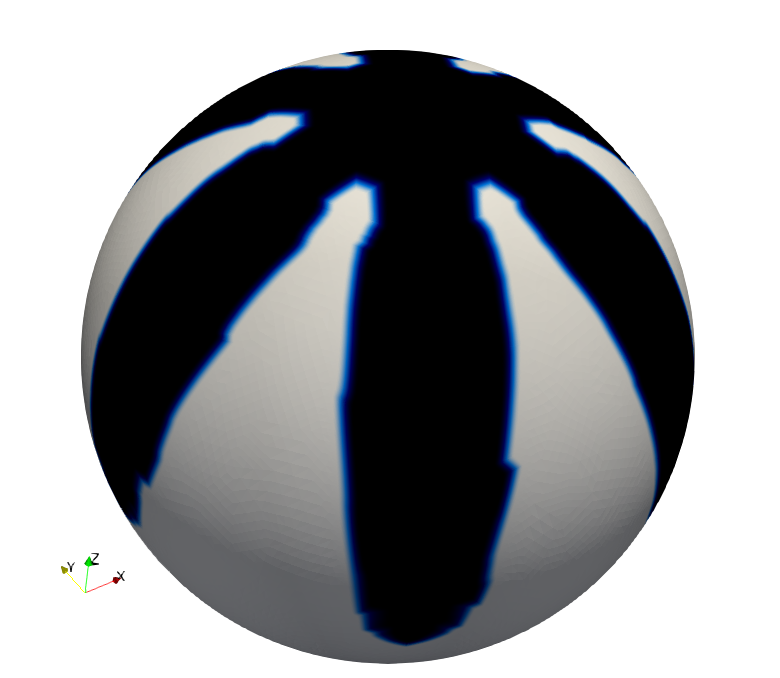} 
\includegraphics[angle=-0,width=0.22\textwidth]{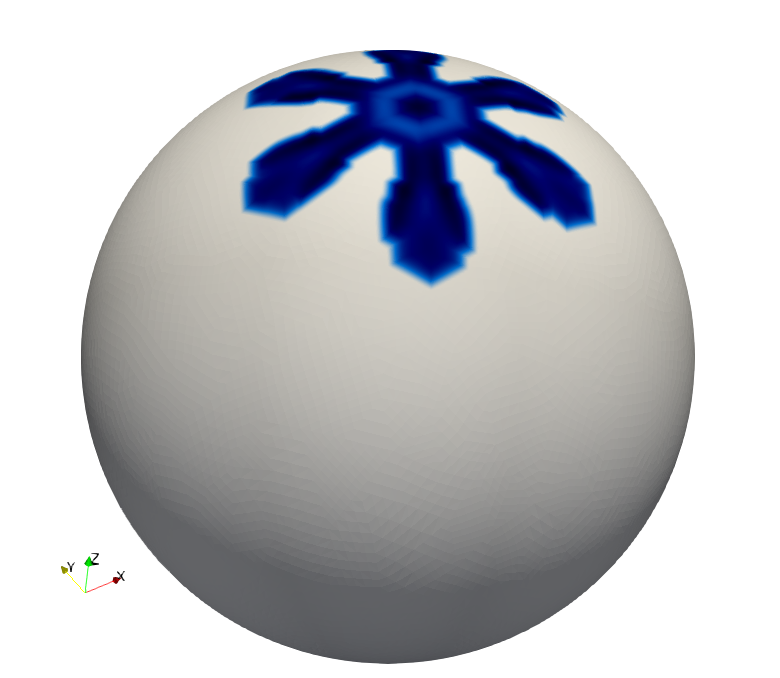}
\includegraphics[angle=-0,width=0.22\textwidth]{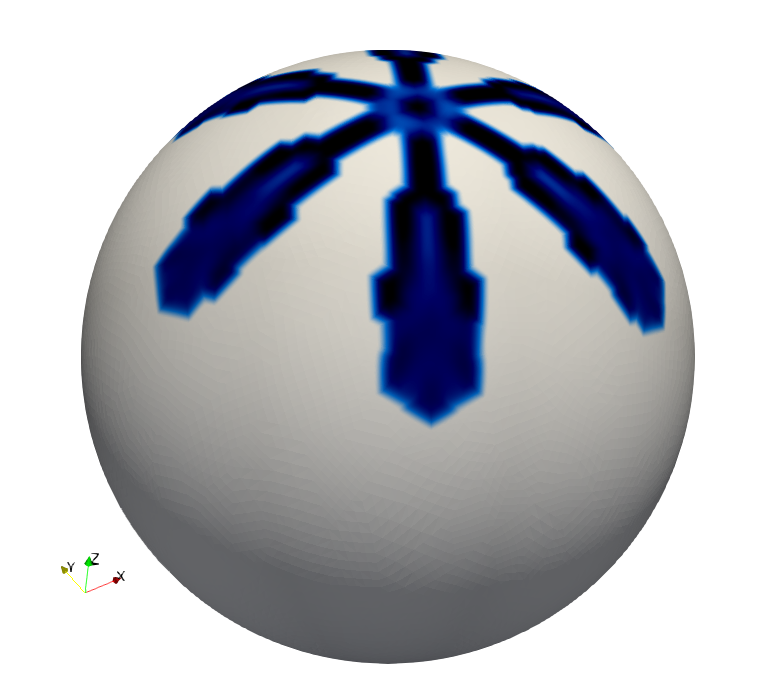}
\includegraphics[angle=-0,width=0.22\textwidth]{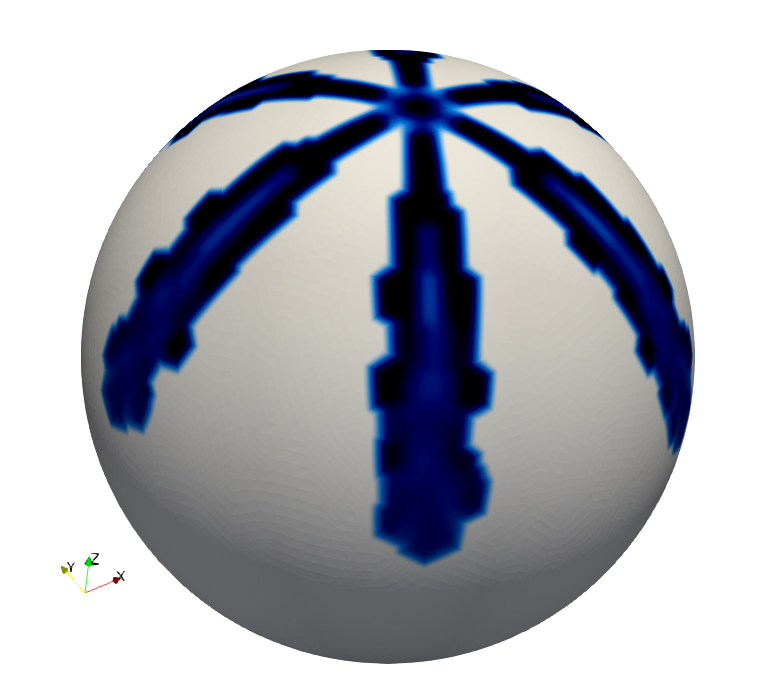}
\includegraphics[angle=-0,width=0.22\textwidth]{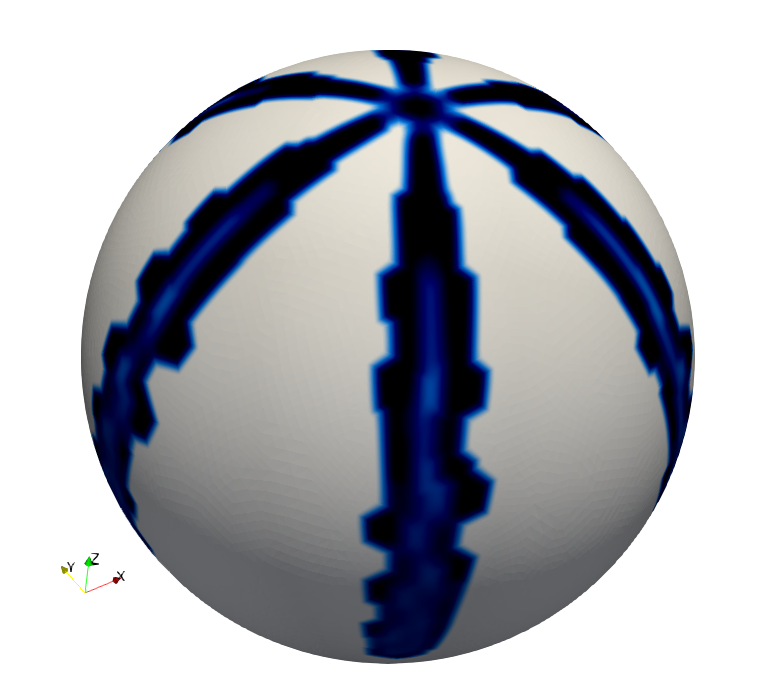}
\caption{($\epsilon=(32\pi)^{-1}$) 
Parameters as in Figure~\ref{fig:spherecap-10more12f}, but with the
anisotropy \eqref{eq:gammasphere} with \eqref{eq:L3}.
}
\label{fig:spherecap-10L103more12f}
\end{figure}

With this arguably more realistic anisotropy for ice crystal growth on a sphere
in place, we can repeat the simulations in Figure~\ref{fig:spherecap-10more12f}
now for this new anisotropy. The results are shown in 
Figure~\ref{fig:spherecap-10L103more12f}, where for the last row we observe 
the development of some mushy regions, suggesting that $\epsilon$ was not
chosen sufficiently small to resolve the physics of the underlying sharp
interface problem.

\begin{figure}
\center
\includegraphics[angle=-0,width=0.22\textwidth]{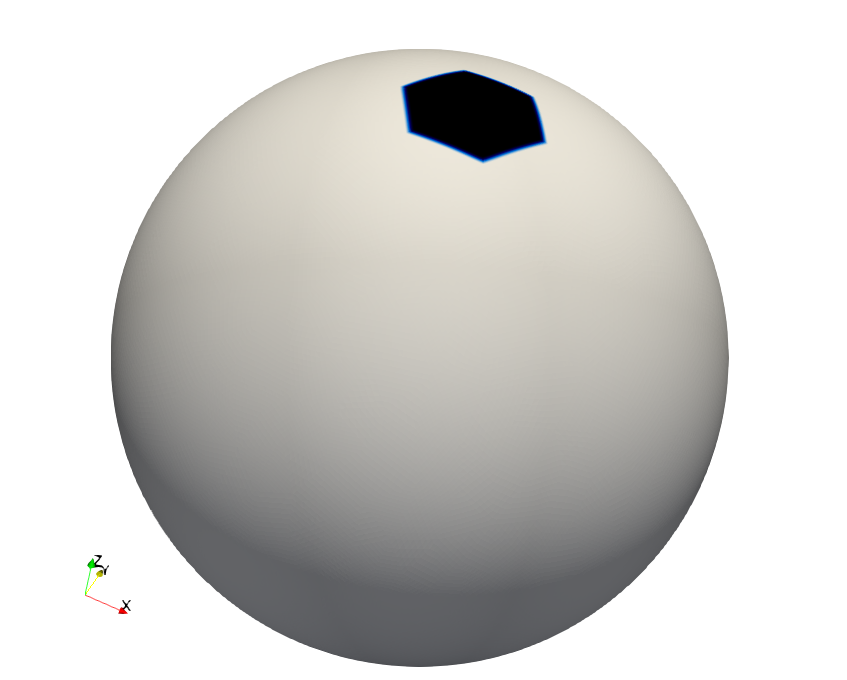}
\includegraphics[angle=-0,width=0.22\textwidth]{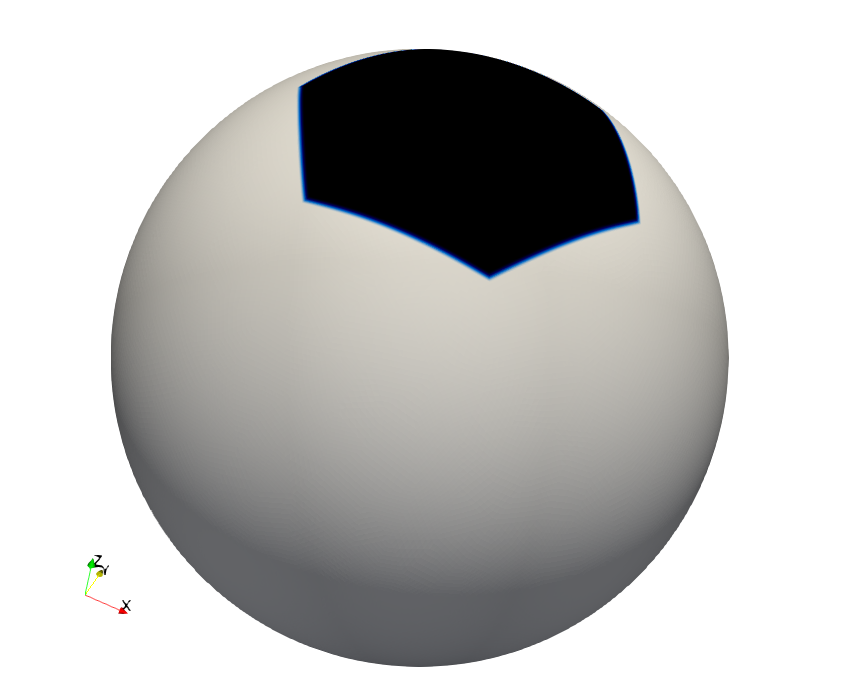}
\includegraphics[angle=-0,width=0.22\textwidth]{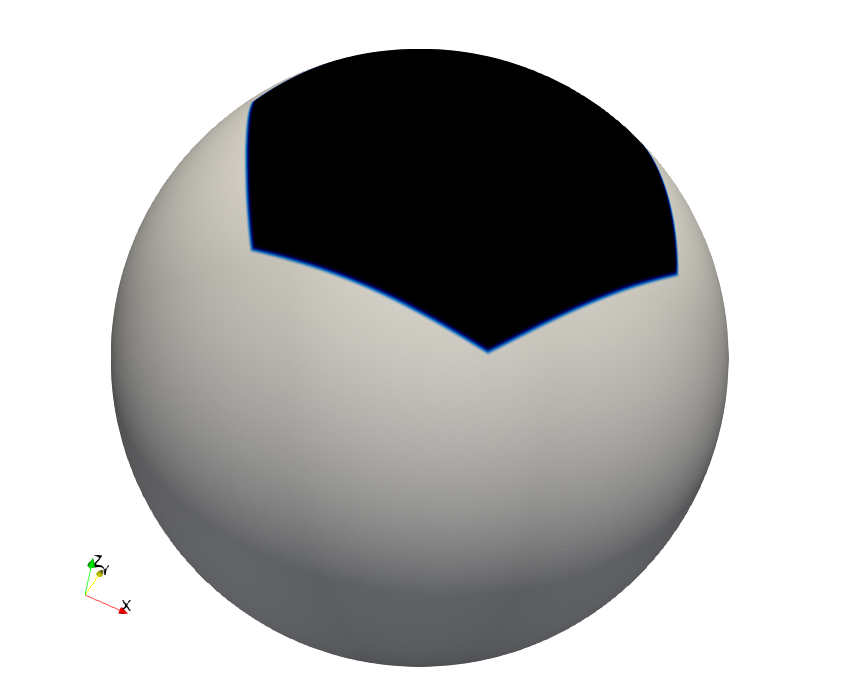}
\includegraphics[angle=-0,width=0.22\textwidth]{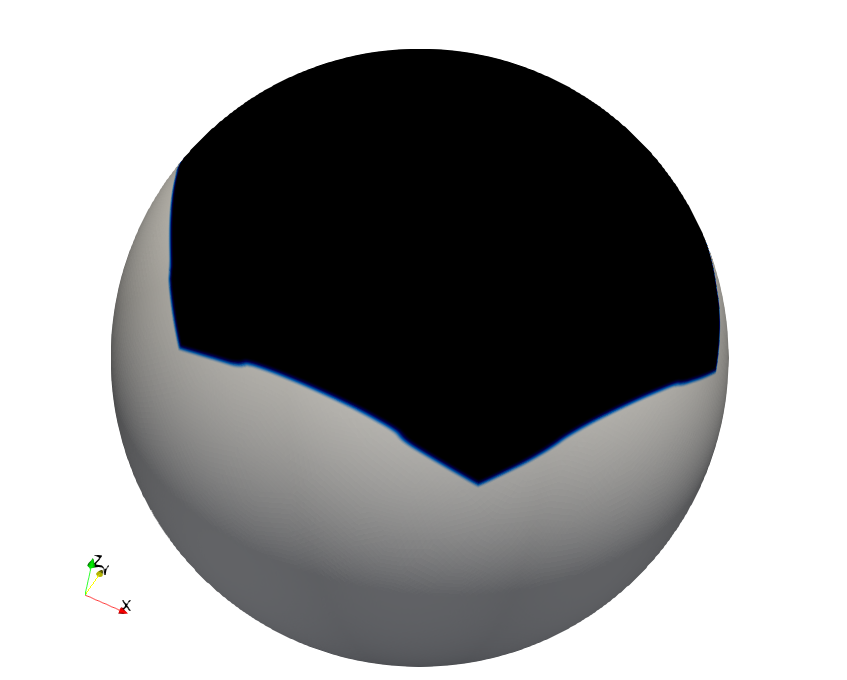} 
\includegraphics[angle=-0,width=0.22\textwidth]{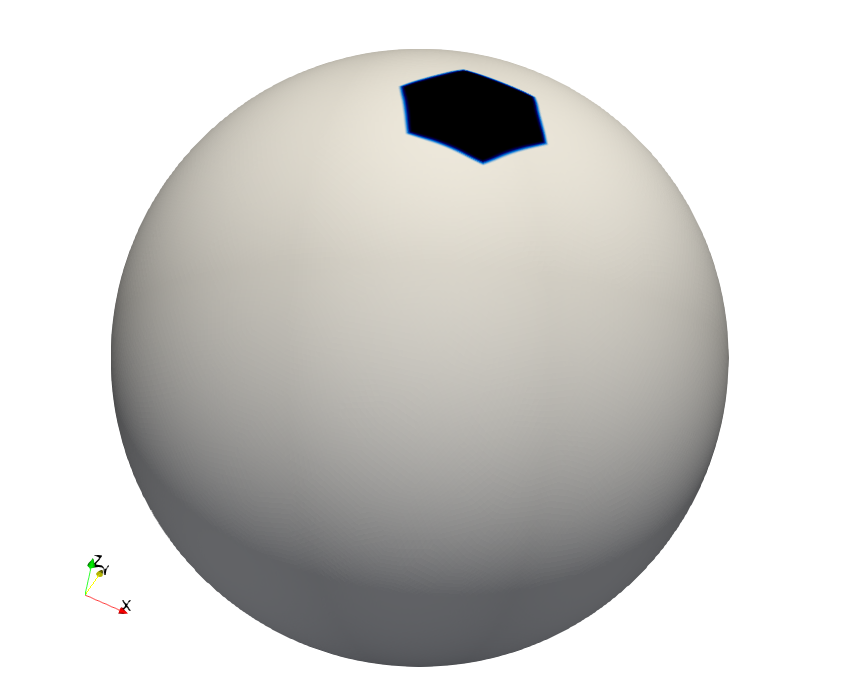}
\includegraphics[angle=-0,width=0.22\textwidth]{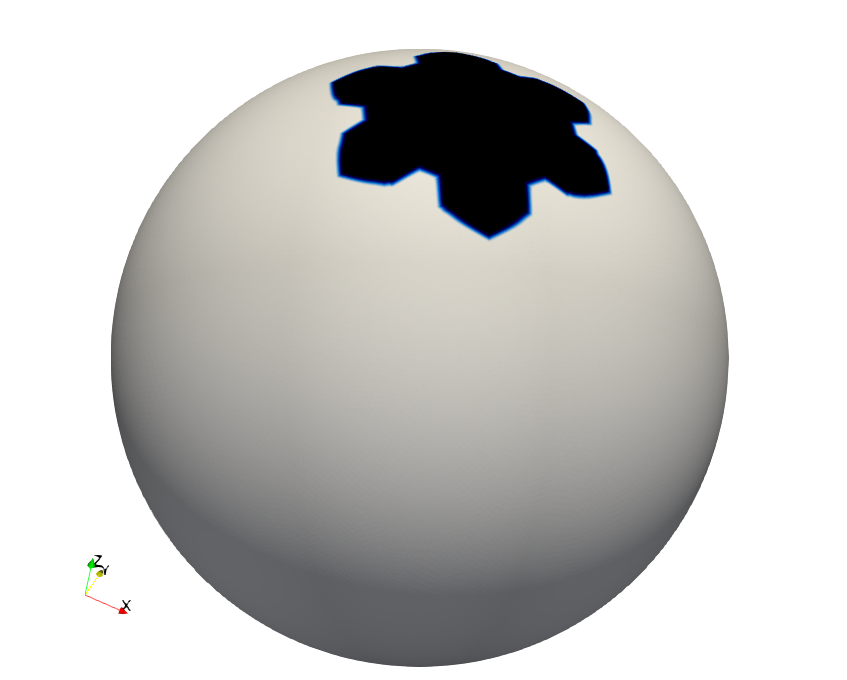}
\includegraphics[angle=-0,width=0.22\textwidth]{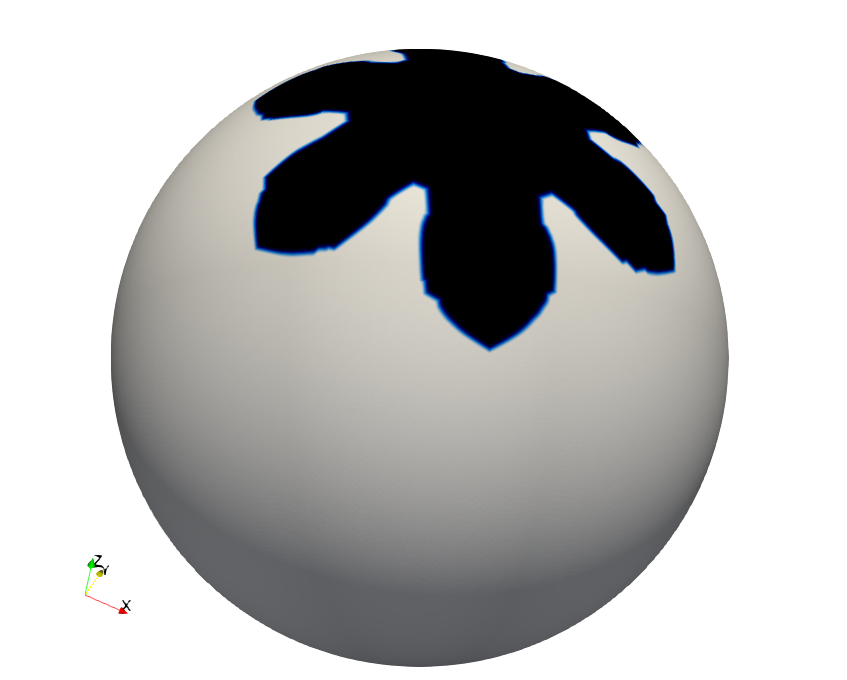}
\includegraphics[angle=-0,width=0.22\textwidth]{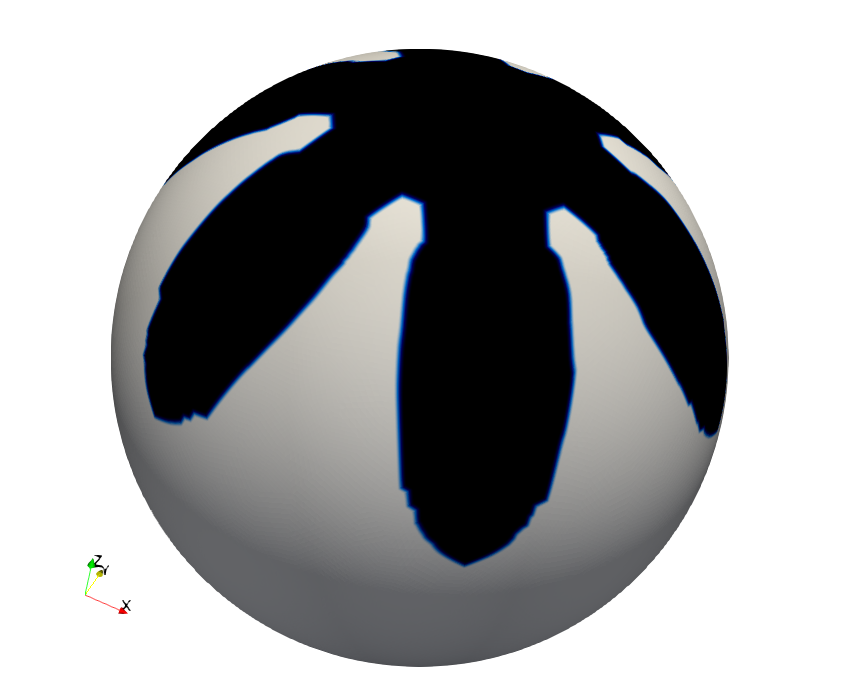} 
\includegraphics[angle=-0,width=0.22\textwidth]{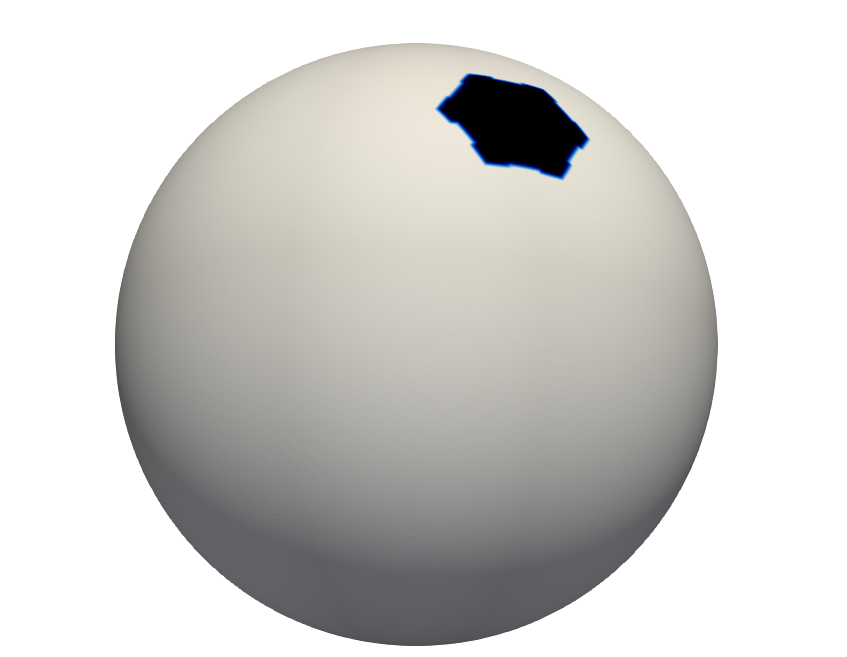}
\includegraphics[angle=-0,width=0.22\textwidth]{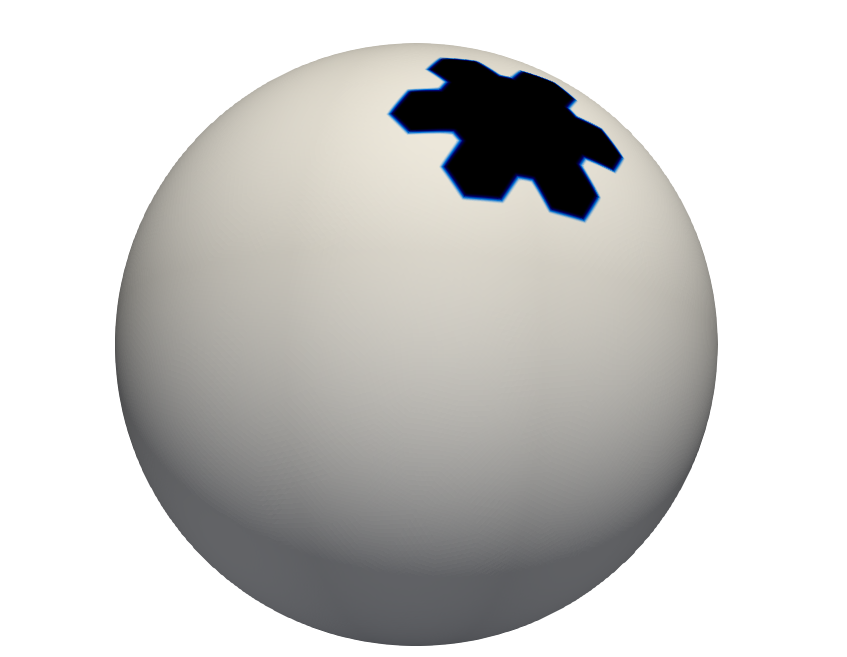}
\includegraphics[angle=-0,width=0.22\textwidth]{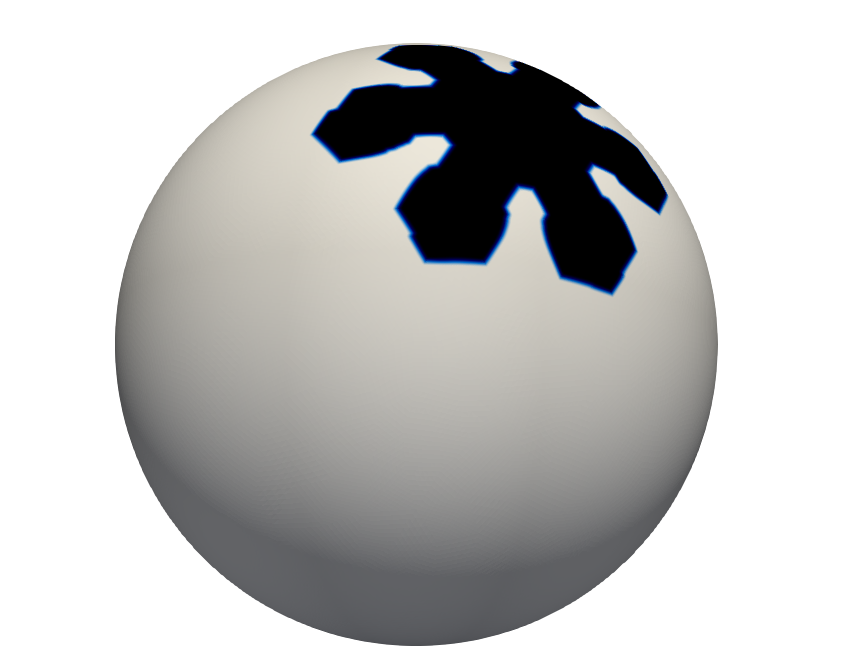}
\includegraphics[angle=-0,width=0.22\textwidth]{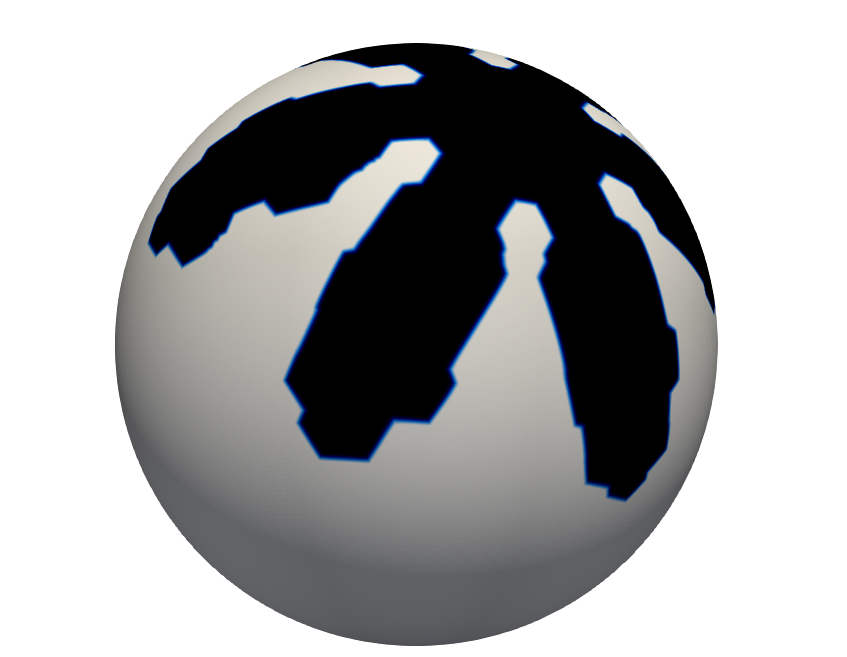}
\caption{($\epsilon=(64\pi)^{-1}$) 
Parameters as in Figure~\ref{fig:spherecap-10more8r0}, but with the
anisotropy \eqref{eq:gammasphere} with \eqref{eq:L3}.
}
\label{fig:spherecap-10L103more8r0}
\end{figure}

\begin{figure}
\center
\includegraphics[angle=-0,width=0.32\textwidth]{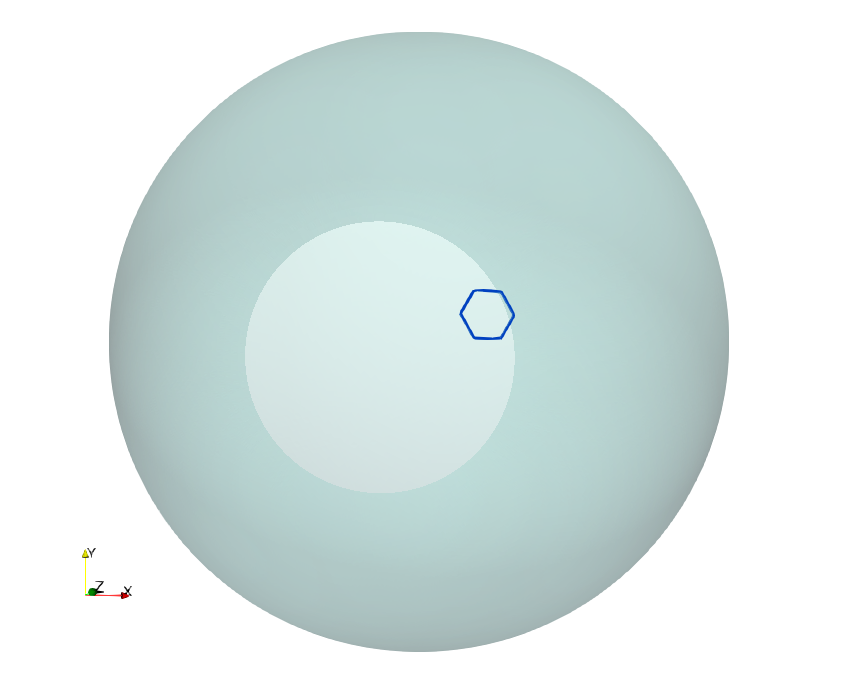}
\includegraphics[angle=-0,width=0.32\textwidth]{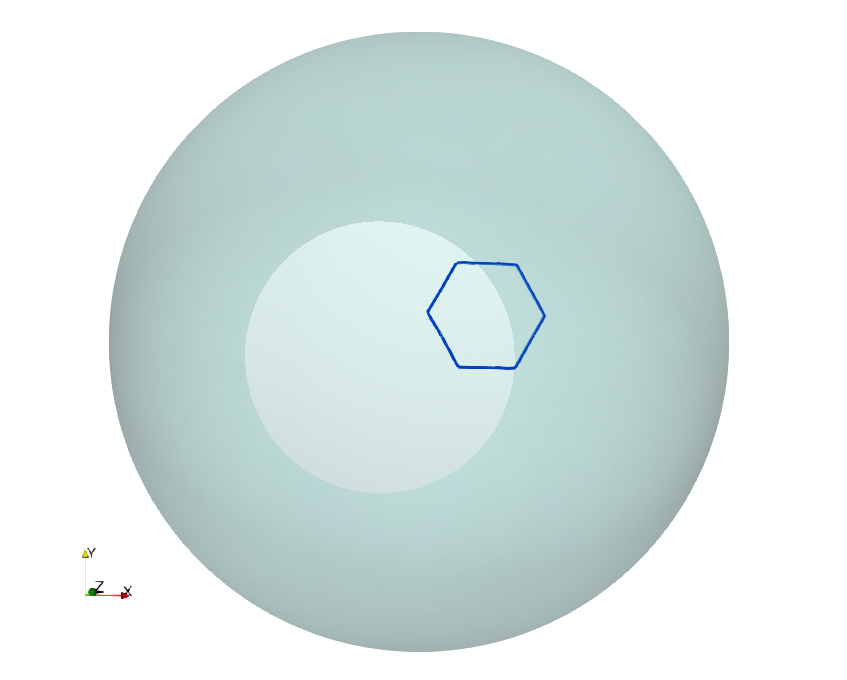}
\includegraphics[angle=-0,width=0.32\textwidth]{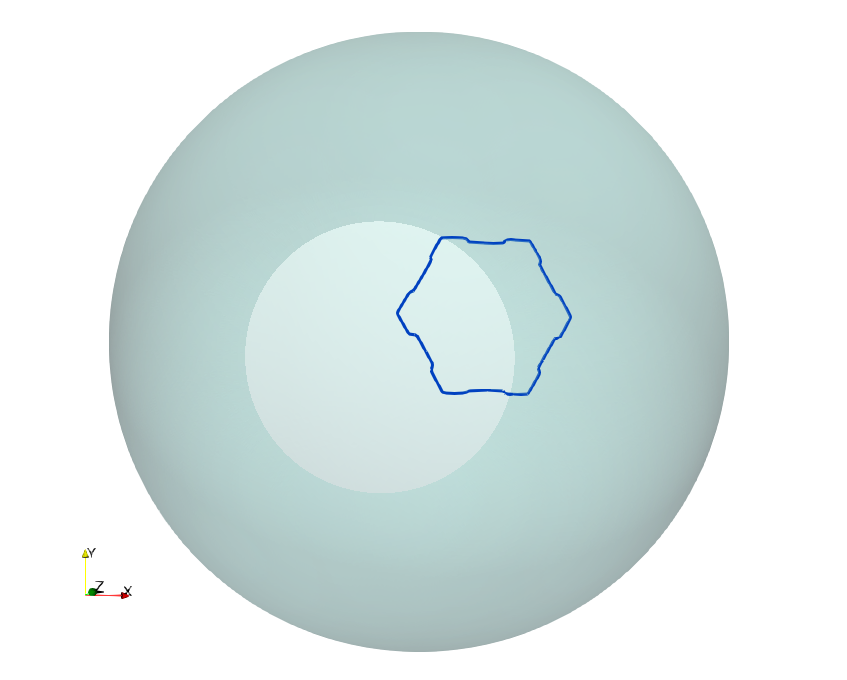}
\includegraphics[angle=-0,width=0.32\textwidth]{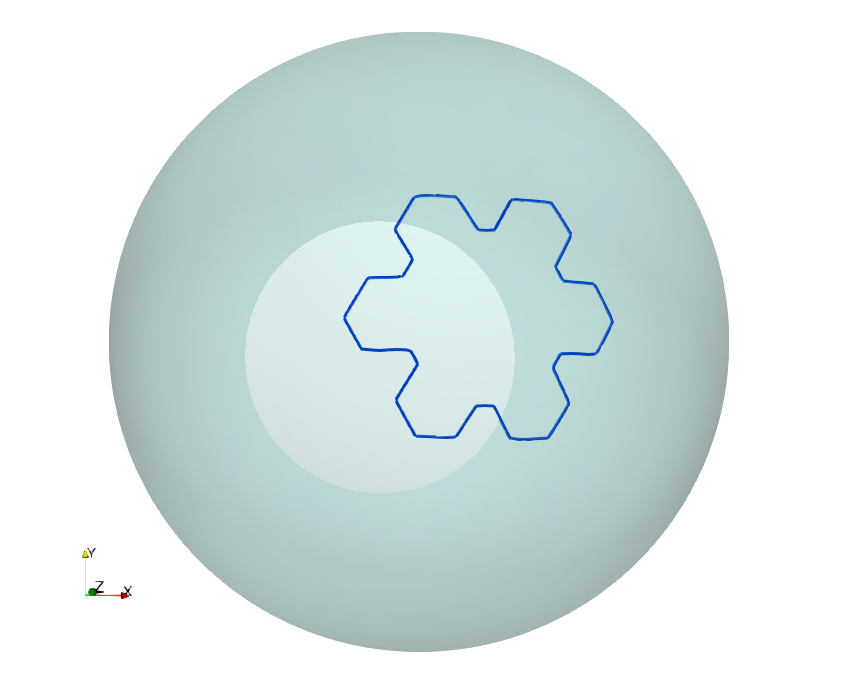}
\includegraphics[angle=-0,width=0.32\textwidth]{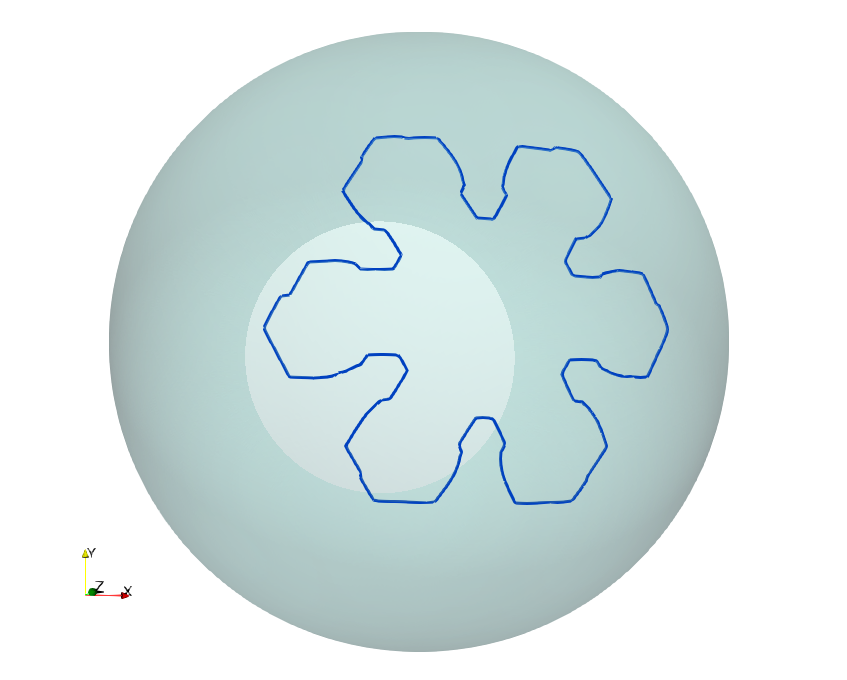}
\includegraphics[angle=-0,width=0.32\textwidth]{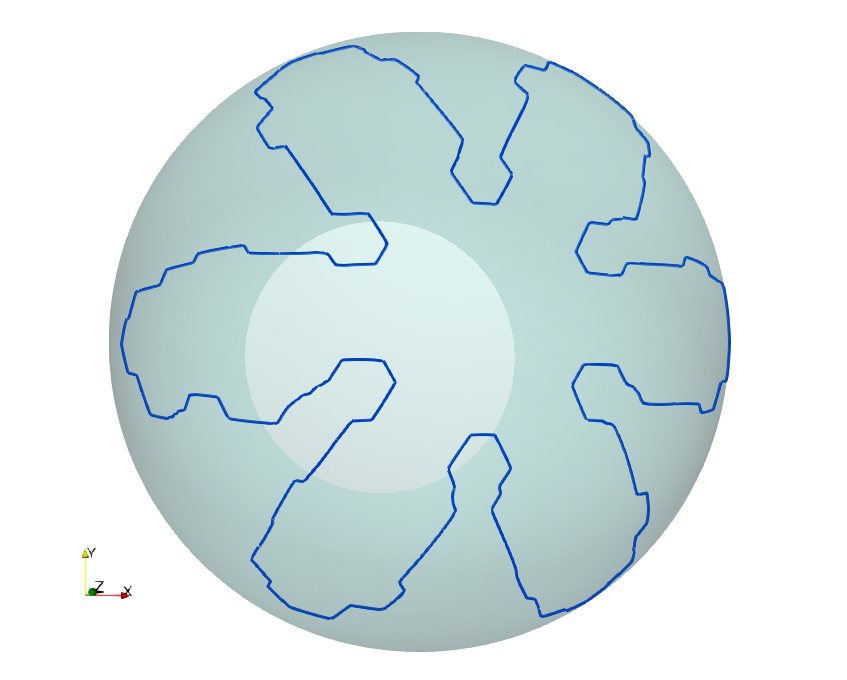}
\caption{($\epsilon=(64\pi)^{-1}$) 
The results from the final row in Figure~\ref{fig:spherecap-10L103more8r0}
displayed from a different point of view, and at times
$t = 0.001, 0.005, 0.01, 0.02, 0.04, 0.1$.
}
\label{fig:spherecap-10L103more8r0_ctop}
\end{figure}

When we repeat the simulations in Figure~\ref{fig:spherecap-10more8r0}
for the new anisotropy these mushy regions disappear, see
Figure~\ref{fig:spherecap-10L103more8r0}. Observe that compared to the
dendritic growth in the former case, the new simulations show the natural
six-fold growth also way beyond the equator. In fact, the interface growths
in Figure~\ref{fig:spherecap-10L103more8r0_ctop} look very similar to the 2d
simulations in \cite{dendritic,jcg,crystal}.


\begin{figure}
\center
\includegraphics[angle=-0,width=0.5\textwidth]{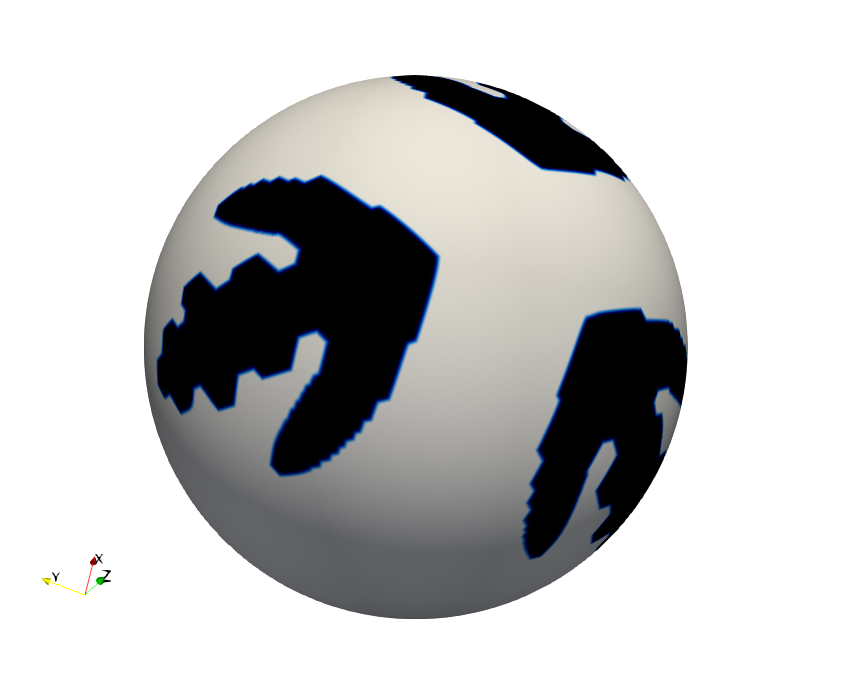}
\caption{An experiment with three initial seeds. Each crystal grows
symmetrically to the others.}
\label{fig:3seeds}
\end{figure}

In Figure~\ref{fig:3seeds} we show a snapshot of a simulation that started
from three initial seeds. Due to the symmetric arrangement of the seeds, the
three crystals continue to grow symmetrically.

Overall the simulations in this subsections underline the ability of the
anisotropies introduced in \S\ref{sec:cons2d} to model a hexagonal crystal 
growth on all parts of the unit sphere segment, in contrast the the results
for the simpler global anisotropies seen in, e.g., \S\ref{sec:cap}.

\def\soft#1{\leavevmode\setbox0=\hbox{h}\dimen7=\ht0\advance \dimen7
  by-1ex\relax\if t#1\relax\rlap{\raise.6\dimen7
  \hbox{\kern.3ex\char'47}}#1\relax\else\if T#1\relax
  \rlap{\raise.5\dimen7\hbox{\kern1.3ex\char'47}}#1\relax \else\if
  d#1\relax\rlap{\raise.5\dimen7\hbox{\kern.9ex \char'47}}#1\relax\else\if
  D#1\relax\rlap{\raise.5\dimen7 \hbox{\kern1.4ex\char'47}}#1\relax\else\if
  l#1\relax \rlap{\raise.5\dimen7\hbox{\kern.4ex\char'47}}#1\relax \else\if
  L#1\relax\rlap{\raise.5\dimen7\hbox{\kern.7ex
  \char'47}}#1\relax\else\message{accent \string\soft \space #1 not
  defined!}#1\relax\fi\fi\fi\fi\fi\fi}

\end{document}